\documentclass[11pt,reqno]{amsart}

\usepackage{amsmath,amssymb,latexsym}
\usepackage{mathtools}
\usepackage{color}
\usepackage{xcolor}
\usepackage{graphicx}
\usepackage{cite}
\usepackage[colorlinks=true]{hyperref}
\hypersetup{urlcolor=blue, citecolor=red, linkcolor=blue}

\newcommand{\I}{\mathrm{i}}

\newcommand{\diff}{\mathrm{d}}

\DeclareMathOperator{\supp}{supp}
\newcommand{\abs}[1]{\lvert#1\rvert}
\newcommand\norm[1]{\left\lVert#1\right\rVert}
\newtheorem{theorem}{Theorem}[section]
\newtheorem{lemma}{Lemma}

\newtheorem{corollary}{Corollary}
\newtheorem*{main-theorem}{Main Theorem}
\newtheorem*{remark*}{Remark}
\newtheorem*{lemma*}{Lemma A.1}

\numberwithin{equation}{section}

\DeclareSymbolFont{largesymbol}{OMX}{yhex}{m}{n}
\DeclareMathAccent{\Widehat}{\mathord}{largesymbol}{"62}

\begin{document}

	\title[]{Global Well-posedness of the Cauchy Problem for the modified Whitham Equations}

	\author{Han Cui \and Yuexun Wang \and Zhouping Xin}
	
	\address{The Institute of Mathematical Sciences, The Chinese University of Hong
Kong, Hong Kong}
	\email{hcui@math.cuhk.edu.hk}
	
	\address{School of Mathematics and Statistics, Lanzhou University, Lanzhou 730000,
P. R. China}
	\email{yuexunwang@lzu.edu.cn}
	
	\address{The Institute of Mathematical Sciences, The Chinese University of Hong
Kong, Hong Kong}
	\email{zpxin@ims.cuhk.edu.hk}

\begin{abstract}
This paper aims to show global existence and modified scattering for the solutions
 of the Cauchy problem to the modified Whitham equations for small, smooth and localized initial data. The main difficulties come from slow decay and non-homogeneity of the Fourier multiplier $(\sqrt{\tanh \xi/\xi})\xi$, which will be overcome by introducing an interaction multiplier theorem and estimating the weighted norms in the frequency space. When estimating the weighted norms, due to loss of derivatives, the energy estimate will be performed in the frequency space, and the absence of time resonance will be effectively utilized by extracting some good terms arising from integration by parts in time before the energy estimate.
 
\end{abstract}

\maketitle

\section{Introduction}
The following shallow water wave model 
\begin{equation}
u_t+uu_x-Lu_x=0
\end{equation}
with the Fourier multiplier operator $L$ defined by \footnote{It usually writes $L=\sqrt{\frac{\tanh D}{D}}$ in convention.}
\begin{equation*}
\begin{aligned}
\widehat{Lf}(\xi)=\sqrt{\frac{\tanh \xi}{\xi}}
\Widehat{f}(\xi)
\end{aligned}
\end{equation*}
was originally proposed by Whitham in \cite{MR0671107} as an alternative to the KdV equation, and is usually referred as the Whitham equation. Note that the dispersion of the Whitham equation reads as $|\xi|^{-1/2}\xi$ 
at high frequencies which is the exact dispersion of the gravity water waves system in infinite depth \cite{MR3314514}, whereas behaves like $\xi-\frac{1}{6}\xi^3$ at low frequencies which is same to the KdV equation (up to a translation). 
Thanks to this different behaviors of the dispersion at low and high frequencies, the Whitham equation displays rich dynamics such as highest waves \cite{MR4002168}, wave breaking \cite{MR4409228, MR3682673, MR4752990}
and long wave limit to the KdV equation \cite{MR3763731}.
The Whitham equation can also be seen as one-dimensional version of the full dispersion KP equation introduced rigorously in \cite{MR3060183} and studied in \cite{MR3177639}.
Moreover, the Whitham equation shares the same dispersive relation with the system for 2D gravity water waves with flat bottoms \cite{MR3914181}.
The Whitham equation is usually compared with another toy model-the fractional KdV (fKdV) equation with the symbol of $L$ taking the form $|\xi|^\alpha$, which not only shares similar dynamics like wave breaking and enhanced lifespan when $\alpha\in (-1,0)$ as well, but also differs on the existence of solitary waves (in fact the fKdV equation does not possess any solitary waves \cite{MR3188389}, whereas the Whitham equation does have \cite{MR2979975}).

The phenomena of both wave breaking and enhanced lifespan reflect
the `competition' between nonlinearity and dispersion, the former scenario usually happens to somehow large solutions while the second scenario occurs for small solutions. When varying the nonlinearity (say $u^pu_x$ with $p=2,3,...$), wave breaking can still happen to both the generalized Whitham equation and the generalized fKdV equation \cite{MR4374000}. Considering the issue of long time existence of solutions to the Cauchy problem, it is usually expected that one can get global existence and modified scattering for $p=2$. In fact, such kind of problems are usually referred to as `long-range' scattering issues, for dispersive models with cubic nonlinearity. This phenomena was first pointed out by Ozawa \cite{MR1121130} in the context of the cubic one-dimensional Schr\"{o}dinger equation. Many results of this type were obtained for various dispersive models, see for instance \cite{MR3519470,MR3462131,MR1613646,MR2850346,MR3382579,MR2199392,MR3121725,MR1491867,MR1687327}.  

Recently, it is shown that the modified fKdV equation admits a unique global solution and a 
scattering with a logarithmic correction in \cite{MR4201835,MR4309883}. The approach is based on the space-time resonance argument and its variants developed in different contexts of PDEs \cite{MR1613646,MR2482120,MR2993751,MR2360438,MR3121725,MR3314514,MR0803256,MR3519470}. 
This work aims to establish results analogous to those in \cite{MR4201835,MR4309883} for the modified Whitham equation which reads\footnote{The argument works for both the nonlinearity of the type $u^2u_x$ (defocusing case) and $-u^2u_x$ (focusing case).}  
\begin{equation}\label{dpb}
u_t-\sqrt{\frac{\tanh D}{D}}u_x=(u^3)_x. 
\end{equation}
The initial data is given by
\begin{equation}\label{initial}
u(0,x)=u_0(x).
\end{equation}
Compared to relevant works \cite{MR3519470,MR3121725,MR4201835,MR4309883}, slow decay and non-homogeneity of the dispersive relation of the Whitham equation appear as new difficulties, which will be overcome by introducing a multiplier theorem allowing interaction estimates and estimating the weighted norms in the frequency space. When estimating weighted norms, due to loss of derivatives, we will perform energy estimate in the frequency space, and the absence of time resonance will be effectively utilized by extracting some good terms arising from integration by parts in time before the energy estimate.
Throughout the paper, we will use $\Lambda(\xi)=(\sqrt{\tanh \xi/\xi}) \xi$ to denote the dispersive relation and $f=e^{-\mathrm{i} t \Lambda(D)}u$ to denote the profile of $u$. We also define the $Z$-norm of $g$ by $\|g\|_Z=\big\|(1+\abs \xi^{30})\Widehat{g}\big\|_{L^\infty}$.

The main result of this paper can be stated as follows:
\begin{theorem}\label{th:main}
Assume  that $N_0=10^7$, $p_0=10^{-5}$ and $p_2=50p_0$\footnote{The parameters $N_0$, $p_0$, $p_2$ could be optimized.}. Suppose that
\begin{equation}\label{smallcondi}
\norm{u_0}_{H^{N_0}}+\norm{xu_0}_{H^1}+\|u_0\|_Z=\epsilon_0 \leq \bar{\epsilon}
\end{equation}
for some constant \(\bar{\epsilon}\) sufficiently small. Then the
Cauchy problem \eqref{dpb}-\eqref{initial} admits a unique global solution $u\in C(\mathbb{R};H^{N_0}(\mathbb{R}))$
satisfying the following uniform bounds
\begin{equation}\label{smallness assumption}
\begin{aligned}
\langle t\rangle^{-p_0}\norm{u}_{H^{N_0}}+\langle t\rangle^{-1/6}\norm{xf}_2+\langle t\rangle^{-p_2} \norm{\partial_x(xf)}_2+\|f\|_Z \lesssim \epsilon_0
\end{aligned}
\end{equation}
for all $t \in [0,\infty)$. Moreover, let
\[
H(\xi,t)=-\frac{6\pi \xi}{\abs{\Lambda''(\xi)}} \int_0^t \frac{\abs{\Widehat f(\xi,s)}^2}{s} \varphi_{>1}(\abs \xi s^{1/3})\, \diff s,
\]
then there exists $w_\infty\in L^\infty(\mathbb{R})$ such that for $t \geq 2$ and $\abs \xi \geq t^{-1/3+\alpha}$ with $\alpha \in (0,10^{-3}]$,
\begin{equation}\label{m-scatter}
\begin{aligned}
\left\lvert (1+\abs \xi^{30})(e^{\I H(\xi,t)}\Widehat f(\xi,t)-w_\infty) \right\rvert \lesssim \epsilon_0 t^{-\kappa},
\end{aligned}
\end{equation}
where $\kappa$ depends only on $\alpha$.
\end{theorem}

We next explain the strategy and difficulties in showing  Theorem \ref{th:main}. 
The local well-posedness of \eqref{dpb}-\eqref{initial} on the time interval \([0,1]\) is standard provided $\|u_0\|_{H^2}$ is sufficiently small, in particular  under the smallness assumption \eqref{smallcondi}. 
Then the existence and uniqueness of the global solutions will be constructed by a bootstrap argument which allows to extend the local solution to a global one.
More precisely, the bootstrap assumptions are
\begin{subequations}
\begin{align}
\norm{f}_{H^{N_0}} &\leq \epsilon_1\langle t\rangle^{p_0}; \label{energyasp}\\
\norm{xf}_{2} &\leq \epsilon_1 \langle t\rangle^{1/6}; \label{weightasp}\\
\norm{\partial_x(xf)}_2 &\leq \epsilon_1 \langle t\rangle^{p_2}; \label{weightasp2}\\
\|f\|_Z &\leq \epsilon_1, \label{Zasp}
\end{align}
\end{subequations}
where $\epsilon_1=100\epsilon_0$. We aim to prove that the assumptions \eqref{energyasp}- \eqref{Zasp} can be improved to
\begin{subequations}
\begin{align}
\norm{f}_{H^{N_0}} &\leq \frac{\epsilon_1}2 \langle t\rangle^{p_0};\label{energyasp-2}\\
\norm{xf}_{2} &\leq \frac{\epsilon_1}2 \langle t\rangle^{1/6};\label{weightasp-2}\\
\norm{\partial_x(xf)}_2 &\leq \frac{\epsilon_1}2 \langle t\rangle^{p_2};\label{weightasp2-2}\\
\|f\|_Z &\leq \frac{\epsilon_1}2.\label{Zasp-2}
\end{align}
\end{subequations}

The difficulties in this paper are different behaviors of high frequencies and low frequencies, slow decay and non-homogeneity of the dispersive relation. The first difficulty appears thoroughly in this paper. As aforementioned the dispersive relation behaves like 2D gravity water waves at high frequencies and like KdV at low frequencies and we will decompose frequencies into high and low frequencies to overcome it. Slow decay will be encountered when energy and weighted norms are estimated. Indeed, since $\Lambda(\xi)$ behaves like KdV in low frequencies, we need the interaction estimate $\norm{u_xu}_\infty \lesssim t^{-1}$ to obtain barely enough decay to close the estimates. The non-homogeneity will be an obstacle when estimating the weighted norms.

Next we sketch the main ideas in the proof. The energy estimate \eqref{energyasp-2} could be proved relatively easily. The $Z$-norm estimate \eqref{Zasp-2} mainly follows from the method of \cite{MR3519470,MR3121725,MR4309883}. The main estimates are divided into high-low, low-low and high-high frequency interactions, in which the low-low frequency interactions are most difficult since they match the scaling. The non-resonant regions will be dealt with by utilizing absence of space resonances and then absence of time resonances. Next, in the resonant region, the non-oscillating terms will be extracted and will lead to modified scattering.

For the weighted norms estimates \eqref{weightasp-2}-\eqref{weightasp2-2}, since the non-homogeneity of $\Lambda(\xi)$ prevents using the standard commuting vector field method, we will estimate the weighted norms in the frequency space as follows. Suppose that the nonlinear term to be estimated has the following generic form:
\[
\int e^{\I s \Phi}m_0 \Widehat f(\xi-\eta-\sigma,s) \Widehat f(\eta,s) \Widehat f(\sigma,s) \, \diff \eta \diff \sigma \diff s,
\]
where $\Phi=\Lambda(\xi-\eta-\sigma)+\Lambda(\eta)+\Lambda(\sigma)-\Lambda(\xi)$ is the resonant function and $m_0$ is the multiplier. We first add a cut-off function $\chi$ on it to make $\xi-\eta-\sigma$ to be a high frequency, which is similar to the process of paralinearization and we may absorb $\chi$ into $m_0$. Then differentiating with respect to $\xi$ will produce the following terms:
\begin{gather*}
I_1=\int e^{\I s \Phi}m_0 \partial_\xi \Widehat f(\xi-\eta-\sigma,s) \Widehat f(\eta,s) \Widehat f(\sigma,s) \, \diff \eta \diff \sigma \diff s;\\
I_2=\int e^{\I s \Phi}\partial_\xi m_0 \Widehat f(\xi-\eta-\sigma,s) \Widehat f(\eta,s) \Widehat f(\sigma,s) \, \diff \eta \diff \sigma \diff s;\\
I_3=\int e^{\I s \Phi} (\I s \Phi_\xi)m_0 \Widehat f(\xi-\eta-\sigma,s) \Widehat f(\eta,s) \Widehat f(\sigma,s) \, \diff \eta \diff \sigma \diff s.
\end{gather*}
For $I_1$, since the equation \eqref{dpb} is quasi-linear, direct estimate on $I_1$ will lead to loss of derivatives. A standard way to resolve it is to use the symmetry of $m_0$. We will utilize the symmetry in the frequency space, which was performed in some literature such as \cite{MR3962880,MR3730012}. However, if we estimate $(\partial_\xi \Widehat f,\partial_\xi \Widehat f)_t$, then the term
\[
\int e^{\I s \Phi} (\I s \Phi_\xi)m \Widehat f(\xi-\eta-\sigma,s) \Widehat f(\eta,s) \Widehat f(\sigma,s) \partial_\xi \Widehat f(-\xi,s) \, \diff \xi \diff \eta \diff \sigma \diff s
\]
will make integration by parts in time unfavorable since if the time derivative hits on $\partial_\xi \Widehat f$, the $I_3$ contained in $\partial_s \partial_\xi \Widehat f$ will make the integration by parts ineffective. Inspired by some works such as \cite{MR2482120,MR3862598}, we integrate by parts in time first and extract good terms (the boundary terms from integration by parts), leading to
\[
\partial_\xi \Widehat f=G+H.
\]
Then energy estimate is performed for $H$ and $I_1$ can be dealt with. $I_2$ is easy to estimate. For $I_3$, we use the space-time method to eliminate $s^2$ from the integral. On the resonant region, we write
\[
\Phi_\xi=m_1 \Phi_\eta+m_2\Phi_\sigma+m_3 \Phi,
\]
which was observed in some works such as \cite{MR3121725,MR2993751,MR3084697}. Hence the space-time resonance method can be performed. Indeed, this process is similar to that of the commuting vector field method, however, in the frequency space, the symmetry requirement is relaxed and the key is that $\Phi_\xi$ vanishes on the resonant set. Another new ingredient is the way to deal with the slow decay rate. Slow decay entails the bilinear interaction estimate $\norm{u_xv}_\infty \lesssim t^{-1}$. However, standard multiplier theorems use the bounds of inputs seperately. Hence we need to establish a multiplier theorem allowing interaction estimates such as
\[
\norm{T_m(f_1,f_2,f_3)}_p \lesssim \norm{f_1}_p \norm{f_2f_3}_\infty.
\]
Inspired by the localized multiplier theorem, we note that the multilinear operators with appropriate multipliers are essentially local, which allows the above kind of multiplier theorem (see Lemma \ref{itmulti}).

The rest of this paper is organized as follows.  Section \ref{Preliminaries} contains some
preliminaries: multiplier theorems, resonance analysis and decay estimates. Section \ref{$H^{N_0}$ Norms} is devoted to closing the energy estimate \eqref{energyasp-2}. We will finish the proof of  $Z$-norm estimate \eqref{Zasp-2} and weighted norms \eqref{weightasp-2}-\eqref{weightasp2-2} in Section \ref{$Z$ Norms} and Section \ref{Weighted Norms}, respectively.

\section{Preliminaries}\label{Preliminaries}
\subsection{Notations}
Let $a \wedge b=\min \{a,b\}$ and $a \vee b=\max \{ a,b \}$. The Fourier transform is defined as
\[
\mathcal{F}(f)(\xi)=\Widehat f(\xi)=\frac {1} {2\pi} \int_{\mathbb{R}} f(x) e^{-\mathrm{i} x \cdot \xi}\, \diff x
\]
and the inverse Fourier transform is defined as
\[
\check f(\xi)= \int_{\mathbb{R}} \hat f(\xi) e^{\mathrm{i} x \cdot \xi}\, \diff \xi.
\]
Set 
\[
\mathcal F (m(D)f)=m(\xi) \widehat f(\xi)
\]
and $\abs \partial^\alpha$ denotes the Fourier multiplier $\abs \xi^\alpha$. The evolution of a profile $f$ by the dispersive semigroup $e^{\mathrm{i} t \Lambda(D)} f$ is denoted by $Ef$.
Let $\Phi=\Lambda(\xi)-\Lambda(\xi-\eta-\sigma)-\Lambda(\eta)-\Lambda(\sigma)$ be the resonant function. Take a bump function $\varphi$ that is supported in $[-3/2,3/2]$ and is equal to $1$ in $[-5/4,5/4]$ and let $$\psi(\xi)=\varphi(\xi)-\varphi(2\xi),\ \varphi_l(\xi)=\varphi(\xi/2^l),\ \psi_l(\xi)=\psi(\xi/2^l).$$
 Moreover, let
\[
\varphi^l_k(\xi)=
\begin{cases}
\psi_k(\xi) &\quad \text{if}\  k>l;\\
\varphi_l(\xi) &\quad \text{if}\  k=l.
\end{cases}
\]
Notations such as $f_k$, $f_{>100}$ denote $\psi_k(D)f$ and $\varphi_{>100}(D)f$. Fixing $t$, we take functions $q_0$, $\dots$, $q_{L+1}$ from $\mathbb{R}$ to $[0,1]$ such that
\begin{itemize}
\item $\abs{L-\log_2(2+t)} \leq 2$;
\item $\sum_{m=0}^{L+1} q_m(s)=1_{[0,t]}(s)$;
\item $\supp q_0 \subset [0,2]$, $\supp q_{L+1} \subset [t-2,t]$, $\supp q_m \subset [2^{m-1},2^{m+1}]$, $q_m \in C^1(\mathbb{R})$ and $\int \abs{q_m'} \lesssim 1$ for $1 \leq m \leq L$.
\end{itemize}

\subsection{Multiplier Theorems}
Set $\norm{m}_S=\norm{\check m}_{L^1}$. Then the following multiplier theorem holds.
\begin{theorem}\label{multi}
Let
\[
\mathcal{F}\big(T_m(f_1,\dots, f_k)\big)(\xi)=\int_{\xi_1+\cdots+\xi_k=\xi} m(\xi_1,\dots,\xi_k) \Widehat{f_1}(\xi_1) \cdots \Widehat{f_k}(\xi_k)\, \diff \xi_1 \cdots \diff \xi_k.
\]
Then
\[
\norm{T_m(f_1,\dots,f_k)}_{p_0} \lesssim \norm{m}_S \norm{f_1}_{p_1} \cdots \norm{f_k}_{p_k},
\]
where
\[
\frac {1} {p_0}=\frac {1} {p_1}+\cdots+\frac {1} {p_k}
\]
and $1 \leq p_1,\dots, p_k \leq \infty$.
\end{theorem}

Sometimes, we use the dual form:
\begin{equation*}
\begin{aligned}
&\bigg|\int m(\xi_1,\ldots,\xi_k) \Widehat{f_1}(\xi_1) \cdots \Widehat{f_k}(\xi_k) \Widehat f\big(-(\xi_1+\cdots+\xi_k)\big)\, \diff \xi_1 \cdots \diff \xi_k\bigg|\\
&\lesssim \norm{f}_{p_0} \norm{f_1}_{p_1} \cdots \norm{f_k}_{p_k},
\end{aligned}
\end{equation*}
where
\[
\frac {1}{p_0}+\frac {1}{p_1}+\cdots +\frac {1}{p_k}=1
\]
and $1 \leq p_0, p_1, \dots, p_k \leq \infty$.

When $m$ is localized, the following lemma is useful.

\begin{lemma}\label{multicri}
If $\supp m \subset B(0,2^{n_1}) \times \cdots \times B(0,2^{n_k})$ and $\abs{\partial^l_{\xi_i} m} \leq A2^{-l\lambda_i}$ for $0 \leq l \leq k+1$ and $1 \leq i \leq k$, then
\[
\norm{m}_S \lesssim A \prod_i 2^{n_i-\lambda_i}.
\]
\end{lemma}

The proofs of Theorem \ref{multi} and Lemma \ref{multicri} are given in \cite[Lemma 5.2]{MR3121725} and \cite[Lemma 2.4]{MR3730012}.
Set
\[
\norm{m}_{S_{k_1,k_2,k_3}}=\norm{m(\xi_1,\xi_2,\xi_3) \varphi_{k_1}(\xi_1) \varphi_{k_2}(\xi_2) \varphi_{k_3}(\xi_3)}_S.
\]
Then by Lemma \ref{multicri}, it holds that
\[
\norm{m}_{S_{k_1,k_2,k_3}} \lesssim A \prod_i 2^{(n_i-\lambda_i)_+}.
\]

Let
\[
I_\lambda(u_1,u_2)=\sup_{\abs{y_1-y_2} \leq \lambda} \abs{u_1(y_1)u_2(y_2)}
\]
be the interaction function of $u_1$ and $u_2$. Then the following multiplier theorem holds.

\begin{lemma}\label{itmulti}
If $\supp m \subset B(0,2^{k_1}) \times B(0,2^{k_2}) \times B(0,2^{k_3})$ and $\abs{\partial^l_{\xi_i} m} \leq A2^{-l\lambda_i}$ for $0 \leq l \leq N+4$ and $1 \leq i \leq 3$, then
\begin{equation*}
\begin{aligned}
\norm{T_m(f_1,f_2,f_3)}_p 
&\lesssim A\bigg(\prod_i 2^{k_i-\lambda_i}\bigg)( \norm{f_1}_p I_\lambda(f_2,f_3)\\
&\quad+\langle 2^{\lambda_{23}}\lambda \rangle^{-N} \norm{f_1}_{p_1} \norm{f_2}_{p_2} \norm{f_3}_{p_3}),
\end{aligned}
\end{equation*}
where $\lambda_{23}=\lambda_2 \wedge \lambda_3$ and $p^{-1}={p_1}^{-1}+{p_2}^{-1}+{p_3}^{-1}$.
\end{lemma}
Indeed, since some similar estimates will be used later, so we will prove Lemma \ref{itmulti} as an illustration. Note that the second term of the right hand side is easy to be estimated, so one could regard it as an error.
\begin{proof}
By the proofs of \cite[Lemma 5.2]{MR3121725} and \cite[Lemma 2.4]{MR3730012},
\[
T_m(f_1,f_2,f_3)=\int K(y_1,y_2,y_3) f_1(x-y_1) f_2(x-y_2) f_3(x-y_3)\, \diff y_1 \diff y_2 \diff y_3,
\]
where
\[
\abs{K(y_1,y_2,y_3)} \lesssim A \langle 2^{\lambda_1}\abs{y_1}+2^{\lambda_2} \abs{y_2}+2^{\lambda_3}\abs{y_3}\rangle^{-N-4} 2^{k_1} 2^{k_2} 2^{k_3}.
\]
Divide $K$ into $K_1=K 1_{\abs{y_2-y_3} \leq \lambda}$ and $K_2=K 1_{\abs{y_2-y_3}>\lambda}$. Then
\begin{align*}
\norm{T_{K_1}(f_1,f_2,f_3)}_p &\lesssim \norm{K_1}_1 \norm{f_1}_p I_\lambda(f_2,f_3)\\
&\lesssim A\bigg(\prod_i 2^{k_i-\lambda_i}\bigg) \norm{f_1}_p I_\lambda(f_2,f_3)
\end{align*}
and
\begin{equation*}
\begin{aligned}
\norm{T_{K_2}(f_1,f_2,f_3)}_p 
&\lesssim \norm{K_2}_1 \norm{f_1}_{p_1} \norm{f_2}_{p_2} \norm{f_3}_{p_3}\\ 
&\lesssim A\bigg(\prod_i 2^{k_i-\lambda_i}\bigg) \langle 2^{\lambda_{23}}\lambda \rangle^{-N} \norm{f_1}_{p_1} \norm{f_2}_{p_2} \norm{f_3}_{p_3}.
\end{aligned}
\end{equation*}
Thus the conclusion follows.
\end{proof}

\subsection{Resonance Analysis}
First, we present some asymptotic properties of $\Lambda(\xi)$. 
\begin{lemma}\label{resoasym}
\begin{enumerate}
\item When $\xi \to 0$, we have
\[
\Lambda(\xi) \simeq \xi, \quad 1-\Lambda'(\xi) \simeq \xi^2, \quad \Lambda''(\xi) \simeq -\xi.
\]
\item When $\xi \to +\infty$, we have
\[
\Lambda^{(k)}(\xi) \simeq \xi^{1/2-k}.
\]
\end{enumerate}
\end{lemma}

Next, we list some resonance inequalities involving $\Lambda(\xi)$.
\begin{lemma}\label{resoin} When $0 \leq a \leq b$, we have
\begin{equation}\label{resoin1}
\begin{aligned}
\Lambda(a)+\Lambda(b)-\Lambda(a+b) \simeq a^{1/2} (a \wedge 1)^{1/2} (b \wedge 1)^2. 
\end{aligned}
\end{equation}
When $0 \leq c \leq a \leq b$, we have
\begin{equation}\label{resoin2}
\begin{aligned}
\Lambda(a)+\Lambda(b)+\Lambda(c)-\Lambda(a+b+c) \gtrsim a^{1/2} (a \wedge 1)^{1/2} (b \wedge 1)^2. 
\end{aligned}
\end{equation}

\end{lemma}
The proof of \eqref{resoin1} is given in \cite[Lemma A.4]{MR3914181}. \eqref{resoin2} follows from convexity of $\Lambda(\xi)$ and \eqref{resoin1}.

The following Corollary is a direct consequence of Lemma \ref{resoin}. 
\begin{corollary}\label{reso32}
If $2^k \simeq \abs{\xi_1} \simeq \abs{\xi_2} \simeq \abs{\xi_3} \gg \abs{\xi_4}$ and $\xi_1+\xi_2+\xi_3+\xi_4=0$, then
\begin{equation}\label{resoin-add}
\abs{\Lambda(\xi_1)+\Lambda(\xi_2)+\Lambda(\xi_3)+\Lambda(\xi_4)} \gtrsim
\begin{cases}
2^{k/2} & k \geq 100,\\
2^{3k} & k<100.
\end{cases}
\end{equation}
\end{corollary}
\begin{proof}
If the signs of $(\xi_1,\xi_2,\xi_3,\xi_4)$ are $(+,+,+,-)$, then the conclusion follows from \eqref{resoin2}. Consider now that the signs of $(\xi_1,\xi_2,\xi_3,\xi_4)$ are $(+,+,-,-)$, and note that
\begin{equation*}
\begin{aligned}
\Lambda(a)+\Lambda(b)-\Lambda(c)-\Lambda(d)
&=\Lambda(a)+\Lambda(b)-\Lambda(a+b)\\
&\quad-(\Lambda(c)+\Lambda(d)-\Lambda(c+d))
\end{aligned}
\end{equation*}
if $a+b=c+d$ and $a$, $b$, $c$, $d \geq 0$. Then this, combined with \eqref{resoin1}, yields the desired estimate.
\end{proof}

\subsection{Decay Estimates}
\begin{lemma}\label{decaylemma}
Let $t \ge 1$. Then it holds that
\begin{subequations}
\begin{gather}
\abs{\abs{\partial}^\beta Ef_\gamma(x,t)} \lesssim t^{-1/3-\beta/3} \langle (x+t)/t^{1/3} \rangle^{-1/4+\beta/2} (\norm{f}_{Z}+t^{-1/6}\norm{xf}_2) \label{decayl}
\intertext{for $\beta \in [0,1]$, $\gamma \in \{ <100,k,\lesssim k \}$ where $k<100$;}
\abs{\abs{\partial}^\beta Ef_{>-100}(x,t)} \lesssim t^{-1/2} (\norm{f}_{Z}+t^{-1/6}\norm{\partial_x(xf)}_2)+t^{-10} \norm{f}_{H^{N_0}}  \label{decayh}
\intertext{for $\beta \in [0,20]$;}
\abs{\abs{\partial}^\beta Ef_{\geq 100}(x,t)} \lesssim
t^{-99/100} (\lVert f \rVert_Z+t^{-p_2}\norm{\partial_x(xf)}_2)
+t^{-10} \norm{f}_{H^{N_0}} 
\label{decayh2}
\end{gather}
for $\beta \in [0,20]$, if $\abs{x/t+1} \ll 1$.
\end{subequations}
\end{lemma}
\begin{proof}
Set 
\begin{equation*}
\begin{aligned}
\Phi(\xi)=\Phi(\xi;x,t):=t^{-1}x\xi+\sqrt{\frac{\tanh \xi}{\xi}}\xi.
\end{aligned}
\end{equation*}	
Then
\begin{equation*}
e^{\I t \Lambda(D)} |D|^\beta g(t,x)
=\int_{-\infty}^\infty e^{\mathrm{i}t\Phi(\xi)}|\xi|^\beta\widehat{g}(t,\xi)\,\diff \xi.
\end{equation*}

\noindent \textbf{Proof of \eqref{decayl}.}
Let
\[
I_\gamma^\beta=\int_{-\infty}^\infty e^{\mathrm{i}t\Phi(\xi)}|\xi|^\beta\widehat{g}(t,\xi) \varphi_\gamma(\xi) \,\diff \xi.
\]
It suffices to prove the case $I^\beta_{<100}$ here since other cases are similar. We will show
	\begin{equation}\label{9}
	\begin{aligned}
	\abs{I_{<100}^\beta} \lesssim t^{-1/3-\beta/3}(1+|x/t^{1/3}|)^{-1/4+\beta/2},
	\end{aligned}
	\end{equation}
for any \(t\geq1\) and any function \(g\) satisfying
		\begin{align}\label{9.5}
		\|g\|_{Z}
		+t^{-\frac{1}{6}}\|xg\|_{L^2}\leq 1.
		\end{align}	
		
	Since the proof is close to \cite{MR3519470}, we focus mainly on the differences between the inhomogeneous Whitham symbol and homogeneous KdV symbol and the new terms,  for the sake of completeness. 	
	To estimate \(I_{<100}^\beta\), we shall study the stationary point of \(\Phi(\xi)\). It is easy to see that on the support of \(\varphi_{100}(\xi)\) \(\partial_\xi\Phi(\xi)=0\) has no root or two roots with opposite signs (corresponding to \(x>-t\)). We only consider the latter case since the first case is much easier to handle by similar calculations. 
	By Lemma \ref{resoasym}, one has
	\begin{equation*}
	\begin{aligned}
	\xi_0 \simeq \sqrt{1+t^{-1}x},
	\end{aligned}
	\end{equation*}	
where \(\xi_0\) denotes the positive root of \(\partial_\xi\Phi(\xi_0)=0\).
	To verify \eqref{9}, by taking complex conjugates, one needs only to show 
	\begin{align}\label{11}
	\left|\int_0^\infty e^{\mathrm{i}t\Phi(\xi)}|\xi|^\beta\widehat{g}(t,\xi)\varphi_{100}(\xi)\,\diff \xi\right|\lesssim t^{-\frac{1}{3}-\frac{\beta}{3}}\max(t^{\frac{1}{3}}\xi_0,1)^{-\frac{1}{2}+\beta}.
	\end{align}

	There are two cases to be considered depending on the size of \(\xi_0\). 
	
	\noindent \textbf{Case 1: $\xi_0\leq t^{-1/3}$.} It suffices to bound the left hand side of \eqref{11} by \(t^{-\frac{1}{3}-\frac{\beta}{3}}\). To this end, we split the integral in \eqref{11} by small and large frequencies:
	\begin{equation*}
	\begin{aligned}
	&\int_0^\infty e^{\mathrm{i}t\Phi(\xi)}|\xi|^\beta\widehat{g}(\xi)\varphi_{100}(\xi)\,\diff \xi\\
	&=\underbrace{\int_0^\infty e^{\mathrm{i}t\Phi(\xi)}|\xi|^\beta\widehat{g}(\xi)\varphi_{10}(t^{1/3}\xi)\varphi_{100}(\xi)\,\diff \xi}_{A_1}\\
	&\quad+\underbrace{\int_0^\infty e^{\mathrm{i}t\Phi(\xi)}|\xi|^\beta\widehat{g}(\xi)\big(1-\varphi_{10}(t^{1/3}\xi)\big)\varphi_{100}(\xi)\,\diff \xi}_{A_2}.
	\end{aligned}
	\end{equation*}
	Using \eqref{9.5}, one immediately gets 
	\begin{equation*}
	\begin{aligned}
	|A_1|\lesssim \|g\|_{Z}\int_0^\infty |\xi|^\beta\varphi_{10}(t^{1/3}\xi)\,\diff \xi
	\lesssim t^{-\frac{1}{3}-\frac{\beta}{3}}.
	\end{aligned}
	\end{equation*}
	To handle the second term \(A_2\),  one integrates by parts to deduce
	\begin{equation*}
	\begin{aligned}
	|A_2|&\lesssim 
	\underbrace{t^{-1}\int_0^\infty \left|\partial_\xi\big[|\xi|^\beta(\partial_\xi\Phi)^{-1}\big(1-\varphi_{10}(t^{1/3}\xi)\big)\varphi_{100}(\xi)\big]\widehat{g}(\xi)\right|\,\diff \xi}_{A_{21}}\\
	&\quad+\underbrace{t^{-1}\int_0^\infty \left||\xi|^\beta(\partial_\xi\Phi)^{-1}\big(1-\varphi_{10}(t^{1/3}\xi)\big)\varphi_{100}(\xi)\partial_\xi\widehat{g}(\xi)\right|\,\diff \xi}_{A_{22}}.
	\end{aligned}
	\end{equation*}
	The key in estimating \(A_{21}\) and \(A_{22}\) is to bound \(|\partial_\xi\Phi|\) from below. 
	Indeed, it follows from Lemma \ref{resoasym} that
	\begin{equation*}
	|\partial_\xi\Phi|\gtrsim \xi^2,
	\end{equation*}
where one has used \(\xi\gg \xi_0\) on the support of the integral $A_2$. Since \(\varphi_{100}(\xi)\) is bounded, 
\(A_{22}\) and \(A_{21}\) can be estimated exactly as that in \cite{MR3519470}, except the term involving \(\partial_\xi \varphi_{100}(\xi)\), that is,
\begin{equation*}
\begin{aligned}
&t^{-1}\int_0^\infty \left||\xi|^\beta(\partial_\xi\Phi)^{-1}\big(1-\varphi_{10}(t^{1/3}\xi)\big)\partial_\xi\varphi_{100}(\xi)\widehat{g}(\xi)\right|\,\diff \xi\\
&\lesssim t^{-1}\|g\|_{Z}\int_0^\infty |\xi|^{\beta-3}|1-\varphi_{10}(t^{1/3}\xi)|\,\diff \xi\\
&\lesssim t^{-\frac{1}{3}-\frac{\beta}{3}}.
\end{aligned}
\end{equation*}
Hence one has
	\begin{equation*}
	\begin{aligned}
	\abs{A_2}\lesssim t^{-\frac{1}{3}-\frac{\beta}{3}}.
	\end{aligned}
	\end{equation*}

\noindent{\textbf{Case 2: \(\xi_0\geq t^{-1/3}\)}.} We will show that the left hand side of \eqref{11} is bounded by \(t^{-\frac{1}{2}}\xi_0^{-\frac{1}{2}+\beta}\). 
	Since the resonant contributions concentrate on \(\xi\approx \xi_0\), one can split the integral in \eqref{11} as follows:
	\begin{equation*}
	\begin{aligned}
	&\int_0^\infty e^{\mathrm{i}t\Phi(\xi)}|\xi|^\beta\widehat{g}(\xi)\varphi_{100}(\xi)\,\diff \xi\\
	&=\underbrace{\int_0^\infty e^{\mathrm{i}t\Phi(\xi)}|\xi|^\beta\widehat{g}(\xi)\big(1-\psi(\xi/\xi_0)\big)\varphi_{100}(\xi)\,\diff \xi}_{A_3}\\
	&\quad+\underbrace{\int_0^\infty e^{\mathrm{i}t\Phi(\xi)}|\xi|^\beta\widehat{g}(\xi)\psi(\xi/\xi_0)\varphi_{100}(\xi)\,\diff \xi}_{A_4}.
	\end{aligned}
	\end{equation*}
	To tackle \(A_3\), one integrates by parts to bound 
	\begin{equation*}
	\begin{aligned}
	|A_3|&\lesssim \underbrace{t^{-1}\int_0^\infty \left|\partial_\xi\big[|\xi|^\beta(\partial_\xi\Phi)^{-1}\big(1-\psi(\xi/\xi_0)\big)\varphi_{100}(\xi)\big]\widehat{g}(\xi)\right|\,\diff \xi}_{A_{31}}\\
	&\quad+\underbrace{t^{-1}\int_0^\infty \left||\xi|^\beta(\partial_\xi\Phi)^{-1}\big(1-\psi(\xi/\xi_0)\big)\varphi_{100}(\xi)\partial_\xi\widehat{g}(\xi)\right|\,\diff \xi}_{A_{32}}.
	\end{aligned}
	\end{equation*}
	It follows from Lemma \ref{resoasym} that
	\begin{equation*}
	|\partial_\xi\Phi| \gtrsim \max(\xi,\xi_0)^2,
	\end{equation*}
on the support of the integral $A_3$. Again, \(A_{32}\) can be handled exactly as in \cite{MR3519470}, 
and the new term in \(A_{31}\) is the term involving \(\partial_\xi \varphi_{100}(\xi)\), that is,
	\begin{equation*}
	\begin{aligned}
	&t^{-1}\int_0^\infty \big|\max(\xi,\xi_0)^{-2}|\xi|^{\beta}\big(1-\psi(\xi/\xi_0)\big)\partial_\xi\varphi_{100}(\xi)\widehat{g}(\xi)\big|\,\diff \xi\\
	&\lesssim t^{-1}\xi_0^{\beta-2},
	\end{aligned}
	\end{equation*}
which is stronger than the desired bound
	\(t^{-\frac{1}{2}}\xi_0^{-\frac{1}{2}+\beta}\) due to \(\xi_0\geq t^{-1/3}\). One concludes that
	\begin{equation*}
	\begin{aligned}
	\abs{A_3}\lesssim t^{-\frac{1}{3}-\frac{\beta}{3}}.
	\end{aligned}
	\end{equation*}

We next study the stationary contributions term \(A_4\).
  Let \(l_0\) be the smallest integer with the property that \(2^{l_0}\geq (t\xi_0)^{-1/2}\). Then it follows that
	\begin{equation*}
	\begin{aligned}
	A_4&= \underbrace{\int_{-\infty}^\infty e^{\mathrm{i}t\Phi(\xi)}|\xi|^\beta\widehat{g}(t,\xi)\psi(\xi/\xi_0)\varphi_{l_0}(\xi-\xi_0)\varphi_{100}(\xi)\,\diff \xi}_{A_{4l_0}}\\
	&\quad+\sum_{l\geq l_0+1}\underbrace{\int_{-\infty}^\infty e^{\mathrm{i}t\Phi(\xi)}|\xi|^\beta\widehat{g}(t,\xi)\psi(\xi/\xi_0)\psi_l(\xi-\xi_0)\varphi_{100}(\xi)\,\diff \xi}_{A_{4l}}.
	\end{aligned}
	\end{equation*}
	The desired bound for \(A_{4l_0}\) is immediate from the definition of \(l_0\):
	\begin{equation*}
	\begin{aligned}
	|A_{4l_0}|&\lesssim \|g\|_{Z}\int_{-\infty}^\infty \left| |\xi|^\beta\psi(\xi/\xi_0)\varphi_{l_0}(\xi-\xi_0) \right| \,\diff \xi\\
	&\lesssim |\xi_0|^{\beta}2^{l_0}\lesssim t^{-\frac{1}{2}}\xi_0^{-\frac{1}{2}+\beta}.
	\end{aligned}
	\end{equation*}
	It remains to handle \(A_{4l}\) for \(l\geq l_0+1\). Integration by parts yields
	\begin{equation*}
	\begin{aligned}
	|A_{4l}|&\lesssim 
	\underbrace{t^{-1}\int_{-\infty}^\infty \left|\partial_\xi\big[|\xi|^\beta(\partial_\xi\Phi)^{-1}\psi(\xi/\xi_0)\psi_l(\xi-\xi_0)\varphi_{100}(\xi)\big]\widehat{g}(\xi)\right|\,\diff \xi}_{A_{4l,1}}\\
	&\quad+\underbrace{t^{-1}\int_{-\infty}^\infty \left||\xi|^\beta(\partial_\xi\Phi)^{-1}\psi(\xi/\xi_0)\psi_l(\xi-\xi_0)\varphi_{100}(\xi)\partial_\xi\widehat{g}(\xi)\right|\,\diff \xi}_{A_{4l,2}}.
	\end{aligned}
	\end{equation*}
	Note first that 
	\begin{equation*}
	|\partial_\xi\Phi| \simeq   |\Lambda''(\xi_0)(\xi_0-\xi)| \gtrsim 2^l\xi_0
	\end{equation*}
on the support of the integral $A_{4l}$. Then, comparing to \cite{MR3519470}, one needs only to focus on the following new term
\begin{equation*}
\begin{aligned}
&t^{-1}\int_{-\infty}^\infty \left||\xi|^\beta(\partial_\xi\Phi)^{-1}\psi(\xi/\xi_0)\psi_l(\xi-\xi_0)
\partial_\xi\varphi_{100}(\xi)\widehat{g}(\xi)\right|\,\diff \xi\\
&\lesssim t^{-1}2^{-l}\xi_0^{-1}\|g\|_{Z}\int_{-\infty}^\infty |\xi|^{\beta-1}\psi(\xi/\xi_0)\psi_l(\xi-\xi_0)\,\diff \xi\\
&\lesssim t^{-1}2^{-l}\xi_0^{-1+\beta},
\end{aligned}
\end{equation*}
which leads to a bound \(t^{-\frac{1}{2}}\xi_0^{-\frac{1}{2}+\beta}\) by taking summation over \(l\geq l_0+1\) and using \(2^{l_0}\geq (t\xi_0)^{-1/2}\). Collecting all the estimates leads to
\begin{equation*}
\begin{aligned}
\abs{A_4} \lesssim t^{-\frac{1}{2}}\xi_0^{-\frac{1}{2}+\beta}.
\end{aligned}
\end{equation*}

	\noindent\textbf{Proof of \eqref{decayh}.}
	We follow the idea of \cite{SSWZ}.
	Let
\[
I_{>-100}^\beta=\int_{-\infty}^\infty e^{\mathrm{i}t\Phi(\xi)}|\xi|^\beta\widehat{g}(t,\xi) \varphi_{>-100}(\xi) \,\diff \xi.
\]
	We claim that
	\begin{equation}\label{12}
	\begin{aligned}
	|I_{>-100}^\beta|\lesssim t^{-1/2},
	\end{aligned}
	\end{equation}
for any \(t\geq1\) and any function \(g\) satisfying
	\begin{align}\label{13}
	\|g\|_{Z}
	+t^{-\frac{1}{6}}\|\partial_x(xg)\|_{L^2}+t^{1/2}t^{-10}\norm{g}_{H^{N_0}} \leq 1.
	\end{align}
Decompose
	\begin{equation*}
	\begin{aligned}
	I_{>-100}^\beta:=\sum_{k\in \mathbb{Z}}\underbrace{\int_{-\infty}^\infty e^{\mathrm{i}t\Phi(\xi)}|\xi|^\beta\widehat{P_kg}(\xi)\varphi_{>-100}(\xi)\,\diff \xi}_{I_{>-100, k}^\beta}\approx \sum_{k\in \mathbb{N}}I_{>-100, k}^\beta,
	\end{aligned}
	\end{equation*}
where one has used the property of the support of the integral $I^\beta_{>-100}$. For \eqref{12}, it suffices to show 	  
	\begin{equation}\label{14}
	\begin{aligned}
	|I_{>-100, k}^\beta|
	\lesssim t^{-\frac{1}{2}}2^{\frac{3k}{4}}2^{\beta k}\|\widehat{g}\|_{L^\infty}
	+t^{-\frac{3}{4}}2^{\frac{k}{8}}2^{\beta k}(\|\widehat{g}\|_{L^2}+2^k\|\partial\widehat{g}\|_{L^2}).
	\end{aligned}
	\end{equation}
We assume \eqref{14} for a moment and use it to show \eqref{12}. If \(2^k\leq t^{1/100}\), then one can use \eqref{14} to deduce that 
\begin{equation}\label{16}
\begin{aligned}
2^k|I_{>-100, k}^\beta|
&\lesssim t^{-\frac{1}{2}}2^{\frac{3k}{4}}2^{\beta k}2^k\|\widehat{P_k^\prime g}(\xi)\|_{L^\infty}\\
&\quad+t^{-\frac{3}{4}}2^{\frac{k}{8}}2^{\beta k}2^k(\|\widehat{P_k^\prime g}\|_{L^2}+2^k\|\partial\widehat{P_k^\prime g}\|_{L^2})\\
&\lesssim t^{-\frac{1}{2}}+t^{-\frac{3}{4}}2^{\frac{k}{8}}2^{\beta k}2^kt^{\frac{1}{6}}\lesssim t^{-\frac{1}{2}}.
\end{aligned}
\end{equation}
While for \(2^k\geq t^{1/100}\), one can get directly
\begin{equation}\label{17}
2^k|I_{>-100, k}^\beta| \lesssim 2^k 2^{\beta k} 2^{\frac k2}\norm{P_k g}_2 \lesssim 2^{\beta k} 2^{\frac{3k}2} 2^{-N_0k}t^{-\frac 12}t^{10} \lesssim t^{-\frac 12}.
\end{equation}
Due to the decay factor \(2^k\), \eqref{12} follows from \eqref{16} and \eqref{17}.

We now show \eqref{14}.  	
Suppose that $c^{-1}_0 \abs{\xi}^{-1/2} \le  \Lambda'(\xi) \le c_0 \abs{\xi}^{-1/2}$ for $\abs{\xi} \ge 2^{-150}$.	
	Set
	\begin{align*}
	\mathcal{I}:=\left\{k\in \mathbb{N}:  \frac{c^{-1}_0}{4} |tx^{-1}|\leq 2^{k/2}\leq 4c_0|tx^{-1}|\right\}.
	\end{align*}

	\noindent{\bf{Case 1: \(k\in \mathbb{N}\setminus\mathcal{I}\).}} 
	Integration by parts yields
	\begin{equation*}
	\begin{aligned}
	|I_{>-100, k}^\beta|
	&\lesssim \underbrace{t^{-1}\int_{-\infty}^\infty \left|\partial_\xi\big[|\xi|^\beta(\partial_\xi\Phi)^{-1}\varphi_{>-100}(\xi)\psi_k(\xi)\big]\widehat{g}(t,\xi)\right|\,\diff \xi}_{B_1}\\
	&\quad+\underbrace{t^{-1}\int_{-\infty}^\infty \left||\xi|^\beta(\partial_\xi\Phi)^{-1}\varphi_{>-100}(\xi)\psi_k(\xi)\partial_\xi\widehat{g}(t,\xi)\right|\,\diff \xi}_{B_2}.
	\end{aligned}
	\end{equation*}
	Noting that 
	\begin{align*}
	|\partial_\xi\Phi(\xi)|\gtrsim \big||t^{-1}x|-\Lambda'(\xi)\big|\gtrsim 2^{-k/2},
	\end{align*}
	one then may estimate directly to get
	\begin{equation*}
	\begin{aligned}
	B_1&\lesssim t^{-1}\|\widehat{g}\|_{L^2}\bigg(\int_{-\infty}^\infty \left|\partial_\xi\big[|\xi|^\beta(\partial_\xi\Phi)^{-1}\varphi_{>-100}(\xi)\psi_k(\xi)\big]\right|^2\,\diff \xi\bigg)^{1/2}\\
	&\lesssim t^{-1}2^{\beta k}\|\widehat{g}\|_{L^2},
	\end{aligned}
	\end{equation*}
	and
	\begin{equation*}
	\begin{aligned}
	B_2&\lesssim t^{-1}\|\partial\widehat{g}\|_{L^2}\bigg(\int_{-\infty}^\infty \left||\xi|^\beta(\partial_\xi\Phi)^{-1}\varphi_{>-100}(\xi)\psi_k(\xi)\right|^2\,\diff \xi\bigg)^{1/2}\\
	&\lesssim t^{-1}2^{k}2^{\beta k}\|\partial\widehat{g}\|_{L^2}.
	\end{aligned}
	\end{equation*}
	Thus, it holds that
	\begin{equation*}
	\begin{aligned}
	|I_{>-100, k}^\beta|
	\lesssim t^{-3/4}2^{k/8}2^{\beta k}(2^k\|\partial\widehat{g}\|_{L^2}+\|\widehat{g}\|_{L^2}).
	\end{aligned}
	\end{equation*}

	\noindent{\bf{Case 2: \(k\in \mathcal{I}\).}} It is easy to check that on the support of $\varphi_{>-100}(\xi)$
	\(\partial_\xi\Phi(\xi)=0\) has no root or two roots with opposite signs (corresponding to \(x<0\)). 
	We only consider the latter case since the first case is much easier. Denote by \(\xi_0\) the positive root of \(\partial_\xi\Phi(\xi)=0\).
	 Let \(l_0\) be the smallest integer satisfying  \(2^{l_0}\geq t^{-1/2}2^{3k/4}\).  Then, one has
	\begin{align*}
	|I_{>-100, k}^\beta|
	&\leq \underbrace{\left| \int_{-\infty}^\infty e^{\mathrm{i}t\Phi(\xi)}|\xi|^\beta\widehat{P_kg}(\xi)\big(1-\varphi_{>-100}(\xi)\big)\varphi_{l_0}\big(\xi-\xi_0\big)\,\diff \xi\right|}_{J_{l_0}}\\
	&\quad+\sum_{l\geq l_0+1}\underbrace{\left|\int_{-\infty}^\infty e^{\mathrm{i}t\Phi(\xi)}|\xi|^\beta\widehat{P_kg}(\xi)\varphi_{>-100}(\xi)\psi_l(\xi-\xi_0)\,\diff \xi\right|}_{J_l}.
	\end{align*}
It is easy to see that
	\begin{align*}
	J_{l_0}
	\leq t^{-1/2}2^{3k/4}2^{\beta k}\|\widehat{g}\|_{L^\infty}.
	\end{align*}
It remains to bound \(J_l\) for \(l\geq l_0+1\). Integration by parts yields
\begin{equation*}
\begin{aligned}
J_l
&\lesssim \underbrace{t^{-1}\int_{-\infty}^\infty \left|\partial_\xi\big[|\xi|^\beta(\partial_\xi\Phi)^{-1}\varphi_{>-100}(\xi)\psi_l(\xi-\xi_0)\psi_k(\xi)\big]\widehat{g}(\xi)\right|\,\diff \xi}_{B_3}\\
&\quad+\underbrace{t^{-1}\int_{-\infty}^\infty \left||\xi|^\beta(\partial_\xi\Phi)^{-1}\varphi_{>-100}(\xi)\psi_l(\xi-\xi_0)\psi_k(\xi)\partial_\xi\widehat{g}(\xi)\right|\,\diff \xi}_{B_4}.
\end{aligned}
\end{equation*}
Noting 
\begin{align*}
|\partial_\xi\Phi(\xi)|\gtrsim 2^{l-\frac{3k}{2}},
\end{align*}
one then may estimate that
\begin{equation*}
\begin{aligned}
B_3&\lesssim t^{-1}\|\widehat{g}\|_{L^\infty}\int_{-\infty}^\infty \big|\partial_\xi\big[|\xi|^\beta(\partial_\xi\Phi)^{-1}\varphi_{>-100}(\xi)\\
&\qquad\qquad\qquad\qquad\qquad\qquad\times\psi_l(\xi-\xi_0)\psi_k(\xi)\big]\big|\,\diff \xi\\
&\lesssim t^{-1}2^{-l+\frac{3k}{2}}2^{\beta k}\|\widehat{g}\|_{L^\infty},
\end{aligned}
\end{equation*}
and
\begin{equation*}
\begin{aligned}
B_4&\lesssim t^{-1}\|\partial\widehat{g}\|_{L^2}\bigg(\int_{-\infty}^\infty \left||\xi|^\beta(\partial_\xi\Phi)^{-1}\varphi_{>-100}(\xi)\psi_l(\xi-\xi_0)\psi_k(\xi)\right|^2\,\diff \xi\bigg)^{1/2}\\
&\lesssim t^{-1}2^{-\frac{l}{2}+\frac{3k}{2}}2^{\beta k}\|\partial\widehat{g}\|_{L^2}.
\end{aligned}
\end{equation*}
Hence 
\begin{align*}
\sum_{l\geq l_0+1}J_l\lesssim t^{-1/2}2^{3k/4}2^{\beta k}\|\widehat{g}\|_{L^\infty}
+t^{-3/4}2^{k/8}2^k2^{\beta k}\|\partial\widehat{g}\|_{L^2}.
\end{align*}
Thus we have shown that
\begin{equation*}
\begin{aligned}
|I_{>-100,k}^\beta|
\lesssim t^{-1/2}2^{3k/4}2^{\beta k}\|\widehat{g}\|_{L^\infty}+t^{-3/4}2^{k/8}2^{\beta k}(2^k\|\partial\widehat{g}\|_{L^2}+\|\widehat{g}\|_{L^2}).
\end{aligned}
\end{equation*}

The proof of \eqref{decayh2} is similar to {\bf{Case 1}} of \eqref{decayh} since $\abs{\partial _\xi \Phi} \gtrsim 1$ for this case.
\end{proof}
Let
\[
\norm{f}_X=\norm{f}_Z+t^{-1/6}\norm{xf}_2+t^{-p_2}\norm{\partial_x(xf)}_2+t^{-p_0}\norm{f}_{H^{N_0}}.
\]
Combining \eqref{decayl} and \eqref{decayh}, we have the following corollary.
\begin{corollary}
Assume that $t \ge 1$ and $\norm{f}_{X} \lesssim 1$. Then
\begin{equation}\label{decayl6}
\norm{Ef}_6 \lesssim t^{-1/3+1/18}.
\end{equation}
\end{corollary}
The following decay estimates concerning interactions will be used, whose proof follows from Lemma \ref{decaylemma}.
\begin{lemma}\label{decayest2}
\mbox{}\\
Assume that $t \ge 1$ and $\norm{u}_{X}$, $\norm{v}_{X} \lesssim 1$. Then 
\begin{equation}\label{decayinter}
\begin{aligned}
I_{t^{1/3}}(u_x,v) &\lesssim t^{-1} 
\end{aligned}
\end{equation}
and
\begin{equation}\label{decayhl}
\begin{aligned}
I_{t}(\partial^\alpha u_{\geq 100},v_{\leq -100}) &\lesssim t^{-1}
\end{aligned}
\end{equation}
for $0 \leq \alpha \leq 20$.
\end{lemma}

\subsection{Estimates for Time Derivatives}
The following estimate on the $L^2$-norm of $\partial_t \Widehat{f}$ will be useful in the absence of time resonances.
\begin{lemma} Under the bootstrap assumptions \eqref{energyasp}-\eqref{Zasp}, for $t \ge 1$, it holds that
\begin{equation}\label{psl}
\lVert\partial_t \Widehat{f_k}\rVert_2 \lesssim \epsilon_1^3 2^{k/2}t^{-1}.
\end{equation}
\end{lemma}
\begin{proof} Notice that
\[
\lVert \partial_t \Widehat{f_k}\rVert_2=\norm{P_k[(u^3)_x]}_2.
\]
One may assume $k_1 \geq k_2 \geq k_3$ and estimate
\begin{align*}
\norm{P_k[(u^3)_x]}_2 
&\lesssim 2^k \sum_{k_1,k_2,k_3} \norm{P_k(u_{k_1}u_{k_2}u_{k_3})}_2\\
&\lesssim 2^k \sum_{k_1,k_2,k_3} \norm{u_{k_1}u_{k_2}}_\infty \norm{u_{k_3}}_2\\
&\lesssim \epsilon_1^3 2^k \sum_{k_1,k_2,k_3} 2^{-k_1}t^{-1} 2^{k_3/2}\\
&\lesssim \epsilon_1^3 2^k \sum_{k_1,k_2} 2^{-k_1} 2^{k_2/2} t^{-1}\\
&\lesssim \epsilon_1^3 2^k \sum_{k_1} 2^{-k_1/2} t^{-1}
\lesssim \epsilon_1^3 2^{k/2}t^{-1}.
\end{align*}
Hence \eqref{psl} follows.

\end{proof}

\subsection{Consequences of \eqref{energyasp}-\eqref{Zasp}}

According to the bootstrap assumptions \eqref{energyasp}-\eqref{Zasp}, one has
\begin{subequations}
\begin{align}
\norm{f}_2 &\lesssim \epsilon_1; \label{estl2}\\
\norm{Ef}_\infty &\lesssim \epsilon_1 t^{-1/3}; \label{estli}\\
\norm{Ef_k}_\infty &\lesssim \epsilon_1 t^{-1/2}2^{-15k}\ \text{ if } k \geq 0; \label{estlih}\\
\norm{Ef_k}_\infty &\lesssim \epsilon_1 t^{-1/2} 2^{-k/2}\ \text{ if } k \leq 0; \label{estlil}\\
\norm{Ef_k}_\infty, \norm{Ef_{<k}}_\infty &\lesssim \epsilon_1 2^k; \label{estlilow}\\
\norm{f_k}_2, \norm{f_{<k}}_2 &\lesssim \epsilon_1 2^{k/2}. \label{estl2l}
\end{align}
\end{subequations}
\eqref{estl2} follows from energy conservation. \eqref{estli} follows from Lemma \ref{decaylemma} and \eqref{estlih} and \eqref{estlil} follow from applying \eqref{decayh} and \eqref{decayl} to $\norm{\abs{\partial}^{15}Ef}_\infty$ and $\norm{\abs{\partial}^{1/2}Ef}_\infty$ respectively. \eqref{estlilow} and \eqref{estl2l} could be obtained in the frequency space using the $Z$-norm.

\section{$H^{N_0}$ Norms}\label{$H^{N_0}$ Norms}
To prove \eqref{energyasp-2}, it suffices to show
\begin{equation*}
\begin{aligned}
\norm{u}_{H^{N_0}} \lesssim \epsilon_0+\epsilon_1^2 t^{p_0}
\end{aligned}
\end{equation*}
for $t \geq 1$ (the cases $t \leq 1$ are easy to estimate). Following the standard procedure of the energy estimate, one needs to estimate
\begin{equation*}
\begin{aligned}
\mathrm{Re} (\partial^{N_0} (u_x u^2), \partial^{N_0} u).
\end{aligned}
\end{equation*}
First, for the term 
\begin{equation*}
\begin{aligned}
\mathrm{Re} (\partial^{N_0} u_x u^2, \partial^{N_0}u),
\end{aligned}
\end{equation*}
integrating by parts leads to
\begin{equation*}
\begin{aligned}
\int (\partial^{N_0} u)^2 (u^2)_x,
\end{aligned}
\end{equation*}
which is easy to deal with. Therefore, it remains to estimate the commutators
\begin{equation*}
\begin{aligned}
\partial^{\alpha_1} u \partial^{\alpha_2} u \partial^{\alpha_3}u,
\end{aligned}
\end{equation*}
where $\alpha_1+\alpha_2+\alpha_3=N_0+1$ and $0 \le \alpha_1,\alpha_2,\alpha_3 \leq N_0$, and show
\begin{equation*}
\begin{aligned}
\norm{\partial^{\alpha_1} u \partial^{\alpha_2} u \partial^{\alpha_3}u}_2 \lesssim \epsilon_1^3 t^{p_0}t^{-1}.
\end{aligned}
\end{equation*}

Decompose $u$ into $u_{\geq 50}=u_h$ and $u_{<50}=u_l$. For the case $(h,h,h)$ (the case $\partial^{\alpha_1} u_h \partial^{\alpha_2} u_h \partial^{\alpha_3}u_h$), by the Coifman-Meyer multiplier theorem, one has
\begin{equation*}
\begin{aligned}
\norm{\partial^{\alpha_1} u_h \partial^{\alpha_2} u_h \partial^{\alpha_3}u_h}_2 \lesssim \norm{\partial^{N_0}u_h}_2 \norm{\partial u_h}_\infty \norm{u_h}_\infty \lesssim \epsilon_1^3 t^{p_0} t^{-1}.
\end{aligned}
\end{equation*}
The case $(l,l,l)$ can be handled easily. For the case $(h,l,l)$, we consider
\begin{equation*}
\begin{aligned}
\partial^{\alpha_1} u_h \partial^{\alpha_2} u_l \partial^{\alpha_3} u_l.
\end{aligned}
\end{equation*}
Since $\alpha_2+\alpha_3 \ge 1$, one has
\begin{equation*}
\begin{aligned}
\norm{\partial^{\alpha_1} u_h \partial^{\alpha_2} u_l \partial^{\alpha_3} u_l}_2 \lesssim \norm{\partial^{\alpha_1} u_h}_2 \norm{\partial^{\alpha_2} u_l \partial^{\alpha_3} u_l}_\infty \lesssim \epsilon_1^3 t^{p_0} t^{-1}.
\end{aligned}
\end{equation*}

Finally, we deal with the case $(h,h,l)$, namely,
\begin{equation*}
\begin{aligned}
\partial^{\alpha_1} u_h \partial^{\alpha_2} u_h \partial^{\alpha_3} u_l,
\end{aligned}
\end{equation*}
which can be decomposed in terms of frequencies as
\begin{equation*}
\begin{aligned}
\sum_{k_1,k_2} \partial^{\alpha_1} u_{k_1} \partial^{\alpha_2} u_{k_2} \partial^{\alpha_3} u_{<50}.
\end{aligned}
\end{equation*}
We assume $k_1 \geq k_2$ and write $\partial^{\alpha_1} u_{k_1} \partial^{\alpha_2} u_{k_2} \partial^{\alpha_3} u_{<50}$ as the multiplier form
\begin{equation*}
\begin{aligned}
\int (\xi-\eta-\sigma)^{\alpha_1} \eta^{\alpha_2} \sigma^{\alpha_3} \Widehat{u_{k_1}}(\xi-\eta-\sigma) \Widehat{u_{k_2}}(\eta) \Widehat{u_{<50}}(\sigma) \, \diff \eta \diff \sigma.
\end{aligned}
\end{equation*}
Then it follows from  Lemma \ref{itmulti}, \eqref{energyasp}, \eqref{decayhl}, \eqref{estlih} and \eqref{estli} that
\begin{equation*}
\begin{aligned}
&\quad\sum_{k_1,k_2}\norm{\partial^{\alpha_1} u_{k_1} \partial^{\alpha_2} u_{k_2} \partial^{\alpha_3} u_{<50}}_2 \\
&\lesssim \sum_{k_1,k_2}2^{\alpha_1 k_1} 2^{\alpha_2 k_2} (\norm{u_{k_1}}_2 I_{t}(u_{k_2},u_{<50})+\langle t\rangle^{-N}\norm{u_{k_1}}_2\norm{u_{k_2}}_\infty\norm{u}_\infty)\\
&\lesssim \sum_{k_1,k_2}(\epsilon_1^3 2^{(\alpha_1-N_0)k_1} 2^{(\alpha_2-5) k_2} t^{p_0} t^{-1}+\epsilon_1^32^{(\alpha_1-N_0)k_1}2^{(\alpha_2-5) k_2}t^{p_0}t^{-5/6-N})\\
&\lesssim \epsilon_1^3 t^{p_0} t^{-1}
\end{aligned}
\end{equation*}
for $\alpha_1<N_0$. The case $\alpha_1=N_0$ can be estimated as follows:
\begin{equation*}
\begin{aligned}
\norm{\partial^{N_0} u_h \partial^{\alpha_2}u_h \partial^{\alpha_3}u_l}_2 \lesssim \norm{\partial^{N_0}u_h}_2 \norm{\partial^{\alpha_2}u_h \partial^{\alpha_3}u_l}_\infty \lesssim \epsilon_1^3 t^{p_0} t^{-1}.
\end{aligned}
\end{equation*}

\section{$Z$ Norms}\label{$Z$ Norms}

We now prove \eqref{Zasp-2} and \eqref{m-scatter}. It suffices to show 
\begin{equation}\label{eq:g-equation}
\begin{aligned}
\norm{(1+\abs \xi^{30})[g(\cdot, t_1)-g(\cdot, t_2)]}_{\infty} \lesssim \epsilon_1^3,
\end{aligned}
\end{equation}
where $g$ is defined by
\[
g(\xi,t)=e^{\mathrm{i} H(\xi,t)} \Widehat f(\xi,t)
\]
with
\[
H(\xi,t)=-\frac{6\pi \xi}{\abs{\Lambda''(\xi)}} \int_0^t \frac{\abs{\Widehat f(\xi,s)}^2}{s} \varphi_{>1}(\abs \xi s^{1/3})\, \diff s.
\]

By \eqref{dpb}, basic calculations yield 
\begin{align*}
\partial_t g(\xi,t)&=e^{\mathrm{i} H(\xi,t)} \big(\partial_t  \Widehat f(\xi,t)+\mathrm{i} \partial_t  H(\xi,t) \Widehat f(\xi,t)\big)\\
&=\mathrm{i} e^{\mathrm{i} H(\xi,t)} \bigg(I(\xi,t)- \frac {6\pi\xi} {\abs{\Lambda''(\xi)}}\frac{\abs{\Widehat f(\xi,t)}^2\hat f(\xi,t)}{t} \varphi_{>1}(\abs \xi t^{1/3})\bigg),
\end{align*}
where
\[
I(\xi,t)=\int e^{-\mathrm{i} t \Phi} \xi \Widehat f(\xi-\eta-\sigma,t) \Widehat f(\eta,t)\Widehat f(\sigma,t) \, \, \diff \eta \diff \sigma.
\]
Therefore to check \eqref{eq:g-equation}, one will show 
\begin{equation}\label{eq:0}
(1+\abs \xi^{30}) \bigg|\int e^{\mathrm{i} H(\xi,s)}\bigg(I(\xi,s)- \frac {6\pi\xi} {\abs{\Lambda''(\xi)}}\frac{\abs{\Widehat f(\xi,s)}^2\hat f(\xi,s)}{s} \varphi_{>1}(\abs \xi s^{1/3})\bigg)\, \diff s\bigg| \lesssim \epsilon_1^3.
\end{equation}

We may assume $\xi>0$ and $\xi \in [2^k,2^{k+1})$ and decompose \eqref{eq:0} into
\begin{equation}\label{eq:0.1}
\sum_m \int q_m(s) e^{\mathrm{i} H(\xi,s)}\bigg(I(\xi,s)- \frac {6\pi\xi} {\abs{\Lambda''(\xi)}}\frac{\abs{\Widehat f(\xi,s)}^2\hat f(\xi,s)}{s}\varphi_{>1}(\abs \xi s^{1/3})\bigg)\, \diff s.
\end{equation}
Except the resonant cases, the correction term
\[
 \frac {6\pi\xi} {\abs{\Lambda''(\xi)}}\frac{\abs{\Widehat f(\xi,t)}^2\hat f(\xi,t)}{t} \varphi_{>1}(\abs \xi t^{1/3})
\]
is not involved. Except using absence of time resonance, $e^{\mathrm{i} H(\xi,t)}$ is directly bounded by $1$. If $m \leq D$ or $m=L+1$, where $D=100$, the terms are easy to be estimated since the difficulty is the summation issue for $m$. Thus we may assume $D<m<L+1$.

First, we can treat easily the case that the highest frequencies are greater than $p_1m=(p_0/10)m$ in \eqref{eq:0.1}. Let
\[
I_{k_1,k_2,k_3,m} \triangleq \int q_m(s) e^{-\mathrm{i} s \Phi+\mathrm{i} H} \xi \Widehat{f_{k_1}}(\xi-\eta-\sigma,s) \Widehat {f_{k_2}}(\eta,s)\Widehat {f_{k_3}}(\sigma,s) \, \, \diff \eta \diff \sigma \diff s.
\]
We aim to estimate
\[
2^{30k_+}\sum_{\max \{k_1,k_2,k_3\} \geq p_1m,m} \abs{I_{k_1,k_2,k_3,m}}.
\]
Assuming that $k_1 \geq k_2 \geq k_3$, one can get
\begin{align*}
&\quad 2^{30k_+}\sum_{k_1 \geq p_1m,m} \abs{I_{k_1,k_2,k_3,m}} \\
&\lesssim 2^{30k_+}\sum_{k_1 \geq p_1m,m}  2^k 2^m \norm{f_{k_1}}_2 \norm{f_{k_2}}_2 \norm{E f_{k_3}}_\infty\\
&\lesssim \epsilon_1^3 \sum_{k_1\geq p_1m,m}  2^{30k_+}2^k 2^m 2^{-N_0{k_1}_+} 2^{k_3}\\
&\lesssim \epsilon_1^3 \sum_{k_1\geq p_1m,m}  2^{31{k_1}_+} 2^m 2^{-(N_0-2){k_1}_+} 2^{-{k_2}_+}2^{k_2/2}2^{-{k_3}_+}2^{k_3/2}
\lesssim \epsilon_1^3,
\end{align*}
where one has used \eqref{energyasp}, \eqref{estl2} and \eqref{estlilow}.

Next, to handle the remaining cases, we need to decompose $I(\xi,s)$ into high frequencies and low   frequencies and analyse the corresponding interactions carefully. 
When $\xi \leq 2^{50}$, we decompose frequencies as follows:
\begin{itemize}
\item
low-low: $I_{\leq 100,\leq 100,\leq 100}$;
\item
high-low: $I_{\alpha,\beta,\gamma}$ where $(100,p_1m)$, $\leq -100 \in \{\alpha,\beta,\gamma\}$ and $\alpha$, $\beta$, $\gamma \in \{(100,p_1m),$ $(-100,100],\leq -100\}$;
\item
high-high: $I_{\alpha,\beta,\gamma}$ where $\alpha$, $\beta$, $\gamma \in \{(100,p_1m),(-100,100]\}$ and at least one of them is $(100,p_1m)$.
\end{itemize}
When $\xi> 2^{50}$, we decompose frequencies as follows:
\begin{itemize}
\item
high-high: $I_{>-100,>-100,>-100}$;
\item
high-low: $I_{\alpha,\beta,\gamma}$ where $(100,p_1m)$, $\leq -100 \in \{\alpha,\beta,\gamma\}$ and $\alpha$, $\beta$, $\gamma \in \{(100,p_1m),$ $(-100,100],\leq -100\}$;
\item
low-low: $I_{\alpha,\beta,\gamma}$ where $\alpha$, $\beta$, $\gamma \in \{\leq -100,(-100,100]\}$ and at least one of them is $\leq -100$.
\end{itemize}

\subsection{High-Low Frequency Interactions}
We only estimate 
\begin{equation*}
\begin{aligned}
I_{(100,p_1m),\beta,\leq -100}&=\int q_m(s) e^{-\mathrm{i} s \Phi+\mathrm{i} H} \xi \Widehat{f_{(100,p_1m)}}(\xi-\eta-\sigma,s) \\
&\quad\times\Widehat{f_{\beta}}(\eta,s) 
\Widehat{f_{\leq -100}}(\sigma,s) \, \, \diff \eta \diff \sigma \diff s
\end{aligned}
\end{equation*}
and the others cases can be handled similarly. Integrating by parts with respect to $\sigma$, one has
\begin{equation*}
\begin{aligned}
|I_{(100,p_1m),\beta,\leq -100}|
 \lesssim \abs{I_1}+\abs{I_2}+\abs{I_3},
\end{aligned}
\end{equation*}
where
\begin{equation*}
\begin{aligned}
I_1&=\int q_m(s) m_1(\xi,\eta,\sigma)e^{-\mathrm{i} s \Phi+\mathrm{i} H} \xi  \partial_\sigma \Widehat{f_{(100,p_1m)}} (\xi-\eta-\sigma,s) \\
&\quad\times\Widehat{f_{\beta}}(\eta,s) \Widehat{f_{\leq -100}}(\sigma,s) \, \, \diff \eta \diff \sigma \diff s, \\
I_2&=\int q_m(s) m_1(\xi,\eta,\sigma)e^{-\mathrm{i} s \Phi+\mathrm{i} H} \xi  \Widehat{f_{(100,p_1m)}} (\xi-\eta-\sigma,s) \\
&\quad\times\Widehat{f_\beta}(\eta,s) \partial_\sigma \Widehat{f_{\leq-100}}(\sigma,s) \, \, \diff \eta \diff \sigma \diff s, \\
I_3&=\int q_m(s) \partial_\sigma m_1(\xi,\eta,\sigma)e^{-\mathrm{i} s \Phi+\mathrm{i} H} \xi  \Widehat{f_{(100,p_1m)}} (\xi-\eta-\sigma,s)\\ 
&\quad\times\Widehat{f_{\beta}}(\eta,s) \Widehat{f_{\leq -100}}(\sigma,s) \, \, \diff \eta \diff \sigma \diff s
\end{aligned}
\end{equation*}
with
\[
m_1=\frac {1} {\mathrm{i} s \partial_\sigma \Phi}.
\]

First we consider $I_1$. By Lemma \ref{resoasym},
\[
\abs{\partial_\sigma \Phi} \gtrsim 1,
\]
which together with Lemma \ref{multicri} yields
\[
\norm{m_1}_{S_{\lesssim p_1m,\lesssim p_1 m}} \lesssim 2^{-m} 2^{3 p_1m}.
\]
Then one may use Theorem \ref{multi}, \eqref{weightasp}, \eqref{estli} and \eqref{estl2} to estimate
\begin{align*}
\abs{I_1} &\lesssim 2^{3 p_1m} 2^k \big\|\partial_\sigma \Widehat{f_{(100,p_1m)}}\big\|_2 \norm{Ef}_\infty \norm{f}_2\\
&\lesssim \epsilon_1^3 2^{3 p_1m} 2^{p_1m} 2^{m/6} 2^{-m/3}.
\end{align*}
$I_2$ can be estimated similarly as $I_1$. 

We finally consider $I_3$. By  Lemma \ref{multicri},
\[
\norm{\partial_\sigma m_1}_{S_{\lesssim p_1m, \lesssim p_1m}} \lesssim 2^{-m} 2^{3 p_1m}.
\]
This, together with  Theorem \ref{multi}, \eqref{estl2} and \eqref{estli}, yields
\begin{align*}
\abs{I_3} &\lesssim 2^{3 p_1m} 2^k \norm{f}_2\norm{Ef}_\infty \norm{f}_2\\
&\lesssim \epsilon_1^3 2^{3 p_1m} 2^{p_1m} 2^{-m/3}.
\end{align*}

\subsection{Low Frequency Interactions}
In this case, one needs to handle the following terms:  
\begin{equation*}
\begin{aligned}
I_{\leq -100,\leq -100,\leq -100}\quad \text{or}\quad I_{(-100,100],\leq -100,\leq -100}.
\end{aligned}
\end{equation*}
Notice that 
\begin{equation*}
\begin{aligned}
I_{(-100,100],\leq -100,\leq -100}=I_{\leq 100,\leq- 100,\leq -100}-I_{\leq -100,\leq -100,\leq -100}.
\end{aligned}
\end{equation*}
It suffices to estimate the generic terms $I_{\leq c_1,\leq c_2,\leq c_3}$.

We may add a cut-off function $\chi$ supported in $\{\abs{\xi-\eta-\sigma} \gtrsim \abs \eta, \abs \sigma\}$ in $I_{\leq c_1,\leq c_2,\leq c_3}$. When $\abs{\xi-\eta-\sigma} \leq 2^{-m/3+D}$, one may write
\begin{equation*}
\begin{aligned}
I_{\lesssim -m/3,\leq c_2,\leq c_3}&=\sum_{k_1 \lesssim -m/3,m} \int q_m(s) e^{-\mathrm{i} s \Phi+\mathrm{i} H} \xi\chi \varphi_{\leq c_2}(\eta)\varphi_{\leq c_3}(\sigma) \\
&\quad\times\Widehat{f_{k_1}}(\xi-\eta-\sigma,s)\Widehat{f_{\lesssim k_1}}(\eta,s) \Widehat{f_{\lesssim k_1}}(\sigma,s) \, \, \diff \eta \diff \sigma \diff s.
\end{aligned}
\end{equation*}
Hence, by \eqref{estl2l} and \eqref{estlilow}, one may estimate
\begin{equation*}
\begin{aligned}
\abs{I_{\lesssim -m/3,\leq c_2,\leq c_3}} &\lesssim \sum_{k_1 \lesssim -m/3,m} 2^m 2^{k}\norm{f_{k_1}}_2 \big\|f_{\lesssim k_1}\big\|_2 \norm{Ef_{\lesssim k_1}}_\infty \\
&\lesssim \epsilon_1^3 \sum_{k_1 \lesssim -m/3,m} 2^m 2^{k} 2^{k_1/2} 2^{k_1/2} 2^{k_1}\\
&\lesssim \epsilon_1^3 \sum_{m} 2^{m/3} 2^{k}
\lesssim \epsilon_1^3.
\end{aligned}
\end{equation*}

Take a partition of unity $\theta_1$, $\theta_2$, $\theta_3$, $\theta_4$ such that 

\begin{itemize}
\item[{\emph {Case 1}}:] $\supp \theta_1\subset \{\abs{\xi-\eta-\sigma} \gg \abs \eta, \abs \sigma\}$,
\item[{\emph {Case 2}}:] $\supp \theta_2\subset \{\abs{\xi-\eta-\sigma}\simeq \abs \eta \gg \abs \sigma\}$,
\item[{\emph {Case 3}}:] $\supp \theta_3\subset \{\abs{\xi-\eta-\sigma}\simeq \abs \sigma \gg \abs \eta\}$,
\item[{\emph {Case 4}}:] $\supp \theta_4\subset \{\abs{\xi-\eta-\sigma}\simeq \abs \eta \simeq \abs \sigma\}$.
\end{itemize}

\subsubsection{Estimates for Case 1}

In this case, $k_1 \simeq k$, and we need to estimate
\begin{equation}\label{eq:2}
\begin{aligned}
&\sum_{m} \int q_m(s) e^{-\mathrm{i} s \Phi+\mathrm{i} H} \xi \chi \theta_1 \varphi_{c_2}(\eta) \varphi_{c_3}(\sigma) \Widehat{f_{k}}(\xi-\eta-\sigma,s)\\
&\quad\times \Widehat{f_{\ll k}}(\eta,s) \Widehat{f_{\ll k}}(\sigma,s)  \, \, \diff \eta \diff \sigma \diff s.
\end{aligned}
\end{equation}
Integrating by parts with respect to $\eta$, one has
\[
\abs{\eqref{eq:2}} \lesssim \abs{I_1}+\abs{I_2}+\abs{I_3},
\]
where
\begin{equation*}
\begin{aligned}
I_1&=\sum_m\int q_m(s) m_1(\xi,\eta,\sigma)e^{-\mathrm{i} s \Phi+\mathrm{i} H} \xi  \partial_\eta \Widehat{f_{k}} (\xi-\eta-\sigma,s)\\
&\quad\times \Widehat{f_{\ll k}}(\eta,s) \Widehat{f_{\ll k}}(\sigma,s)  \, \, \diff \eta \diff \sigma \diff s, \\
I_2&=\sum_m\int q_m(s) m_1(\xi,\eta,\sigma)e^{-\mathrm{i} s \Phi+\mathrm{i} H} \xi  \Widehat{f_{k}} (\xi-\eta-\sigma,s) \\
&\quad\times\partial_\eta \Widehat{f_{\ll k}}(\eta,s) \Widehat{f_{\ll k}}(\sigma,s)  \, \, \diff \eta \diff \sigma \diff s, \\
I_3&=\sum_m\int q_m(s) \partial_\eta m_1(\xi,\eta,\sigma)e^{-\mathrm{i} s \Phi+\mathrm{i} H} \xi  \Widehat{f_{k}} (\xi-\eta-\sigma,s)\\
&\quad\times \Widehat{f_{\ll k}}(\eta,s) \Widehat{f_{\ll k}}(\sigma,s) \, \, \diff \eta \diff \sigma \diff s
\end{aligned}
\end{equation*}
with
\[
m_1=\frac {\chi \theta_1 \varphi_{c_2}\varphi_{c_3}} {\mathrm{i} s \partial_\eta \Phi}.
\]

For $I_1$, by Lemma \ref{resoasym},
\[
\abs{\partial_\eta \Phi} \gtrsim 2^{2k},
\]
which combining with Lemma \ref{multicri} yields
\[
\norm{m_1}_{S_{\ll k,\ll k}} \lesssim 2^{-m} 2^{-2k}.
\]
It follows from Theorem  \ref{multi}, \eqref{estl2l}, \eqref{weightasp} and \eqref{estli} that
\begin{align*}
\abs{I_1} &\lesssim \sum_m 2^k2^{-2k} (2^{-k} \norm{f_k}_2+\norm{xf}_2) \norm{f_{\ll k}}_2 \norm{Ef}_\infty\\
&\lesssim \epsilon_1^3 \sum_m 2^{-m/3} 2^{-k/2} (2^{-k/2}+2^{m/6})
\lesssim \epsilon_1^3.
\end{align*}
$I_2$ can be estimated similarly as $I_1$. 

Considering $I_3$, by Lemma \ref{multicri},
\[
\norm{\partial_\eta m_1}_{S_{\ll k,\ll k}} \lesssim 2^{-m}2^{-3k}.
\]
Using this fact, Theorem  \ref{multi}, \eqref{estl2l} and \eqref{estli}, one has
\begin{align*}
\abs{I_3} &\lesssim \sum_m 2^k 2^{-3k} \norm{f_k}_2 \norm{f_{\ll k}}_2 \norm{Ef}_\infty\\
&\lesssim \epsilon_1^3 \sum_m 2^k 2^{-3k} 2^{k/2} 2^{k/2} 2^{-m/3}
\lesssim \epsilon_1^3.
\end{align*}

\subsubsection{Estimates for Case 2 and Case 3}
These two cases can be handled in a similar fashion, so it suffices to consider {\emph {Case 3}}. We need to estimate
\begin{equation}\label{eq:3}
\begin{aligned}
&\sum_{k_1,m} \int q_m(s) e^{-\mathrm{i} s \Phi+\mathrm{i} H} \xi \chi \theta_3 \varphi_{\leq c_2}(\eta)\Widehat{f_{k_1}}(\xi-\eta-\sigma,s) \\
&\quad\times\Widehat{f_{\ll k_1}}(\eta,s) \Widehat{f_{k_1}}(\sigma,s) \, \, \diff \eta \diff \sigma \diff s.
\end{aligned}
\end{equation}
Integrating by parts with respect to $\eta$  yields
\[
\abs{\eqref{eq:3}} \lesssim \abs{I_1}+\abs{I_2}+\abs{I_3},
\]
where
\begin{equation*}
\begin{aligned}
I_1&=\sum_{k_1,m}\int q_m(s) m_1(\xi,\eta,\sigma)e^{-\mathrm{i} s \Phi+\mathrm{i} H} \xi  \partial_\eta \Widehat{f_{k_1}} (\xi-\eta-\sigma,s)\\
&\quad\times \Widehat{f_{\ll k_1}}(\eta,s) \Widehat{f_{k_1}}(\sigma,s) \, \, \diff \eta \diff \sigma \diff s, \\
I_2&=\sum_{k_1,m}\int q_m(s) m_1(\xi,\eta,\sigma)e^{-\mathrm{i} s \Phi+\mathrm{i} H} \xi  \Widehat{f_{k_1}} (\xi-\eta-\sigma,s) \\
&\quad\times \partial_\eta \Widehat{f_{\ll k_1}}(\eta,s) \Widehat{f_{k_1}}(\sigma,s) \, \, \diff \eta \diff \sigma \diff s, \\
I_3&=\sum_{k_1,m}\int q_m(s) \partial_\eta m_1(\xi,\eta,\sigma)e^{-\mathrm{i} s \Phi+\mathrm{i} H} \xi  \Widehat{f_{k_1}} (\xi-\eta-\sigma,s) \\
&\quad\times \Widehat{f_{\ll k_1}}(\eta,s) \Widehat{f_{k_1}}(\sigma,s) \, \, \diff \eta \diff \sigma \diff s,
\end{aligned}
\end{equation*}
and
\[
m_1=\frac {\chi \theta_3 \varphi_{\leq c_2}} {\mathrm{i} s \partial_\eta \Phi}.
\]

 $I_1$ and $I_2$ can be handled in a similar manner.  By Lemma \ref{resoasym}, 
\[
\abs{\partial_\eta \Phi} \gtrsim 2^{2k_1}.
\]
This together with Lemma \ref{multicri} leads to
\[
\norm{m_1}_{S_{\ll k_1,k_1}} \lesssim 2^{-m} 2^{-2k_1}.
\]
Hence, one uses Theorem \ref{multi}, \eqref{estl2l}, \eqref{weightasp} and \eqref{estlil} to estimate
\begin{align*}
\abs{I_1}&\lesssim \sum_{k_1,m} 2^k 2^{-2k_1} (2^{-k_1}\norm{f_{k_1}}_2+\norm{xf}_2) \norm{f_{\ll k_1}}_2 \norm{Ef_{k_1}}_\infty\\
&\lesssim \epsilon_1^3\sum_{k_1 \gtrsim k \vee -m/3,m} 2^k 2^{-2k_1}(2^{-k_1/2}+2^{m/6}) 2^{-m/2}\\
&\lesssim \epsilon_1^3\sum_{m} (2^{-3k/2} \wedge 2^k2^{5m/6}) 2^{-m/2}+(2^{-k} \wedge 2^k 2^{2m/3}) 2^{-m/3}
\lesssim \epsilon_1^3.
\end{align*}
For $I_3$, by Lemma \ref{multicri},
\[
\norm{\partial_\eta m_1}_{S_{\ll k_1,k_1}} \lesssim 2^{-m}2^{-3k_1}.
\]
It follows from Theorem  \ref{multi}, \eqref{estl2l} and \eqref{estlil} that
\begin{align*}
\abs{I_3} &\lesssim \sum_{k_1,m} 2^k 2^{-3k_1} \norm{f_{k_1}}_2 \norm{f_{\ll k_1}}_2\norm{Ef_{k_1}}_\infty\\
&\lesssim \epsilon_1^3 \sum_{k_1 \gtrsim k \vee -m/3,m} 2^k 2^{-3k_1} 2^{k_1/2} 2^{k_1/2} 2^{-m/2}2^{-k_1/2}\\
&\lesssim \epsilon_1^3 \sum_m (2^{-3k/2} \wedge 2^k 2^{5m/6}) 2^{-m/2}
\lesssim \epsilon_1^3.
\end{align*}

\subsubsection{Estimates for Case 4}
We decompose it into
\[
\sum_{k_1,m} \int q_m(s) e^{-\mathrm{i} s \Phi+\mathrm{i} H} \xi \chi \theta_4 \Widehat{f_{k_1}}(\xi-\eta-\sigma,s) \Widehat{f_{k_1}}(\eta,s) \Widehat{f_{k_1}}(\sigma,s) \, \, \diff \eta \diff \sigma \diff s.
\]
For $k_1 \gg k$, by the elementary inequality
\[
\big|\abs a-\abs b\big| \geq \abs{a-b} \wedge \abs{a+b},
\]
one has
\begin{equation}\label{case2reso}
\big|\abs{\xi-\eta-\sigma}-\abs{\eta}\big| \gtrsim 2^{k_1}\quad \text{if}\ \abs \eta/\abs \sigma >0.9
\end{equation}
and
\begin{equation}
\big|\abs{\xi-\eta-\sigma}-\abs \sigma\big| \gtrsim 2^{k_1}\quad \text{if}\ \abs \sigma/\abs \eta>0.9.
\end{equation}
Take a partition of unity $\chi_1$, $\chi_2$ for these two regions. We only estimate
\begin{equation}\label{eq:4}
\sum_{k_1,m} \int q_m(s) e^{-\mathrm{i} s\Phi+\mathrm{i} H} \chi\theta_4\chi_1 \xi \Widehat{f_{k_1}}(\xi-\eta-\sigma,s) \Widehat{f_{k_1}}(\eta,s) \Widehat{f_{k_1}}(\sigma,s) \, \, \diff \eta \diff \sigma \diff s,
\end{equation}
since the one for $\chi_2$ is similar. Integrating by parts with respect to $\eta$, one has
\[
\abs{\eqref{eq:4}} \lesssim \abs{I_1}+\abs{I_2}+\abs{I_3},
\]
where
\begin{equation*}
\begin{aligned}
I_1&=\sum_{k_1,m}\int q_m(s) m_1(\xi,\eta,\sigma)e^{-\mathrm{i} s \Phi+\mathrm{i} H} \xi  \partial_\eta \Widehat{f_{k_1}} (\xi-\eta-\sigma,s) \\
&\quad\times\Widehat{f_{k_1}}(\eta,s) \Widehat{f_{k_1}}(\sigma,s) \, \diff \eta \diff \sigma \diff s, \\
I_2&=\sum_{k_1,m}\int q_m(s) m_1(\xi,\eta,\sigma)e^{-\mathrm{i} s \Phi+\mathrm{i} H} \xi  \Widehat{f_{k_1}} (\xi-\eta-\sigma,s) \\
&\quad\times\partial_\eta \Widehat{f_{k_1}}(\eta,s) \Widehat{f_{k_1}}(\sigma,s)  \,  \diff \eta \diff \sigma \diff s, \\
I_3&=\sum_{k_1,m}\int q_m(s) \partial_\eta m_1(\xi,\eta,\sigma)e^{-\mathrm{i} s \Phi+\mathrm{i} H} \xi  \Widehat{f_{k_1}} (\xi-\eta-\sigma,s)\\
&\quad\times\Widehat{f_{k_1}}(\eta,s) \Widehat{f_{k_1}}(\sigma,s) \,  \diff \eta \diff \sigma \diff s
\end{aligned}
\end{equation*}
with
\[
m_1=\frac {\chi\theta_4\chi_1} {\mathrm{i} s \partial_\eta \Phi}.
\]

For $I_1$, by \eqref{case2reso}, Lemma \ref{resoasym} and Lemma \ref{multicri},
\[
\norm{m_1}_{S_{k_1,k_1}} \lesssim 2^{-m} 2^{-2k_1}.
\]
It then follows from Theorem  \ref{multi}, \eqref{estl2l}, \eqref{weightasp} and \eqref{estlil} that
\begin{align*}
\abs{I_1}&\lesssim \sum_{k_1,m}2^k2^{-2k_1} (2^{-k_1}\norm{f_{k_1}}_2+\norm{xf}_2) \norm{f_{k_1}}_2 \norm{Ef_{k_1}}_\infty\\
&\lesssim \epsilon_1^3 \sum_{k_1 \gtrsim k \vee -m/3,m} 2^k 2^{-2k_1} (2^{-k_1/2}+2^{m/6}) 2^{k_1/2} 2^{-m/2} 2^{-k_1/2}
\lesssim \epsilon_1^3.
\end{align*}
$I_2$ can be handled similarly as for $I_1$. For $I_3$, notice by \eqref{case2reso}, Lemma \ref{resoasym} and  Lemma \ref{multicri} that
\[
\norm{\partial_\eta m_1}_{S_{k_1,k_1}} \lesssim 2^{-m}2^{-3k_1}.
\]
The remaining arguments are similar to the above ones.

\subsubsection{Resonant Analysis}
It remains to deal with the case $k_1 \simeq k$ that is the resonant region where we will extract the main part
\[
\frac \xi {\Lambda''(\xi)} \frac{\abs{\Widehat f(\xi,t)}^2\hat f(\xi,t)}{t} \varphi_{>1}(\abs \xi t^{1/3})
\]
step by step. Note that this case only happens for $I_{\leq 100,\leq 100,\leq 100}$. We decompose it further into
\[
I^{\iota_1,\iota_2,\iota_3}_{k_1,k_2,k_3,m}=\int q_m(s) e^{-\mathrm{i} s \Phi+\mathrm{i} H} \chi\theta_4\xi \Widehat {f^{\iota_1}_{k_1}}(\xi-\eta-\sigma,s) \Widehat {f^{\iota_2}_{k_2}}(\eta,s)\Widehat {f^{\iota_3}_{k_3}}(\sigma,s) \, \, \diff \eta \diff \sigma \diff s,
\]
where $\Widehat{f^{+}}(\xi)=\Widehat f(\xi) 1_{>0}(\xi)$ and $\Widehat{f^{-}}(\xi)=\Widehat f(\xi) 1_{\leq 0}(\xi)$. Except the final resonant region, we may regard $k_i$ as $k$ and ignore the summation for $k_i$.\\

\paragraph{{\emph{Time Resonant Analysis}}}
For $(\iota_1,\iota_2,\iota_3)\in \{(+,+,+)$, $(+,-,-)$, $(-,+,-)$, $(-,-,+)$, $(-,-,-)\}$, we can use the absence of time resonance. Indeed, by \eqref{resoin2}, 
\begin{equation}\label{case3resotime}
\abs \Phi \gtrsim 2^{3k}.
\end{equation}
Integrating by parts with respect to $s$, one has
\begin{equation}\label{eq:5}
\sum_m \abs{I^{\iota_1,\iota_2,\iota_3}_{m}} \lesssim I_1+I_2+\text{similar terms},
\end{equation}
where
\begin{equation*}
\begin{aligned}
I_1&=\sum_m \bigg|\int q_m(s) m_1(\xi,\eta,\sigma) e^{-\mathrm{i} s \Phi+\mathrm{i} H} \xi \partial_s \Widehat {f^{\iota_1}_{k}}(\xi-\eta-\sigma,s) \\
&\quad\times\Widehat {f^{\iota_2}_{k}}(\eta,s)\Widehat {f^{\iota_3}_{k}}(\sigma,s)  \,  \diff \eta \diff \sigma \diff s\bigg|\\
I_2&=\sum_m \bigg|\int q'_m(s) m_1(\xi,\eta,\sigma) e^{-\mathrm{i} s \Phi+\mathrm{i} H} \xi \Widehat {f^{\iota_1}_{k}}(\xi-\eta-\sigma,s) \\
&\quad\times\Widehat {f^{\iota_2}_{k}}(\eta,s)\Widehat {f^{\iota_3}_{k}}(\sigma,s) \,  \diff \eta \diff \sigma \diff s\bigg|
\end{aligned}
\end{equation*}
with
\[
m_1=\frac{\chi\theta_4}{-\mathrm{i} \Phi+\mathrm{i} H_s}.
\]

For $I_1$, by \eqref{case3resotime} and Lemma \ref{multicri},
\[
\norm{m_1}_{S_{k,k}} \lesssim 2^{-3k}.
\]
It follows from Theorem  \ref{multi}, \eqref{psl}, \eqref{estl2l} and \eqref{estlil} that
\begin{align*}
I_1 &\lesssim \sum_m 2^m 2^{-3k} 2^k \norm{\partial_s f^{\iota_1}_{k}}_2 \norm{f_{k}}_2 \norm{Ef_{k}}_\infty\\
&\lesssim \epsilon_1^3 \sum_m 2^m 2^{-3k} 2^k 2^{k/2} 2^{-m} 2^{k/2} 2^{-m/2} 2^{-k/2}\\
&\lesssim \epsilon_1^3 \sum_m 2^{-m/2}2^{-3k/2}
\lesssim \epsilon_1^3.
\end{align*}
For $I_2$, one has
\begin{align*}
I_2 &\lesssim \sum_m 2^k 2^{-3k}\norm{f_{k}}_2 \norm{f_{k}}_2 \norm{Ef_{k}}_\infty\\
&\lesssim \epsilon_1^3 \sum_m 2^k 2^{-3k}2^{k/2} 2^{k/2} 2^{-m/2}2^{-k/2}
\lesssim \epsilon_1^3 .
\end{align*}

\paragraph{{\emph{Space Resonant Analysis}}}
It remains to consider $(\iota_1,\iota_2,\iota_3)\in \{(-,+,+)$, $(+,-,+)$, $(+,+,-)\}$. We only handle the case $(-,+,+)$ here, since the others are similar. We decompose it further into
\begin{equation*}
\begin{aligned}
&\sum_{m,l_1,l_2}\int q_m(s) e^{-\mathrm{i} s \Phi+\mathrm{i} H} \chi\theta_4 \xi \Widehat {f^{-}_{k}}(-\xi-\eta-\sigma,s) \Widehat {f^{+}_{k}}(\xi+\eta,s)\\
&\quad\times\Widehat {f^{+}_{k}}(\xi+\sigma,s) \varphi^{\bar l}_{l_1}(\eta) \varphi^{\bar l}_{l_2}(\sigma)\,  \diff \eta \diff \sigma \diff s,
\end{aligned}
\end{equation*}
where $\bar l=-29m/60-9k/20-D$. We start with
\begin{equation}\label{eq:6}
\begin{aligned}
&\sum_{l_1 \vee l_2>\bar l,m} \int q_m(s) e^{-\mathrm{i} s \Phi+\mathrm{i} H} \chi\theta_4\xi \Widehat {f^{-}_{k}}(-\xi-\eta-\sigma,s) \Widehat {f^{+}_{k}}(\xi+\eta,s)\\
&\quad\times\Widehat {f^{+}_{k}}(\xi+\sigma,s) \varphi^{\bar l}_{l_1}(\eta) \varphi^{\bar l}_{l_2}(\sigma)\,  \diff \eta \diff \sigma \diff s.
\end{aligned}
\end{equation}
We may assume $l_2 \geq l_1$ and use the absence of space resonance. Integrating by parts with respect to $\eta$ gives
\[
\abs{\eqref{eq:6}} \lesssim I_1+I_2+I_3
\]
where
\begin{equation*}
\begin{aligned}
I_1&=\sum_{l_2 \geq (\bar l+1) \vee l_1,m}\bigg|\int q_m(s) m_1(\xi,\eta,\sigma) e^{-\mathrm{i} s \Phi+\mathrm{i} H} \xi \partial_\eta \Widehat {f^{-}_{k}}(-\xi-\eta-\sigma,s),\\
&\quad\times \Widehat {f^{+}_{k}}(\xi+\eta,s)\Widehat {f^{+}_{k}}(\xi+\sigma,s) \,  \diff \eta \diff \sigma \diff s\bigg|,\\
I_2&=\sum_{l_2 \geq (\bar l+1) \vee l_1,m}\bigg|\int q_m(s) m_1(\xi,\eta,\sigma) e^{-\mathrm{i} s \Phi+\mathrm{i} H} \xi \Widehat {f^{-}_{k}}(-\xi-\eta-\sigma,s)\\ &\quad\times\partial_\eta \Widehat {f^{+}_{k}}(\xi+\eta,s)\Widehat {f^{+}_{k}}(\xi+\sigma,s)  \, \diff \eta \diff \sigma \diff s\bigg|,\\
I_3&=\sum_{l_2 \geq (\bar l+1) \vee l_1,m}\bigg|\int q_m(s) \partial_\eta m_1(\xi,\eta,\sigma) e^{-\mathrm{i} s \Phi+\mathrm{i} H} \xi \Widehat {f^{-}_{k}}(-\xi-\eta-\sigma,s) \\
&\quad\times\Widehat {f^{+}_{k}}(\xi+\eta,s)\Widehat {f^{+}_{k}}(\xi+\sigma,s)\,  \diff \eta \diff \sigma \diff s\bigg|
\end{aligned}
\end{equation*}
with
\[
m_1=\frac{\chi \theta_4 \varphi^{\bar l}_{l_1}(\eta) \varphi^{\bar l}_{l_2}(\sigma)}{-\mathrm{i} s \partial_\eta \Phi}.
\]

By Lemma \ref{resoasym},
\[
\abs{\partial_\eta \Phi} \gtrsim 2^{l_2} 2^k,
\]
hence using Lemma \ref{multicri} one has
\[
\norm{m_1}_S \lesssim 2^{-m} 2^{-l_2} 2^{-k}.
\]
For $I_3$, it follows from Theorem  \ref{multi}, \eqref{weightasp}, \eqref{estl2l} and \eqref{estlil} that
\begin{align*}
I_1 &\lesssim \sum_{l_2,l_1,m} 2^{-l_2} 2^{-k} 2^k \big(\norm{xf}_2+2^{-k}\big\|f_{B(\xi,2^{l_2+4})}\big\|_2\big) \big\|f_{B(\xi,2^{l_1+4})}\big\|_2 \norm{Ef_k}_\infty\\
&\lesssim \epsilon_1^3 \sum_{l_2,l_1,m} 2^{-l_2} (2^{m/6}+2^{l_2/2-k}) 2^{l_1/2} 2^{-m/2} 2^{-k/2}\\
&\lesssim \epsilon_1^3 \sum_{l_2,m} 2^{-l_2/2} (2^{m/6}+2^{l_2/2-k})2^{-m/2} 2^{-k/2}\\
&\lesssim \epsilon_1^3 \sum_{l_2,m} 2^{-l_2} 2^{-m/2} 2^{-k/2}+2^{-l_2/2}2^{-m/3}2^{-k/2}\\
&\lesssim \epsilon_1^3 \sum_m 2^{-\bar l} 2^{-m/2}2^{-k/2}+2^{-\bar l/2}2^{-m/3}2^{-k/2}
\lesssim \epsilon_1^3.
\end{align*}
$I_2$ can be handled similarly as $I_1$. For $I_3$, by Lemma \ref{multicri},
\[
\norm{\partial_\eta m_1}_S \lesssim 2^{-m} 2^{-l_1} 2^{-l_2} 2^{-k}.
\]
which together with Theorem  \ref{multi}, \eqref{estl2l} and \eqref{estlil} yields
\begin{align*}
I_3 &\lesssim \sum_{l_2,l_1,m} 2^{-l_1} 2^{-l_2} 2^{-k} 2^k \big\|f_{B(\xi,2^{l_2+4})}\big\|_2 \big\|f_{B(\xi,2^{l_1+4})}\big\|_2 \norm{Ef_k}_\infty\\
&\lesssim \epsilon_1^3\sum_{l_2,l_1,m} 2^{-l_1} 2^{-l_2} 2^{l_2/2} 2^{l_1/2} 2^{-m/2} 2^{-k/2}\\
&\lesssim \epsilon_1^3 \sum_m 2^{-\bar l} 2^{-m/2} 2^{-k/2}
\lesssim \epsilon_1^3.
\end{align*}

Finally, we are ready to show 
 \begin{equation}\label{resonant term}
\begin{aligned}
&\bigg|\sum_{k_1,k_2,k_3,m} \int q_m(s) e^{\mathrm{i} H}\bigg( \int e^{-\mathrm{i} s \Phi} \chi\theta_4\xi \Widehat {f^{-}_{k_1}}(-\xi-\eta-\sigma,s) 
 \\ 
&\quad \times \Widehat {f^{+}_{k_2}}(\xi+\eta,s)\Widehat {f^{+}_{k_3}}(\xi+\sigma,s)\phi_{\bar l}(\eta) \varphi_{\bar l}(\sigma)\, \, \diff \eta \diff \sigma\\
 &\quad- \frac {2\pi\xi} {\abs{\Lambda''(\xi)}}\frac{\Widehat{f^-_{k_1}}(-\xi,s)\Widehat{f^+_{k_2}}(\xi,s)\Widehat{f^+_{k_3}}(\xi,s)}{s}  \varphi_{>1}(\abs \xi s^{1/3})\bigg)\, \diff s\bigg|\lesssim \epsilon_1^3,
\end{aligned}
\end{equation}
where $k_1>-m/3+D$, $k_1$, $k_2$, $k_3 \simeq k$, and $k \gtrsim -m/3$.

Notice that
\[
\abs{\Phi(\xi,\eta,\sigma)-\Lambda''(\xi)\eta\sigma} \lesssim (\abs \eta+\abs \sigma)^3.
\]
Then one has
\begin{equation}\label{resonant estimate-1}
\begin{aligned}
&\sum_m \bigg| \int q_m(s)e^{\mathrm{i} H} (e^{-\mathrm{i} s \Phi}-e^{-\mathrm{i} s \Lambda''(\xi)\eta\sigma}) \chi\theta_4\xi \Widehat {f^{-}_{k}}(-\xi-\eta-\sigma,s) \\
&\quad\times \Widehat {f^{+}_{k}}(\xi+\eta,s)\Widehat {f^{+}_{k}}(\xi+\sigma,s) \varphi_{\bar l}(\eta) \varphi_{\bar l}(\sigma)\,  \diff \eta \diff \sigma \diff s \bigg|\\
&\lesssim \sum_m 2^{2m} 2^{3\bar l}2^k \big\|\Widehat f\big\|_\infty^3 2^{2\bar l}
\lesssim \epsilon_1^3 \sum_m 2^{2m} 2^{3\bar l} 2^k 2^{2\bar l}
\lesssim \epsilon_1^3.
\end{aligned}
\end{equation}
By \eqref{weightasp} and \eqref{estl2l}, one notices
\begin{align*}
&\big|\Widehat{f^-_{k_1}}(-\xi-\eta-\sigma)-\Widehat{f^-_{k_1}}(-\xi)\big|
=\bigg|\int_{-\xi}^{-\xi-\eta-\sigma} \partial \Widehat{f_{k_1}} (\zeta)\, \diff \zeta\bigg|\\
&\lesssim 2^{\bar l/2} (\norm{xf}_2+2^{-k_1}\norm{f_{k_1}}_2)
\lesssim \epsilon_1 2^{\bar l/2}(2^{m/6}+2^{-k/2}),
\end{align*}
(the terms for $k_2$ and $k_3$ are similar), and then estimates
\begin{equation}\label{resonant estimate-2}
\begin{aligned}
&\sum_{k_1,k_2,k_3,m} \bigg| \int q_m(s) e^{\mathrm{i} H}e^{-\mathrm{i} s \Lambda''(\xi)\eta\sigma} \chi\theta_4\xi \bigg(\Widehat {f^{-}_{k_1}}(-\xi-\eta-\sigma,s) \Widehat {f^{+}_{k_2}}(\xi+\eta,s)\\
&\quad \times\Widehat {f^{+}_{k_3}}(\xi+\sigma,s)-\Widehat{f^-_{k_1}}(-\xi,s)\Widehat {f^{+}_{k_2}}(\xi,s) \Widehat {f^{+}_{k_3}}(\xi,s)\bigg) \varphi_{\bar l}(\eta) \varphi_{\bar l}(\sigma)\, \, \diff \eta \diff \sigma \diff s \bigg|\\
&\lesssim \epsilon_1^3 \sum_m 2^m 2^k2^{\bar l/2}(2^{m/6}+2^{-k/2}) 2^{2\bar l}
\lesssim \epsilon_1^3.
\end{aligned}
\end{equation}
One calculates 
\begin{equation*}
\begin{aligned}
\int e^{-\mathrm{i} \lambda \eta \sigma} F(\eta,\sigma) \, \, \diff \eta \diff \sigma&=\int e^{\mathrm{i} \frac{\eta\sigma}{\lambda}} \frac{2\pi}{\lambda} \Widehat F(\eta,\sigma) \, \, \diff \eta \diff \sigma\\
&=\frac{2\pi}{\lambda} F(0)+O \bigg(\frac{\big\|\Widehat F(\eta,\sigma) \eta \sigma\big\|_1}{\lambda^2}\bigg),
\end{aligned}
\end{equation*}
and
\begin{equation*}
\begin{aligned}
\int e^{-\mathrm{i} s \Lambda''(\xi) \eta \sigma} \chi\theta_4\phi_{\bar l}(\eta) \varphi_{\bar l}(\sigma) \, \diff \eta \diff \sigma=\frac{2\pi}{s\Lambda''(\xi)}+O(2^{-2\bar l}2^{-2m} 2^{-2k}).
\end{aligned}
\end{equation*}
We then may estimate
\begin{equation}\label{resonant estimate-3}
\begin{aligned}
&\sum_{k_1,k_2,k_3,m} \bigg| \int q_m(s) e^{\mathrm{i} H}\bigg(\int e^{-\mathrm{i} s \Lambda''(\xi) \eta\sigma} \chi\theta_4\phi_{\bar l}(\eta) \varphi_{\bar l}(\sigma) \, \diff \eta\diff \sigma  \\
&\quad-\frac{2\pi}{s\Lambda''(\xi)} \varphi_{>1}(\abs \xi s^{1/3})\bigg)\xi\hat{f^-_{k_1}}(-\xi,s)\Widehat {f^{+}_{k_2}}(\xi,s) \Widehat {f^{+}_{k_3}}(\xi,s)\, \diff s \bigg| \lesssim \epsilon_1^3.
\end{aligned}
\end{equation}
Hence \eqref{resonant term} follows from \eqref{resonant estimate-1}-\eqref{resonant estimate-3}.

\subsection{High Frequency Interactions}
In this section, we deal with high-high interactions and only consider generic terms $I_{(c_1,p_1m),(c_2,p_1m),(c_3,p_1m)}$. These terms could be estimated similarly as the low-low interaction terms and are actually easier since they are not  sharp. We decompose it into
\[
(1+\abs{\xi}^{30})\xi\sum_{k_1,k_2,k_3,m} \int q_m(s) e^{-\mathrm{i} s \Phi+\mathrm{i} H} \Widehat{f_{k_1}}(\xi-\eta-\sigma,s) \Widehat{f_{k_2}}(\eta,s) \Widehat{f_{k_3}}(\sigma,s) \, \, \diff \eta \diff \sigma \diff s.
\]
First, one has the resonant region 
\begin{equation*}
\begin{aligned}
\max \{\abs{k_1-k},\abs{k_2-k},\abs{k_3-k}\} \leq 20,
\end{aligned}
\end{equation*}
which only occurs in $I_{>-100,>-100,>-100}$. For the non-resonant region, one can take absolute value for each term. We may assume $k_1 \geq k_2 \geq k_3$ and then divide them into three cases:
\begin{equation*}
\begin{aligned}
k_1 \gg k_2 \geq k_3, \quad k_1 \simeq k_2 \gg k_3,\quad k_1 \simeq k_2 \simeq k_3.
\end{aligned}
\end{equation*}

\subsubsection{$k_1 \gg k_2 \geq k_3$}
In this case, $k_1 \simeq k$. One needs to estimate
\begin{equation}\label{eq:3.3.1}
\begin{aligned}
&(1+\abs \xi^{30})\xi \sum_{k_2,k_3,m} \int q_m(s) e^{-\mathrm{i} s \Phi+\mathrm{i} H} \Widehat{f_{k}}(\xi-\eta-\sigma,s) \\
&\quad\times\Widehat{f_{k_2}}(\eta,s) \Widehat{f_{k_3}}(\sigma,s)  \,  \diff \eta \diff \sigma \diff s.
\end{aligned}
\end{equation}
Integrating by parts with respect to $\eta$, one has
\[
\abs{\eqref{eq:3.3.1}} \lesssim \abs{I_1}+\abs{I_2}+\abs{I_3},
\]
where
\begin{equation*}
\begin{aligned}
I_1&=(1+\abs \xi^{30})\xi\sum_{k_2,k_3,m}\int q_m(s) m_1(\xi,\eta,\sigma)e^{-\mathrm{i} s \Phi+\mathrm{i} H}  \partial_\eta \Widehat{f_{k}} (\xi-\eta-\sigma,s)\\
&\quad\times \Widehat{f_{k_2}}(\eta,s) \Widehat{f_{k_3}}(\sigma,s)  \,  \diff \eta \diff \sigma \diff s, \\
I_2&=(1+\abs \xi^{30})\xi\sum_{k_2,k_3,m}\int q_m(s) m_1(\xi,\eta,\sigma)e^{-\mathrm{i} s \Phi+\mathrm{i} H} \Widehat{f_{k}} (\xi-\eta-\sigma,s) \\
&\quad\times\partial_\eta \Widehat{f_{k_2}}(\eta,s) \Widehat{f_{k_3}}(\sigma,s)   \, \diff \eta \diff \sigma \diff s, \\
I_3&=(1+\abs \xi^{30})\xi\sum_{k_2,k_3,m}\int q_m(s) \partial_\eta m_1(\xi,\eta,\sigma)e^{-\mathrm{i} s \Phi+\mathrm{i} H}  \Widehat{f_{k}} (\xi-\eta-\sigma,s)\\ &\quad\times\Widehat{f_{k_2}}(\eta,s) \Widehat{f_{k_3}}(\sigma,s)  \,  \diff \eta \diff \sigma \diff s
\end{aligned}
\end{equation*}
with
\[
m_1=\frac {1} {\mathrm{i} s \partial_\eta \Phi}.
\]

For $I_1$, by Lemma \ref{resoasym},
\[
|\partial_\eta \Phi| \gtrsim 2^{-k_2/2},
\]
which together with Lemma \ref{multicri} yields
\[
\norm{m_1}_{S_{k,k_2,k_3}} \lesssim 2^{-m} 2^{k_2/2}.
\]
It follows from Theorem \ref{multi}, \eqref{weightasp2}, \eqref{energyasp} and \eqref{estlih} that
\begin{align*}
\abs{I_1} &\lesssim \sum_{k_2,k_3,m} 2^{31k}2^{k_2/2} (\norm{f_k}_2+\norm{(xf)_k}_2) \norm{f_{k_2}}_2 \norm{Ef_{k_3}}_\infty\\
&\lesssim \epsilon_1^3 \sum_{k_2,k_3,m} 2^{31p_1m}2^{p_2m}2^{p_0m}2^{(-N_0+1/2)k_2}2^{-m/2}
\lesssim \epsilon_1^3.
\end{align*}
$I_2$ can be estimated similarly as $I_1$. For $I_3$, by Lemma \ref{multicri},
\[
\norm{\partial_\eta m_1}_{S_{k,k_2,k_3}} \lesssim 2^{-m}2^{-k_2/2}.
\]
Using Theorem \ref{multi}, \eqref{energyasp} and \eqref{estlih}, one has
\begin{align*}
\abs{I_3} &\lesssim \sum_{k_2,k_3,m} 2^{31k}2^{-k_2/2} \norm{f_k}_2 \norm{f_{k_2}}_2 \norm{Ef_{k_3}}_\infty\\
&\lesssim \epsilon_1^3 \sum_{k_2,k_3,m} 2^{31k} 2^{p_0m}2^{-N_0k} 2^{-k_2/2} 2^{-m/2}
\lesssim \epsilon_1^3.
\end{align*}

\subsubsection{$k_1 \simeq k_2 \gg k_3$}
We need to estimate
\begin{equation}\label{eq:3.3.2}
\begin{aligned}
&(1+\abs \xi^{30})\xi\sum_{k_1,k_3,m} \int q_m(s) e^{-\mathrm{i} s \Phi+\mathrm{i} H} \Widehat{f_{k_1}}(\xi-\eta-\sigma,s) \\
&\quad\times\Widehat{f_{k_1}}(\eta,s) \Widehat{f_{k_3}}(\sigma,s) \, \, \diff \eta \diff \sigma \diff s.
\end{aligned}
\end{equation}
Integrating by parts with respect to $\sigma$ gives
\[
\abs{\eqref{eq:3.3.2}} \lesssim \abs{I_1}+\abs{I_2}+\abs{I_3},
\]
where
\begin{equation*}
\begin{aligned}
I_1&=(1+\abs \xi^{30})\xi\sum_{k_1,k_3,m}\int q_m(s) m_1(\xi,\eta,\sigma)e^{-\mathrm{i} s \Phi+\mathrm{i} H} \partial_\sigma \Widehat{f_{k_1}} (\xi-\eta-\sigma,s) \\
&\quad\times\Widehat{f_{k_1}}(\eta,s) \Widehat{f_{k_3}}(\sigma,s) \,  \diff \eta \diff \sigma \diff s, \\
I_2&=(1+\abs \xi^{30})\xi\sum_{k_1,k_3,m}\int q_m(s) m_1(\xi,\eta,\sigma)e^{-\mathrm{i} s \Phi+\mathrm{i} H}  \Widehat{f_{k_1}} (\xi-\eta-\sigma,s)\\ &\quad\times\Widehat{f_{k_1}}(\eta,s) \partial_\sigma \Widehat{f_{k_3}}(\sigma,s) \,  \diff \eta \diff \sigma \diff s, \\
I_3&=(1+\abs \xi^{30})\xi\sum_{k_1,k_3,m}\int q_m(s) \partial_\sigma m_1(\xi,\eta,\sigma)e^{-\mathrm{i} s \Phi+\mathrm{i} H}  \Widehat{f_{k_1}} (\xi-\eta-\sigma,s)\\
&\quad\times \Widehat{f_{k_1}}(\eta,s) \Widehat{f_{k_3}}(\sigma,s) \,  \diff \eta \diff \sigma \diff s
\end{aligned}
\end{equation*}
with
\[
m_1=\frac {1} {\mathrm{i} s \partial_\sigma \Phi}.
\]

For $I_1$, by Lemma \ref{resoasym}, it holds that
\[
\abs{\partial_\sigma \Phi} \gtrsim 2^{-k_3/2}
\]
which combined with Lemma \ref{multicri} yields
\[
\norm{m_1}_{S_{k_1,k_1,k_3}} \lesssim 2^{-m} 2^{k_3/2}.
\]
Using Theorem \ref{multi}, \eqref{weightasp2}, \eqref{energyasp} and \eqref{estlih}, one has
\begin{align*}
\abs{I_1}&\lesssim \sum_{k_1,k_3,m} 2^{30k_++k}2^{k_3/2} (\norm{f_{k_1}}_2+\norm{(xf)_{k_1}}_2) \norm{f_{k_1}}_2 \norm{Ef_{k_3}}_\infty\\
&\lesssim \epsilon_1^3 \sum_{k_1,k_3,m} 2^{31k_+} 2^{k_3/2} 2^{p_2m} 2^{p_0m}2^{-N_0k_1}2^{-m/2}
\lesssim \epsilon_1^3.
\end{align*}
$I_2$ can be handled similarly as $I_1$. For $I_3$, Lemma \ref{multicri} implies
\[
\norm{\partial_\sigma m_1}_{S_{k_1,k_1,k_3}} \lesssim 2^{-m}2^{-k_3/2}.
\]
It follows from Theorem \ref{multi}, \eqref{energyasp} and \eqref{estlih} that
\begin{align*}
\abs{I_3} &\lesssim \sum_{k_1,k_3,m} 2^{30k_++k}2^{-k_3/2} \norm{f_{k_1}}_2 \norm{f_{k_1}}_2\norm{Ef_{k_3}}_\infty\\
&\lesssim \epsilon_1^3 \sum_{k_1,k_3,m} 2^{31k_+}2^{-k_3/2}2^{p_0m}2^{-N_0k_1}2^{-m/2}
\lesssim \epsilon_1^3.
\end{align*}

\subsubsection{$k_1 \simeq k_2 \simeq k_3$}
Similar to low frequency cases, one can get
$$\big|\abs{\xi-\eta-\sigma}-\abs \eta\big| \gtrsim 2^{k_1}\quad \text{if}\ \abs \eta/\abs \sigma >0.9$$ 
and 
$$\big|\abs{\xi-\eta-\sigma}-\abs \sigma\big| \gtrsim 2^{k_1}\quad \text{if}\ \abs \sigma/\abs \eta>0.9.$$
 Take a partition of unity $\chi_1$, $\chi_2$ for these two regions. We estimate
\begin{equation}\label{eq:3.3.3.1}
\begin{aligned}
&(1+\abs \xi^{30})\xi\sum_{k_1,m} \int q_m(s) e^{-\mathrm{i} s\Phi+\mathrm{i} H} \chi_1 \Widehat{f_{k_1}}(\xi-\eta-\sigma,s) \\
&\quad\times\Widehat{f_{k_1}}(\eta,s) \Widehat{f_{k_1}}(\sigma,s)  \, \diff \eta \diff \sigma \diff s,
\end{aligned}
\end{equation}
and the one for $\chi_2$ is similar. Integrating by parts with respect to $\eta$, one has
\[
\abs{\eqref{eq:3.3.3.1}} \lesssim \abs{I_1}+\abs{I_2}+\abs{I_3},
\]
where
\begin{equation*}
\begin{aligned}
I_1&=(1+\abs \xi^{30})\xi\sum_{k_1,m}\int q_m(s) m_1(\xi,\eta,\sigma)e^{-\mathrm{i} s \Phi+\mathrm{i} H} \partial_\eta \Widehat{f_{k_1}} (\xi-\eta-\sigma,s)\\ 
&\quad\times\Widehat{f_{k_1}}(\eta,s) \Widehat{f_{k_1}}(\sigma,s) \, \diff \eta \diff \sigma \diff s, \\
I_2&=(1+\abs \xi^{30})\xi\sum_{k_1,m}\int q_m(s) m_1(\xi,\eta,\sigma)e^{-\mathrm{i} s \Phi+\mathrm{i} H} \Widehat{f_{k_1}} (\xi-\eta-\sigma,s) \\
&\quad\times\partial_\eta \Widehat{f_{k_1}}(\eta,s) \Widehat{f_{k_1}}(\sigma,s)   \, \diff \eta \diff \sigma \diff s, \\
I_3&=(1+\abs \xi^{30})\xi\sum_{k_1,m}\int q_m(s) \partial_\eta m_1(\xi,\eta,\sigma)e^{-\mathrm{i} s \Phi+\mathrm{i} H}  \Widehat{f_{k_1}} (\xi-\eta-\sigma,s)
\\ 
&\quad\times\Widehat{f_{k_1}}(\eta,s) \Widehat{f_{k_1}}(\sigma,s) \, \diff \eta \diff \sigma \diff s
\end{aligned}
\end{equation*}
with
\[
m_1=\frac {\chi_1} {\mathrm{i} s \partial_\eta \Phi}.
\]

For $I_1$, by Lemma \ref{resoasym} and Lemma \ref{multicri}, one gets
\[
\norm{m_1}_{S_{k_1,k_1,k_1}} \lesssim 2^{-m} 2^{k_1/2}.
\]
It follows from Theorem \ref{multi}, \eqref{weightasp2}, \eqref{energyasp}, \eqref{estlih} that
\begin{align*}
\abs{I_1}&\lesssim \sum_{k_1,m}2^{30k_++k} 2^{k_1/2}(\norm{f_{k_1}}_2+\norm{(xf)_{k_1}}_2) \norm{f_{k_1}}_2 \norm{Ef_{k_1}}_\infty\\
&\lesssim \epsilon_1^3 \sum_{k_1,m} 2^{31k_+}2^{k_1/2} 2^{p_2m}2^{p_0m}2^{-N_0k_1}2^{-m/2}
\lesssim \epsilon_1^3.
\end{align*}
$I_2$ is similar. For $I_3$, Lemma \ref{resoasym} and Lemma \ref{multicri} give
\[
\norm{\partial_\eta m_1}_{S_{k_1,k_1,k_1}} \lesssim 2^{-m}2^{-k_1/2}.
\]
The remaining arguments are similar to the above ones.

Similar to low frequency cases, we will extract the main part
\[
\frac \xi {\Lambda''(\xi)} \frac{\abs{\Widehat f(\xi,t)}^2\hat f(\xi,t)}{t} \varphi_{>1}(\abs \xi t^{1/3})
\]
step by step. We also decompose it into $I^{\iota_1,\iota_2,\iota_3}_{k_1,k_2,k_3,m}$. \\

\paragraph{{\emph{Time Resonant Analysis}}}
For $$(\iota_1,\iota_2,\iota_3)=(+,+,+), (+,-,-), (-,+,-), (-,-,+), (-,-,-),$$ we can use the absence of time resonance, since by \eqref{resoin2},
\begin{equation}\label{3.3.4.1}
\abs \Phi \gtrsim 2^{k/2}.
\end{equation}
Integrating by parts with respect to $s$, one has
\[
\sum_m \abs{I^{\iota_1,\iota_2,\iota_3}_{m}} \lesssim I_1+I_2+\text{similar terms},
\]
where
\begin{equation*}
\begin{aligned}
I_1&=\sum_m \bigg|(1+\abs \xi^{30})\xi\int q_m(s) m_1(\xi,\eta,\sigma) e^{-\mathrm{i} s \Phi+\mathrm{i} H} \partial_s \Widehat {f^{\iota_1}_{k}}(\xi-\eta-\sigma,s) \\
&\quad\times\Widehat {f^{\iota_2}_{k}}(\eta,s)\Widehat {f^{\iota_3}_{k}}(\sigma,s)  \, \, \diff \eta \diff \sigma \diff s\bigg|\\
I_2&=\sum_m \bigg|(1+\abs \xi^{30})\xi\int q'_m(s) m_1(\xi,\eta,\sigma) e^{-\mathrm{i} s \Phi+\mathrm{i} H} \xi \Widehat {f^{\iota_1}_{k}}(\xi-\eta-\sigma,s)\\
&\quad\times \Widehat {f^{\iota_2}_{k}}(\eta,s)\Widehat {f^{\iota_3}_{k}}(\sigma,s) \, \, \diff \eta \diff \sigma \diff s\bigg|,
\end{aligned}
\end{equation*}
and
\[
m_1=\frac{1}{-\mathrm{i} \Phi+\mathrm{i} H_s}.
\]
For $I_1$, by \eqref{3.3.4.1} and Lemma \ref{multicri},
\[
\norm{m_1}_{S_{k,k,k}} \lesssim 2^{-k/2}.
\]
It follows from Theorem \ref{multi}, \eqref{psl}, \eqref{estl2}, \eqref{estlih} that
\begin{align*}
I_1 &\lesssim \sum_m 2^{30k_++k}2^m 2^{-k/2} \norm{\partial_s f^{\iota_1}_{k}}_2 \norm{f_{k}}_2 \norm{Ef_{k}}_\infty\\
&\lesssim \epsilon_1^3 \sum_m 2^{30k_+}2^{k/2}2^m 2^{k/2}2^{-m}2^{-m/2}
\lesssim \epsilon_1^3.
\end{align*}
For $I_2$, similarly, one has
\begin{align*}
I_2 &\lesssim \sum_m 2^{30k_++k}2^{-k/2}\norm{f_{k}}_2 \norm{f_{k}}_2 \norm{Ef_{k}}_\infty\\
&\lesssim \epsilon_1^3 \sum_m 2^{30k_+}2^{k/2}2^{-m/2}
\lesssim \epsilon_1^3.
\end{align*}

\paragraph{{\emph{Space Resonant Analysis}}}
For $(\iota_1,\iota_2,\iota_3)=(-,+,+)$, $(+,-,+)$, $(+,+,-)$, we estimate the case $(-,+,+)$ here, since the others are similar. We decompose it further into
\begin{equation*}
\begin{aligned}
&(1+\abs \xi^{30})\xi\sum_{m,l_1,l_2}\int q_m(s) e^{-\mathrm{i} s \Phi+\mathrm{i} H} \Widehat {f^{-}_{k}}(-\xi-\eta-\sigma,s) \Widehat {f^{+}_{k}}(\xi+\eta,s)
\\
&\quad\times\Widehat {f^{+}_{k}}(\xi+\sigma,s) \varphi^{\bar l}_{l_1}(\eta) \varphi^{\bar l}_{l_2}(\sigma)\, \diff \eta \diff \sigma \diff s,
\end{aligned}
\end{equation*}
where $\bar l=-49m/100$. First we estimate
\begin{equation}\label{3.3.4.3}
\begin{aligned}
&(1+\abs \xi^{30})\xi \sum_{l_1 \vee l_2>\bar l,m} \int q_m(s) e^{-\mathrm{i} s \Phi+\mathrm{i} H} \Widehat {f^{-}_{k}}(-\xi-\eta-\sigma,s) \Widehat {f^{+}_{k}}(\xi+\eta,s)\\
&\quad\times\Widehat {f^{+}_{k}}(\xi+\sigma,s) \varphi^{\bar l}_{l_1}(\eta) \varphi^{\bar l}_{l_2}(\sigma)\, \diff \eta \diff \sigma \diff s.
\end{aligned}
\end{equation}
We may assume $l_2 \geq l_1$ and use the absence of space resonance. Integrating by parts with respect to $\eta$, one has
\[
\abs{\eqref{3.3.4.3}} \lesssim I_1+I_2+I_3,
\]
where
\begin{equation*}
\begin{aligned}
I_1&=\sum_{l_2 \geq (\bar l+1) \vee l_1,m}\bigg|(1+\abs \xi^{30})\xi \int q_m(s) m_1 e^{-\mathrm{i} s \Phi+\mathrm{i} H} \partial_\eta \Widehat {f^{-}_{k}}(-\xi-\eta-\sigma,s) \\
&\quad\times\Widehat {f^{+}_{k}}(\xi+\eta,s)\Widehat {f^{+}_{k}}(\xi+\sigma,s)  \, \diff \eta \diff \sigma \diff s\bigg|\\
I_2&=\sum_{l_2 \geq (\bar l+1) \vee l_1,m}\bigg|(1+\abs \xi^{30})\xi\int q_m(s) m_1 e^{-\mathrm{i} s \Phi+\mathrm{i} H} \Widehat {f^{-}_{k}}(-\xi-\eta-\sigma,s)\\ &\quad\times\partial_\eta \Widehat {f^{+}_{k}}(\xi+\eta,s)\Widehat {f^{+}_{k}}(\xi+\sigma,s)  \,  \diff \eta \diff \sigma \diff s\bigg|\\
I_3&=\sum_{l_2 \geq (\bar l+1) \vee l_1,m}\bigg|(1+\abs \xi^{30})\xi\int q_m(s) \partial_\eta m_1 e^{-\mathrm{i} s \Phi+\mathrm{i} H} \Widehat {f^{-}_{k}}(-\xi-\eta-\sigma,s) \\
&\quad\times\Widehat {f^{+}_{k}}(\xi+\eta,s)\Widehat {f^{+}_{k}}(\xi+\sigma,s)\,  \diff \eta \diff \sigma \diff s\bigg|
\end{aligned}
\end{equation*}
with
\[
m_1=\frac{\varphi^{\bar l}_{l_1}(\eta)\varphi^{\bar l}_{l_2}(\sigma)}{-\mathrm{i} s \partial_\eta \Phi}.
\]

By Lemma \ref{resoasym},
\[
\abs{\partial_\eta \Phi} \gtrsim 2^{l_2} 2^{-3k/2},
\]
hence using Lemma \ref{multicri} one has
\[
\norm{m_1}_{S} \lesssim 2^{-m} 2^{-l_2} 2^{3k/2}.
\]
For $I_1$, it follows from Theorem \ref{multi}, \eqref{weightasp2}, \eqref{Zasp}, \eqref{estlih} that
\begin{align*}
I_1 &\lesssim \sum_{l_2,l_1,m} 2^{30k_++k}2^{-l_2} 2^{3k/2} (\norm{f}_2+\norm{(xf)_k}_2) \big\|f_{B(\xi,2^{l_1+4})}\big\|_2 \norm{Ef_k}_\infty\\
&\lesssim \epsilon_1^3 \sum_{l_2,l_1,m} 2^{32.5k_+}2^{-l_2}2^{p_2m}2^{l_1/2}2^{-m/2}\\
&\lesssim \epsilon_1^3 \sum_{m} 2^{-\bar l/2} 2^{32.5p_1m}2^{p_2m}2^{-m/2}
\lesssim \epsilon_1^3.
\end{align*}
The estimate on $I_2$ is similar. For $I_3$, by Lemma \ref{multicri},
\[
\norm{\partial_\eta m_1}_S \lesssim 2^{-m} 2^{-l_1} 2^{-l_2} 2^{3k/2}.
\]
which, together with Theorem \ref{multi},\eqref{Zasp} and \eqref{estlih}, yields
\begin{align*}
I_3 &\lesssim \sum_{l_2,l_1,m} 2^{30k_++k}2^{-l_1} 2^{-l_2} 2^{3k/2} \big\|f_{B(\xi,2^{l_2+4})}\big\|_2 \big\|f_{B(\xi,2^{l_1+4})}\big\|_2 \norm{Ef_k}_\infty\\
&\lesssim \epsilon_1^3\sum_{l_2,l_1,m} 2^{32.5k_+}2^{-l_1} 2^{-l_2} 2^{l_2/2} 2^{l_1/2} 2^{-m/2}\\
&\lesssim \epsilon_1^3 \sum_m 2^{32.5p_1m}2^{-\bar l} 2^{-m/2}
\lesssim \epsilon_1^3.
\end{align*}

Finally, we estimate the resonant region:
\begin{equation*}
\begin{aligned}
&(1+\abs \xi^{30}) \sum_{k_1,k_2,k_3,m} \int q_m(s)e^{\mathrm{i} H} \int e^{-\mathrm{i} s \Phi} \xi \Widehat {f^{-}_{k_1}}(-\xi-\eta-\sigma,s) \Widehat {f^{+}_{k_2}}(\xi+\eta,s)\\
&\quad\times\Widehat {f^{+}_{k_3}}(\xi+\sigma,s) \varphi_{\bar l}(\eta) \varphi_{\bar l}(\sigma)\, \, \diff \eta \diff \sigma,
\end{aligned}
\end{equation*}
where $\max\{\abs{k_1-k},\abs{k_2-k},\abs{k_3-k}\} \leq 20$. We extract the correction term and aim to estimate
\begin{equation*}
\begin{aligned}
&(1+\abs \xi^{30})\sum_{k_1,k_2,k_3,m} \int q_m(s) e^{\mathrm{i} H}\bigg( \int e^{-\mathrm{i} s \Phi} \xi \Widehat {f^{-}_{k_1}}(-\xi-\eta-\sigma,s) \\
&\quad\times\Widehat {f^{+}_{k_2}}(\xi+\eta,s)\Widehat {f^{+}_{k_3}}(\xi+\sigma,s)
\varphi_{\bar l}(\eta) \varphi_{\bar l}(\sigma)\, \diff \eta \diff \sigma\\
&\quad - \frac {2\pi\xi} {\abs{\Lambda''(\xi)}}\frac{\Widehat{f^-_{k_1}}(-\xi,s)\Widehat{f^+_{k_2}}(\xi,s)\Widehat{f^+_{k_3}}(\xi,s)}{s}  \varphi_{>1}(\abs \xi s^{1/3})\bigg)\, \diff s.
\end{aligned}
\end{equation*}
Since
\[
\abs{\Phi(\xi,\eta,\sigma)-\Lambda''(\xi)\eta\sigma} \lesssim (\abs \eta+\abs \sigma)^3,
\]
\begin{equation*}
\begin{aligned}
&\sum_m (1+\abs \xi^{30})\bigg| \int q_m(s) e^{\mathrm{i} H}(e^{-\mathrm{i} s \Phi}-e^{-\mathrm{i} s \Lambda''(\xi)\eta\sigma}) \xi \Widehat {f^{-}_{k}}(-\xi-\eta-\sigma,s)  \\ 
&\quad\times \Widehat {f^{+}_{k}}(\xi+\eta,s)\Widehat {f^{+}_{k}}(\xi+\sigma,s) \varphi_{\bar l}(\eta) \varphi_{\bar l}(\sigma)\, \diff \eta \diff \sigma \diff s \bigg|\\
&\lesssim \sum_m 2^{2m} 2^{3\bar l}2^{30k_++k} \big\|\Widehat f\big\|_\infty^3 2^{2\bar l}
\lesssim \epsilon_1^3 \sum_m 2^{2m} 2^{5\bar l} 2^{31p_1m}
\lesssim \epsilon_1^3.
\end{aligned}
\end{equation*}
By \eqref{weightasp2},
\begin{equation*}
\begin{aligned}
&\big|\Widehat{f^-_{k_1}}(-\xi-\eta-\sigma)-\Widehat{f^-_{k_1}}(-\xi)\big|
=\bigg|\int_{-\xi}^{-\xi-\eta-\sigma} \partial \Widehat{f_{k_1}} (\zeta)\, \diff \zeta\bigg|\\
&\lesssim 2^{\bar l/2} (\norm{(xf)_{k_1}}_2+\norm{f_{k_1}}_2)
\lesssim \epsilon_1 2^{\bar l/2}2^{p_2m},
\end{aligned}
\end{equation*}
(the terms for $k_2$ and $k_3$ are similar), which leads to
\begin{equation*}
\begin{aligned}
&\sum_{k_1,k_2,k_3,m} (1+\abs \xi^{30})\bigg| \int q_m(s) e^{\mathrm{i} H}e^{-\mathrm{i} s \Lambda''(\xi)\eta\sigma}\xi \big[\Widehat {f^{-}_{k_1}}(-\xi-\eta-\sigma,s) \Widehat {f^{+}_{k_2}}(\xi+\eta,s)  \\
&\quad\times \Widehat {f^{+}_{k_3}}(\xi+\sigma,s) -\Widehat{f^-_{k_1}}(-\xi,s)\Widehat {f^{+}_{k_2}}(\xi,s) \Widehat {f^{+}_{k_3}}(\xi,s)\big] \varphi_{\bar l}(\eta) \varphi_{\bar l}(\sigma)\, \diff \eta \diff \sigma \diff s \bigg|\\
&\lesssim \epsilon_1^3 \sum_m 2^m 2^{30k_++k}2^{\bar l/2}2^{p_2m} 2^{2\bar l}
\lesssim \epsilon_1^3.
\end{aligned}
\end{equation*}
Since
\[
\int e^{-\mathrm{i} s \Lambda''(\xi) \eta \sigma} \varphi_{\bar l}(\eta) \varphi_{\bar l}(\sigma) \, \, \diff \eta \diff \sigma=\frac{2\pi}{s\Lambda''(\xi)}+O(2^{-2\bar l}2^{-2m} 2^{3k}),
\]
it follows that
\begin{equation*}
\begin{aligned}
&\sum_{k_1,k_2,k_3,m} (1+\abs \xi^{30})\bigg| \int q_m(s) e^{\mathrm{i} H}\bigg(\int e^{-\mathrm{i} s \Lambda''(\xi) \eta\sigma} \varphi_{\bar l}(\eta) \varphi_{\bar l}(\sigma)\, \diff \eta\diff \sigma\\
&\qquad -\frac{2\pi}{s\Lambda''(\xi)} \varphi_{>1}(\abs \xi s^{1/3})\bigg)  \xi\Widehat{f^-_{k_1}}(-\xi,s)\Widehat {f^{+}_{k_2}}(\xi,s) \Widehat {f^{+}_{k_3}}(\xi,s)\, \diff s \bigg| \lesssim \epsilon_1^3.
\end{aligned}
\end{equation*}

\section{Weighted Norms}\label{Weighted Norms}
To prove \eqref{weightasp-2}-\eqref{weightasp2-2}, it suffices to show
\begin{align*}
\|\partial \Widehat{f}\|_2 &\lesssim \epsilon_0+\epsilon_1^2 \langle t\rangle^{1/6};\\
\|\xi \partial \Widehat{f}\|_2 &\lesssim \epsilon_0+\epsilon_1^2 \langle t\rangle^{p_2}.
\end{align*}
\subsection{Basic decompositions}
Take a partition of unity $\chi_1$, $\chi_2$, $\chi_3$ such that
\begin{align*}
&\supp \chi_1 \subset \{ |\xi-\eta-\sigma| \gtrsim |\eta|, |\sigma| \},\\
&\supp \chi_2 \subset \{ |\eta| \gtrsim |\xi-\eta-\sigma|, |\sigma| \},\\
&\supp \chi_3 \subset \{ |\sigma| \gtrsim |\xi-\eta-\sigma|, |\eta| \}.
\end{align*}
We have
\begin{equation*}
\begin{aligned}
\Widehat f(\xi,t)&=\Widehat f(\xi,0)+\mathrm{i} \int_0^t \bigg(\int e^{-\mathrm{i} s \Phi} \xi \chi_1\Widehat f(\xi-\eta-\sigma) \Widehat f(\eta) \Widehat f(\sigma) \, \diff \eta \diff \sigma\bigg) \diff s\\
&\quad+\text{similar terms}.
\end{aligned}
\end{equation*}
Then
\[
\partial_\xi\widehat f(\xi,t)=\partial_\xi \Widehat f(\xi,0)+g_1+g_2+g_3+\text{similar terms},
\]
where
\begin{gather*}
g_1=\mathrm{i} \int_0^t e^{-\mathrm{i} s \Phi} \xi \chi_1(\xi,\eta,\sigma) \partial_\xi \Widehat f(\xi-\eta-\sigma) \Widehat f(\eta) \Widehat f(\sigma) \, \diff \eta \diff \sigma \diff s;\\
g_2=\mathrm{i} \int_0^t e^{-\mathrm{i} s \Phi} \partial_\xi [\xi \chi_1(\xi,\eta,\sigma)] \Widehat f(\xi-\eta-\sigma) \Widehat f(\eta) \Widehat f(\sigma) \, \diff \eta \diff \sigma \diff s;\\
g_3=\mathrm{i} \int_0^t e^{-\mathrm{i} s \Phi} (-\mathrm{i} s \Phi_\xi) \xi \chi_1(\xi,\eta,\sigma) \Widehat f(\xi-\eta-\sigma) \Widehat f(\eta) \Widehat f(\sigma) \, \diff \eta \diff \sigma \diff s.
\end{gather*}

For $g_3$, we extract some regions and make use of the absence of time resonances at first. 
Let
\begin{equation*}
\begin{aligned}
& \supp \zeta_2\subset \{\abs{\xi} \simeq \abs{\xi-\eta-\sigma} \simeq \abs \eta \gg \abs \sigma\},\\
& \supp \zeta_3\subset \{\abs{\xi} \simeq \abs{\xi-\eta-\sigma} \simeq \abs \sigma \gg \abs \eta\},\\
&\supp \zeta_4\subset \{\abs \xi \simeq \abs{\xi-\eta-\sigma} \simeq \abs \eta \simeq \abs \sigma\},\\
& \zeta_1\ \text{supports\ in\ other\ regions}.
\end{aligned}
\end{equation*}
Adding $\zeta_i$, $1 \le i \le 4$, into $g_3$, we have
\begin{align*}
g_3&=\mathrm{i} \sum_j \int_0^t e^{-\mathrm{i} s \Phi} (-\mathrm{i} s \Phi_\xi) \xi \chi_1\zeta_j \Widehat f(\xi-\eta-\sigma) \Widehat f(\eta) \Widehat f(\sigma) \, \diff \eta \diff \sigma \diff s\\
&=g_{31}+g_{32}+g_{33}+g_{34}.
\end{align*}

For $g_{32}$, integrating by parts with respect to $s$ leads to
\[
g_{32}=G+I_0+I_1+I_2+I_3,
\]
where
\begin{gather*}
G=\mathrm{i} \int e^{-\mathrm{i} t \Phi} \frac{t \Phi_\xi}\Phi \xi \chi_1\zeta_2 \Widehat f(\xi-\eta-\sigma) \Widehat f(\eta) \Widehat f(\sigma) \, \diff \eta \diff \sigma;\\
I_0=\mathrm{i} \int_0^t e^{-\mathrm{i} s \Phi}  \frac{\Phi_\xi}\Phi \xi \chi_1\zeta_2 \Widehat f(\xi-\eta-\sigma) \Widehat f(\eta) \Widehat f(\sigma) \, \diff \eta \diff \sigma \diff s;\\
I_1=\mathrm{i} \int_0^t e^{-\mathrm{i} s \Phi}  \frac{s\Phi_\xi}\Phi \xi \chi_1\zeta_2 \partial_s \Widehat f(\xi-\eta-\sigma) \Widehat f(\eta) \Widehat f(\sigma) \, \diff \eta \diff \sigma \diff s;\\
I_2=\mathrm{i} \int_0^t e^{-\mathrm{i} s \Phi}  \frac{s\Phi_\xi}\Phi \xi \chi_1\zeta_2 \Widehat f(\xi-\eta-\sigma) \partial_s \Widehat f(\eta) \Widehat f(\sigma) \, \diff \eta \diff \sigma \diff s;\\
I_3=\mathrm{i} \int_0^t e^{-\mathrm{i} s \Phi}  \frac{s\Phi_\xi}\Phi \xi \chi_1\zeta_2 \Widehat f(\xi-\eta-\sigma) \Widehat f(\eta) \partial_s \Widehat f(\sigma) \, \diff \eta \diff \sigma \diff s.
\end{gather*}
Similarly, we can perform integration by parts with respect to $s$ for $g_{33}$.

For $g_{34}$, we split it as 
\[
g_{34}=\sum_{\iota_1,\iota_2,\iota_3}\mathrm{i} \int_0^t e^{-\mathrm{i} s \Phi} (-\mathrm{i} s \Phi_\xi) \xi \chi_1\zeta_4 \Widehat{f^{\iota_1}}(\xi-\eta-\sigma) \Widehat{f^{\iota_2}}(\eta) \Widehat{f^{\iota_3}}(\sigma) \, \diff \eta \diff \sigma \diff s.
\]
For the cases $(\iota_1,\iota_2,\iota_3)\in \{(-,-,+), (-,+,-), (+,-,-), (+,+,+), (-,-,-)\}$, integrating by parts with respect to $s$, one can get
\begin{align*}
&\mathrm{i} \int_0^t e^{-\mathrm{i} s \Phi} (-\mathrm{i} s \Phi_\xi) \xi \chi_1\zeta_4 \Widehat{f^{\iota_1}}(\xi-\eta-\sigma) \Widehat{f^{\iota_2}}(\eta) \Widehat{f^{\iota_3}}(\sigma) \, \diff \eta \diff \sigma \diff s\\
={}&G+I_0+I_1+I_2+I_3,	
\end{align*}
where
\begin{gather*}
G=\mathrm{i} \int e^{-\mathrm{i} t \Phi} \frac{t \Phi_\xi}\Phi \xi \chi_1\zeta_4 \Widehat{f^{\iota_1}}(\xi-\eta-\sigma) \Widehat{f^{\iota_2}}(\eta) \Widehat{f^{\iota_3}}(\sigma) \, \diff \eta \diff \sigma;\\
I_0=\mathrm{i} \int_0^t e^{-\mathrm{i} t \Phi} \frac{\Phi_\xi}\Phi \xi \chi_1\zeta_4 \Widehat{f^{\iota_1}}(\xi-\eta-\sigma) \Widehat{f^{\iota_2}}(\eta) \Widehat{f^{\iota_3}}(\sigma) \, \diff \eta \diff \sigma \diff s;\\
I_1=\mathrm{i} \int_0^t e^{-\mathrm{i} s \Phi} \frac{s \Phi_\xi}\Phi \xi \chi_1\zeta_4 \partial_s \Widehat{f^{\iota_1}}(\xi-\eta-\sigma) \Widehat{f^{\iota_2}}(\eta) \Widehat{f^{\iota_3}}(\sigma) \, \diff \eta \diff \sigma \diff s;\\
I_2=\mathrm{i} \int_0^t e^{-\mathrm{i} s \Phi} \frac{s \Phi_\xi}\Phi \xi \chi_1\zeta_4 \Widehat{f^{\iota_1}}(\xi-\eta-\sigma) \partial_s \Widehat{f^{\iota_2}}(\eta) \Widehat{f^{\iota_3}}(\sigma) \, \diff \eta \diff \sigma \diff s;\\
I_3=\mathrm{i} \int_0^t e^{-\mathrm{i} s \Phi} \frac{s \Phi_\xi}\Phi \xi \chi_1\zeta_4 \Widehat{f^{\iota_1}}(\xi-\eta-\sigma) \Widehat{f^{\iota_2}}(\eta) \partial_s \Widehat{f^{\iota_3}}(\sigma) \, \diff \eta \diff \sigma \diff s.
\end{gather*}
For the cases $(\iota_1,\iota_2,\iota_3)\in \{(-,+,+), (+,-,+), (+,+,-)\}$, the key ingredient is the following observation:
\begin{equation*}
\begin{aligned}
&\quad\Lambda'(\xi)-\Lambda'(\xi-\eta-\sigma)\\
&=\frac{\Lambda'(\xi)-\Lambda'(\xi-\eta-\sigma)}{\Lambda(\abs \xi)-\iota_1 \Lambda(\xi-\eta-\sigma)} \Phi\\
&\quad-\iota_2 \frac{\Lambda'(\xi)-\Lambda'(\xi-\eta-\sigma)}{\Lambda(\abs \xi)-\iota_1 \Lambda(\xi-\eta-\sigma)} \cdot \frac{\iota_2 \Lambda(\eta)-\iota_1\Lambda(\xi-\eta-\sigma)}{\Lambda'(\eta)-\Lambda'(\xi-\eta-\sigma)} \Phi_\eta\\
&\quad-\iota_3 \frac{\Lambda'(\xi)-\Lambda'(\xi-\eta-\sigma)}{\Lambda(\abs \xi)-\iota_1 \Lambda(\xi-\eta-\sigma)} \cdot \frac{\iota_3 \Lambda(\sigma)-\iota_1\Lambda(\xi-\eta-\sigma)}{\Lambda'(\sigma)-\Lambda'(\xi-\eta-\sigma)} \Phi_\sigma,
\end{aligned}
\end{equation*}
as $\iota_1+\iota_2+\iota_3=1$, which is denoted by
\begin{equation}\label{eq:phixi}
\Phi_\xi=m_1 \Phi+m_2 \Phi_\eta+m_3 \Phi_\sigma
\end{equation}
where 
\begin{equation*}
\begin{aligned}
&\norm{m_1}_{S_{k,k,k,k}} \lesssim 2^{-k}, \\
&\norm{m_2}_{S_{k,k,k,k}} \lesssim 1, \quad \norm{\partial_\eta m_2}_{S_{k,k,k,k}} \lesssim 2^{-k},\\
&\norm{m_3}_{S_{k,k,k,k}} \lesssim 1,\quad
\norm{\partial_\sigma m_3}_{S_{k,k,k,k}} \lesssim 2^{-k}.
\end{aligned}
\end{equation*}

For the terms in $g_{34}$ with signs $(-,+,+)$, $(+,-,+)$, $(+,+,-)$, we split $\Phi_\xi$ by \eqref{eq:phixi}. For $m_1\Phi$, integrating by parts with respect to $s$, one has
\begin{align*}
&\mathrm{i} \int_0^t e^{-\mathrm{i} s \Phi} (-\mathrm{i} s m_1\Phi) \xi \chi_1\zeta_4 \Widehat{f^{\iota_1}}(\xi-\eta-\sigma) \Widehat{f^{\iota_2}}(\eta) \Widehat{f^{\iota_3}}(\sigma) \, \diff \eta \diff \sigma \diff s\\
={}&G+I_0+I_1+I_2+I_3,	
\end{align*}
where
\begin{gather*}
G=\mathrm{i} \int e^{-\mathrm{i} t \Phi} tm_1 \xi \chi_1\zeta_4 \Widehat{f^{\iota_1}}(\xi-\eta-\sigma) \Widehat{f^{\iota_2}}(\eta) \Widehat{f^{\iota_3}}(\sigma) \, \diff \eta \diff \sigma;\\
I_0=\mathrm{i} \int_0^t e^{-\mathrm{i} s \Phi} m_1 \xi \chi_1\zeta_4 \Widehat{f^{\iota_1}}(\xi-\eta-\sigma) \Widehat{f^{\iota_2}}(\eta) \Widehat{f^{\iota_3}}(\sigma) \, \diff \eta \diff \sigma \diff s;\\
I_1=\mathrm{i} \int_0^t e^{-\mathrm{i} s \Phi} s m_1 \xi \chi_1\zeta_4 \partial_s \Widehat{f^{\iota_1}}(\xi-\eta-\sigma) \Widehat{f^{\iota_2}}(\eta) \Widehat{f^{\iota_3}}(\sigma) \, \diff \eta \diff \sigma \diff s;\\
I_2=\mathrm{i} \int_0^t e^{-\mathrm{i} s \Phi} s m_1 \xi \chi_1\zeta_4 \Widehat{f^{\iota_1}}(\xi-\eta-\sigma) \partial_s \Widehat{f^{\iota_2}}(\eta) \Widehat{f^{\iota_3}}(\sigma) \, \diff \eta \diff \sigma \diff s;\\
I_3=\mathrm{i} \int_0^t e^{-\mathrm{i} s \Phi} s m_1 \xi \chi_1\zeta_4 \Widehat{f^{\iota_1}}(\xi-\eta-\sigma) \Widehat{f^{\iota_2}}(\eta) \partial_s \Widehat{f^{\iota_3}}(\sigma) \, \diff \eta \diff \sigma \diff s.
\end{gather*}
The terms $G$ are favorable since there is no integration with respect to time. We extract all good terms $G$ in such a way and write
\[
\partial_\xi \Widehat f(t)=G+\partial_\xi \Widehat f(0)+g_1+g_2+g_{31}+g'_{32}+g'_{33}+g'_{34}+\text{similar terms}=G+H
\]
(here $g'_{32}$, $g'_{33}$, $g'_{34}$ denote the remaining parts after extracting the good terms).

Next we show that
\begin{subequations}\label{estg}
\begin{align}
\norm{G}_2 &\lesssim \epsilon_1^2 \langle t\rangle^{1/6}; \label{estg1}\\
\norm{\xi^\iota G}_2 &\lesssim \epsilon_1^2 \langle t\rangle^{p_2}, \label{estg2}
\end{align}
\end{subequations}
for $1 \leq \iota \leq 5$.
The proof is relatively easy. Note that $G$ can be written as generic forms:
\begin{subequations}
\begin{equation}\label{eq:gterm1}
\begin{aligned}
\int e^{-\mathrm{i} t \Phi} \frac{t \Phi_\xi}{\Phi} \xi \chi_1\zeta_2 \Widehat f(\xi-\eta-\sigma) \Widehat f(\eta) \Widehat f(\sigma) \, \diff \eta \diff \sigma
\end{aligned}
\end{equation}
and
\begin{equation}\label{eq:gterm2}
\begin{aligned}
\int e^{-\mathrm{i} t \Phi} tm\zeta_4 \Widehat{f^{\iota_1}}(\xi-\eta-\sigma) \Widehat{f^{\iota_2}}(\eta) \Widehat{f^{\iota_3}}(\sigma) \, \diff \eta \diff \sigma,
\end{aligned}
\end{equation}
\end{subequations}
where $\norm{m}_{S_{k,k,k,k}} \lesssim 1$. For \eqref{eq:gterm2}, we decompose it into
\[
\sum_k \int e^{-\mathrm{i} t \Phi} tm\zeta_4 \Widehat{f_k^{\iota_1}}(\xi-\eta-\sigma) \Widehat{f_k^{\iota_2}}(\eta) \Widehat{f_k^{\iota_3}}(\sigma) \, \diff \eta \diff \sigma.
\]
The case $t \lesssim 1$ is clear. So we may assume $t \simeq 2^m \gg 1$. If $k \leq -m/3$, then
\begin{equation*}
\begin{aligned}
\norm{\xi^{\iota}G}_2 &\lesssim \sum_k 2^m 2^{\iota k}\norm{f_k}_2 \norm{f_k}_\infty \norm{f_k}_\infty \\
&\lesssim \sum_k \epsilon_1^3 2^m 2^{\iota k+5k/2}
 \lesssim \epsilon_1^3(2^{m/6}| 1).
\end{aligned}
\end{equation*}
If $k>-m/3$, then
\begin{equation*}
\begin{aligned}
\norm{\xi^{\iota}G}_2 &\lesssim \sum_k 2^m 2^{\iota k} \norm{f_k}_2 \norm{f_k}_\infty \norm{f_k}_\infty \\
&\lesssim \sum_k \epsilon_1^3 2^m 2^{\iota k} (2^{k/2} \wedge 2^{p_0m}2^{-10k_+}) 2^{-m} 2^{-k} \\
&\lesssim \epsilon_1^3 (2^{m/6} | 2^{p_2m}).
\end{aligned}
\end{equation*}
For \eqref{eq:gterm1}, we also decompose it into
\[
\sum_k \int e^{-\mathrm{i} t \Phi} \frac{t \Phi_\xi}\Phi \xi \chi_1\zeta_2 \Widehat{f_k}(\xi-\eta-\sigma) \Widehat{f_k}(\eta) \Widehat{f_{\ll k}}(\sigma) \, \diff \eta \diff \sigma.
\]
The cases $k \leq 100$ are similar to \eqref{eq:gterm2}. For $k>100$, note that
\[
\norm{\Phi^{-1}}_{S_{k,k,k,\ll k}} \lesssim 2^{4k}.
\]
The case $t \lesssim 1$ is easy to deal with. When $t \gg 1$, one has
\begin{equation*}
\begin{aligned}
\norm{\xi^\iota G}_2 &\lesssim \sum_k 2^m 2^{-k/2} 2^{4k} 2^{(1+\iota)k}\norm{f_k}_\infty \norm{f_k}_\infty \norm{f_{\ll k}}_2 \\
&\lesssim \sum_k \epsilon_1^3 2^m 2^{-k/2} 2^{4k} 2^{(1+\iota)k} 2^{-m}2^{-20k} \lesssim \epsilon_1^3.
\end{aligned}
\end{equation*}

It follows from \eqref{weightasp}, \eqref{weightasp}, \eqref{estg} that
\begin{subequations}\label{estH}
\begin{align}
\norm{H}_2 &\lesssim \epsilon_1 \langle t\rangle^{1/6}; \label{estH1}\\
\norm{\xi H}_2 &\lesssim \epsilon_1 \langle t\rangle^{p_2}. \label{estH2}
\end{align}
\end{subequations}
For the remaining parts, we use the method of energy estimate. It suffices to prove that
\begin{align*}
(H,H) &\lesssim \epsilon_0^2+\epsilon_1^4 \langle t\rangle^{1/3};\\
(\xi H, \xi H) &\lesssim \epsilon_0^2+\epsilon_1^4 \langle t\rangle^{2p_2}.
\end{align*}
We have
\begin{equation*}
\begin{aligned}
\frac 12 (H,H)_t=(H_t,H)&=((g_1)_t,H)+((g_2)_t,H)+((g_{31})_t,H)\\
&\quad+((g'_{32})_t,H)+((g'_{34})_t,H)+\text{similar terms}.
\end{aligned}
\end{equation*}
Hence it suffices to prove that
\[
\abs{F^\iota_1}, \, \abs{F^\iota_2}, \, \abs{F^\iota_{31}}, \, \abs{F^\iota_{32}},  \, \abs{F^\iota_{34}} \lesssim \epsilon_1^4 (\langle t\rangle^{1/3}|\langle t\rangle^{2p_2}),\quad \iota=0 ,1,
\]
where
\begin{equation*}
\begin{aligned}
F^\iota_1&=\mathrm{i}\sum_m \int q_m(s) e^{-\mathrm{i} s \Phi} \xi^{1+2\iota} \chi_1 \partial_\xi \Widehat f(\xi-\eta-\sigma) \Widehat f(\eta) \Widehat f(\sigma)H(-\xi)\,\diff \xi  \diff \eta \diff \sigma \diff s;\\
F^\iota_2&=\mathrm{i} \sum_m \int e^{-\mathrm{i} s \Phi} \xi^{2\iota}\partial_\xi [\xi \chi_1] \Widehat f(\xi-\eta-\sigma) \Widehat f(\eta) \Widehat f(\sigma) H(-\xi)\,\diff \xi  \diff \eta \diff \sigma \diff s
\end{aligned}
\end{equation*}
and
\begin{equation*}
\begin{aligned}
F^\iota_{31}&=\mathrm{i} \sum_m \int e^{-\mathrm{i} s \Phi} (-\mathrm{i} s \Phi_\xi) \xi^{1+2\iota} \chi_1\zeta_1 \Widehat f(\xi-\eta-\sigma) \Widehat f(\eta) \Widehat f(\sigma) H(-\xi)\,\diff \xi  \diff \eta \diff \sigma \diff s;\\
F^\iota_{32}&=
\mathrm{i} \sum_m \int e^{-\mathrm{i} s \Phi}  \frac{\Phi_\xi}\Phi \xi^{1+2\iota} \chi_1\zeta_2 \Widehat f(\xi-\eta-\sigma) \Widehat f(\eta) \Widehat f(\sigma) H(-\xi) \, \diff \xi \diff \eta \diff \sigma \diff s\\
&\quad+\mathrm{i} \sum_m \int e^{-\mathrm{i} s \Phi}  \frac{s\Phi_\xi}\Phi \xi^{1+2\iota} \chi_1\zeta_2 \partial_s \Widehat f(\xi-\eta-\sigma) \Widehat f(\eta) \Widehat f(\sigma) H(-\xi) \, \diff \xi \diff \eta \diff \sigma \diff s;\\
F^\iota_{34}&=
\mathrm{i} \sum_m \int e^{-\mathrm{i} s \Phi} \xi^{2\iota} m\zeta_4 \Widehat{f^{\iota_1}}(\xi-\eta-\sigma) \Widehat{f^{\iota_2}}(\eta) \Widehat{f^{\iota_3}}(\sigma) H(-\xi)\, \diff \xi  \diff \eta \diff \sigma \diff s\\
&\quad+\mathrm{i} \sum_m \int e^{-\mathrm{i} s \Phi} \xi^{2\iota}s m \zeta_4 \partial_s \Widehat{f^{\iota_1}}(\xi-\eta-\sigma) \Widehat{f^{\iota_2}}(\eta) \Widehat{f^{\iota_3}}(\sigma) H(-\xi)\, \diff \xi  \diff \eta \diff \sigma \diff s\\
&\quad+\mathrm{i} \sum_m \int e^{-\mathrm{i} s \Phi}  \xi^{2\iota}(-\mathrm{i} s m \Phi_\eta) \zeta_4 \Widehat{f^{\iota_1}}(\xi-\eta-\sigma) \Widehat{f^{\iota_2}}(\eta) \Widehat{f^{\iota_3}}(\sigma) H(-\xi)\, \diff \xi  \diff \eta \diff \sigma \diff s,
\end{aligned}
\end{equation*}
in which the terms in $F^\iota_{32}$ and $F^\iota_{34}$ are generic terms. Later the ``$\I$''s in the front of these terms will be dropped for convenience.

Decompose $F^\iota_1$, $F^\iota_2$, $F^\iota_{31}$, $F^\iota_{32}$, $F^\iota_{33}$, $F^\iota_{34}$ into each time piece by adding $q_m$ into them. For the cases $0 \leq m \leq D$ and $m=L+1$, since there is no summation for time, by dividing frequency regions and utilizing the symmetry of symbols, which will be performed for the main cases $D<m \leq L$, they are easy to be estimated. Hence we mainly deal with the cases $D<m \leq L$.

\subsection{Estimates for $F^\iota_1$}
We decompose $F^\iota_1$ as follows:
\begin{equation*}
\begin{aligned}
 F^\iota_1
&=\sum_m \int  q_m e^{-\mathrm{i} s \Phi} \varphi_{\leq -m/3}(\xi) \xi^{1+2\iota} \chi_1 \partial_\xi \Widehat f(\xi-\eta-\sigma)\\
&\quad\times \Widehat f(\eta) \Widehat f(\sigma)H(-\xi)\,\diff \xi\diff \eta \diff \sigma \diff s\\
&\quad+\sum_m \int q_m  e^{-\mathrm{i} s \Phi}\varphi_{>-m/3}(\xi) \xi^{1+2\iota} \chi_1 \theta_j \partial_\xi \Widehat f(\xi-\eta-\sigma)\\
&\quad\times \Widehat f(\eta) \Widehat f(\sigma) H(-\xi)\,\diff \xi  \diff \eta \diff \sigma \diff s\\
&=L^\iota_{1}+H^\iota_{1}.
\end{aligned}
\end{equation*}

By the Coifman-Meyer multiplier theorem, \eqref{weightasp}, \eqref{weightasp2}, \eqref{estH} and \eqref{estli}, one may get
\begin{equation*}
\begin{aligned}
\abs{L^\iota_{1}} &\lesssim \sum_m 2^m 2^{-m/3} \norm{\partial^\iota (xf)}_2 \norm{\xi^\iota H}_2 \norm{u}_\infty^2 \\
&\lesssim \sum_m \epsilon_1^4 (2^{m/3}|2^{2p_2m}) \lesssim \epsilon_1^4 (2^{L/3}|2^{2p_2L}).
\end{aligned}
\end{equation*}

To handle $H^\iota_{1}$, we choose the following cut-off functions 
\begin{equation*}
\begin{aligned}
& \supp \theta_1\subset \{\abs{\xi-\eta-\sigma} \simeq \abs \eta\},\\
& \supp \theta_2\subset \{\abs{\xi-\eta-\sigma} \simeq \abs \sigma\},\\
&\supp \theta_3\subset \{\abs{\xi-\eta-\sigma} \gg \abs \eta, \abs \sigma\},
\end{aligned}
\end{equation*}
and decompose 
\begin{align*}
H^\iota_{1}=\sum_{j=1}^3H^\iota_{1j}
\end{align*}
with
\begin{align*}
H^\iota_{1j}&=\sum_m \int q_m  e^{-\mathrm{i} s \Phi}\varphi_{>-m/3}(\xi) \xi^{1+2\iota} \chi_1 \theta_j \partial_\xi \Widehat f(\xi-\eta-\sigma)\\
&\quad\times \Widehat f(\eta) \Widehat f(\sigma) H(-\xi)\,\diff \xi  \diff \eta \diff \sigma \diff s,\quad j=1,2,3.
\end{align*}

For $H^\iota_{11}$, we further decompose
\begin{equation*}
\begin{aligned}
H^\iota_{11}&=\sum_{m,k_1 \gtrsim k \gtrsim -m/3} \int q_m  e^{-\mathrm{i} s \Phi} \xi^{1+2\iota} \chi_1 \theta_1 \Widehat{(xf)_{k_1}}(\xi-\eta-\sigma) \\
&\quad\times \Widehat{f_{k_1}}(\eta) \Widehat{f_{\lesssim k_1}}(\sigma) (\psi_k H)(-\xi)\,\diff \xi  \diff \eta \diff \sigma \diff s.
\end{aligned}
\end{equation*}
It follows from Lemma \ref{itmulti}, \eqref{decayinter}, \eqref{estli}, \eqref{weightasp}, \eqref{weightasp2} and \eqref{estH} that
\begin{equation*}
\begin{aligned}
\abs{H^\iota_{11}} &\lesssim \sum_{m,k_1,k} 2^m2^{(1+2\iota)k} (\norm{\psi_k H}_2\norm{(xf)_{k_1}}_2 I_{2^{m/3}}(u_{k_1},u_{\lesssim k_1})\\
&\quad+\langle 2^{k_1}2^{m/3}\rangle^{-N} \norm{\psi_k H}_2\norm{(xf)_{k_1}}_2 \norm{u}_\infty^2)\\
&\lesssim \sum_{m,k_1,k} (\epsilon_1^2 2^m 2^{(1+\iota)k} \norm{\xi^\iota \psi_k H}_2\norm{\partial^\iota (xf)_{k_1}}_2 2^{-m} 2^{-(1+\iota)k_1}\\
&\quad+\epsilon_1^4 2^m 2^{(1+2\iota)k} (2^{k_1}2^{m/3})^{-N} 2^{m/3} 2^{-2m/3})\\
&\lesssim \sum_{m,k_1,k} \epsilon_1^2 2^{(1+\iota)(k-k_1)} \norm{\xi^\iota \psi_k H}_2\norm{\partial^\iota (xf)_{k_1}}_2+\epsilon_1^4 (2^{L/3}|1)\\
&\lesssim \epsilon_1^4 (2^{L/3}|2^{2p_2L}). 
\end{aligned}
\end{equation*}
$H^\iota_{12}$ can be handled similarly as $H^\iota_{11}$. 

For $H^\iota_{13}$, by using $H=\partial_\xi \Widehat f-G$, we write it as 
\begin{equation*}
\begin{aligned}
H^\iota_{13}=H^\iota_{131}+H^\iota_{132}+H^\iota_{133},
\end{aligned}
\end{equation*}
where 
\begin{equation*}
\begin{aligned}
H^\iota_{131}&=\sum_m \int q_m m_0(\xi,\eta,\sigma) \Widehat{E(xf)}(\xi-\eta-\sigma) \Widehat u(\eta) \Widehat u (\sigma) \Widehat{E(xf)}(-\xi)\, \diff \xi  \diff \eta \diff \sigma \diff s,\\
H^\iota_{132}&=\sum_m \int q_m  e^{-\mathrm{i} s \Phi} \varphi_{>-m/3}(\xi) \xi^{1+2\iota} \chi_1 \theta_3 \partial_\xi \Widehat f(\xi-\eta-\sigma) \Widehat f(\eta) \Widehat f(\sigma)\\
&\quad \times s(m_1\zeta_4)(-\xi,\mu,\nu) \Widehat{f^{\iota_1}}(-\xi-\mu-\nu) \Widehat{f^{\iota_2}}(\mu) \Widehat{f^{\iota_3}}(\nu)\, \diff \xi  \diff \eta \diff \sigma \diff \mu \diff \nu \diff s,\\
H^\iota_{133}&=\sum_m \int q_m  e^{-\mathrm{i} s \Phi} \varphi_{>-m/3}(\xi) \xi^{1+2\iota} \chi_1 \theta_3 \partial_\xi \Widehat f(\xi-\eta-\sigma) \Widehat f(\eta) \Widehat f(\sigma)\\
&\quad \times s(m_2\zeta_2)(-\xi,\mu,\nu) \Widehat f(-\xi-\mu-\nu) \Widehat f(\mu) \Widehat f(\nu)\, \diff \xi  \diff \eta \diff \sigma \diff \mu \diff \nu \diff s
\end{aligned}
\end{equation*}
with
\begin{equation*}
\begin{aligned}
 m_0=\varphi_{>-m/3} \xi^{1+2\iota} \chi_1 \theta_3,\quad
 \norm{m_1}_{S_{k,k,k,k}} \lesssim 1,\quad
 m_2=\frac{\Phi_\xi\xi}\Phi \chi_1.
\end{aligned}
\end{equation*}
Here $H^\iota_{132}$ and $H^\iota_{133}$ denote generic terms.
In the following, we handle $H^\iota_{131}, H^\iota_{132}, H^\iota_{133}$ term by term. 

First, we use the symmetry of $m_0$ to write
\begin{equation*}
\begin{aligned}
H^\iota_{131}
&=
\frac {1}{2} \sum_m \bigg( \int q_m m_0(\xi,\eta,\sigma) \Widehat{E(xf)}(\xi-\eta-\sigma) \Widehat u(\eta) \Widehat u (\sigma)  \\
&\quad\times\Widehat{E(xf)}(-\xi)\, \diff \xi  \diff \eta \diff \sigma \diff s+\int q_m m_0(-\xi+\eta+\sigma,\eta,\sigma)\\
&\quad\times \Widehat{E(xf)}(-\xi) \Widehat u(\eta) \Widehat u (\sigma)\Widehat{E(xf)}(\xi-\eta-\sigma)\, \diff \xi  \diff \eta \diff \sigma \diff s \bigg)\\
&=\frac 12 \sum_m \int q_m [m_0(\xi,\eta,\sigma)-m_0(\xi-\eta-\sigma,\eta,\sigma)]
\\
&\quad\times \Widehat{E(xf)}(\xi-\eta-\sigma)\Widehat u(\eta) \Widehat u (\sigma) \Widehat{E(xf)}(-\xi)\, \diff \xi  \diff \eta \diff \sigma \diff s\\
&=
\frac {1}{2}\sum_{m,k} \int q_m [m_0(\xi,\eta,\sigma)-m_0(\xi-\eta-\sigma,\eta,\sigma)]  \\
&\quad\times \Widehat{[E(xf)]_k}(\xi-\eta-\sigma)\Widehat{u_{\ll k}}(\eta) \Widehat{u_{\ll k}} (\sigma) \Widehat{[E(xf)]_k}(-\xi) \, \diff \xi \diff \eta \diff \sigma \diff s.
\end{aligned}
\end{equation*}
Since
\[
m_0(\xi,\eta,\sigma)-m_0(\xi-\eta-\sigma,\eta,\sigma)=(\eta+\sigma) \tilde m,
\]
where $\norm{\tilde m}_{S_{k,\ll k,\ll k,k}} \lesssim 2^{2\iota k}$,
it follows from Lemma \ref{itmulti}, \eqref{decayinter}, \eqref{estli}, \eqref{weightasp} and \eqref{weightasp2} that
\begin{equation*}
\begin{aligned}
\abs{H^\iota_{131}} &\lesssim \sum_{m,k} (2^m \norm{\partial^\iota (xf)_k}_2^2 I_{2^{m/3}}((Du)_{\ll k},u_{\ll k})\\
&\quad+ 2^m2^{(1+2\iota)k} \langle 2^k2^{m/3}\rangle^{-N} \norm{(xf)_k}_2^2 \norm{u}_\infty^2) \\
&\lesssim \epsilon_1^4 (2^{L/3}|2^{2p_2L}).
\end{aligned}
\end{equation*}

Next, we decompose 
\begin{equation*}
\begin{aligned}
H^\iota_{132}&=\mathrm{i} \sum_{m,k} \int q_m e^{-\mathrm{i} s \Phi} \xi^{1+2\iota} s\chi_1 \theta_3m_1\zeta_4 \Widehat{(xf)_k}(\xi-\eta-\sigma) \Widehat{f_{\ll k}}(\eta) \Widehat{f_{\ll k}}(\sigma)\\
&\quad\times \Widehat{f_k^{\iota_1}}(-\xi-\mu-\nu) \Widehat{f_k^{\iota_2}}(\mu) \Widehat{f_k^{\iota_3}}(\nu)\, \diff \xi  \diff \eta \diff \sigma \diff \mu \diff \nu \diff s,
\end{aligned}
\end{equation*}
and then estimate
\begin{equation*}
\begin{aligned}
\abs{H^\iota_{132}} &\lesssim \sum_{m,k} 2^m 2^{(1+2\iota)k} 2^m (\norm{(xf)_k}_2\norm{f_k}_22^{-2k}I^2_{2^{m/3}}((Du)_k,u_{\ll k})\\
&\quad+\langle 2^{m/3} 2^k\rangle^{-N} \norm{(xf)_k}_2 \norm{f_k}_2\norm{u}^4_\infty)\\
&\lesssim
 \sum_{m,k}(\epsilon_1^5 2^m 2^{(-1+\iota)k} 2^m \norm{\partial^\iota (xf)_k}_2 (2^{k/2} \wedge 2^{p_0m}2^{-10k_+}) 2^{-2m}\\
&\quad + \epsilon_1^6 2^m 2^{(1+2\iota)k} 2^m (2^{m/3}2^k)^{-N} 2^{m/6}2^{k/2}2^{-4m/3})\\
&\lesssim \epsilon_1^4 (2^{L/3} | 2^{2p_2L}),
\end{aligned}
\end{equation*}
where one has used Lemma \ref{itmulti}, \eqref{weightasp}, \eqref{weightasp2}, \eqref{estl2l}, \eqref{energyasp}, \eqref{decayinter} and \eqref{estli}.

Finally, we handle $H^\iota_{133}$. When $k \leq 100$, $H^\iota_{133}$ can  be estimated  in a similar way as  $H^\iota_{132}$. For $k>100$, we decompose $H^\iota_{133}$ as follows:
\begin{equation*}
\begin{aligned}
H^\iota_{133}=&\sum_{m,k} \int q_m e^{-\mathrm{i} s \Phi} \xi^{1+2\iota} s\chi_1 \theta_3m_2\zeta_2 \Widehat{(xf)_k}(\xi-\eta-\sigma) \Widehat{f_{\ll k}}(\eta) \Widehat{f_{\ll k}}(\sigma)\\
&\quad\times \Widehat{f_k}(-\xi-\mu-\nu) \Widehat{f_k}(\mu) \Widehat{f_{\ll k}}(\nu)\, \diff \xi  \diff \eta \diff \sigma \diff \mu \diff \nu \diff s,
\end{aligned}
\end{equation*}
and then use
Lemma \ref{itmulti}, \eqref{weightasp}, \eqref{weightasp2}, \eqref{estl2}, \eqref{decayinter}, \eqref{energyasp} and \eqref{estli} to estimate
\begin{equation*}
\begin{aligned}
\abs{H^\iota_{133}} &\lesssim \sum_{m,k} 2^m 2^{(1+2\iota)k} 2^m 2^{10k}(\norm{(xf)_k}_2\norm{f_{\ll k}}_2 I^2_{2^{m/3}}(u_k,u_{\ll k})\\
&\quad+\langle 2^k2^{m/3}\rangle^{-3} \norm{(xf)_k}_2 \norm{f_k}_2\norm{u}^4_\infty)\\
&\lesssim \sum_{m,k} \epsilon_1^5 2^{2m}2^{(1+\iota)k} 2^{10k}\norm{\partial^\iota (xf)_k}_2 2^{-20k}2^{-2m}\\
&\quad+\sum_{m,k} \epsilon_1^6 2^{2m} 2^{(1+2\iota)k} 2^{-3k} 2^{-m} 2^{m/6} 2^{p_0m}2^{-N_0k} 2^{-4m/3}\\
&\lesssim \epsilon_1^4 (2^{L/3}|2^{2p_2L}). 
\end{aligned}
\end{equation*}

\subsection{Estimates for $F^\iota_2$}

For $F^0_2$, by the Coifman-Meyer multiplier theorem, \eqref{decayl6} and \eqref{weightasp}, it holds that
\[
\abs{F^0_2} \lesssim \sum_m 2^m \norm{u}^3_6\norm{H}_2 \lesssim \epsilon_1^4 2^{L/3}.
\]

For $F^1_{2}$, we use the following cut-off functions 
\begin{equation*}
\begin{aligned}
 \supp \theta_1\subset \{\abs \eta \gtrsim \abs \sigma\},\quad
 \supp \theta_2\subset \{\abs \sigma \gtrsim \abs \eta\},
\end{aligned}
\end{equation*}
and then further decompose it as 
\begin{equation*}
\begin{aligned}
F^1_{2}=F^1_{21}+F^1_{22},
\end{aligned}
\end{equation*}
where
\begin{equation*}
\begin{aligned}
F^1_{2j}
&=\sum_m \int q_m e^{-\mathrm{i} s \Phi} \xi^2 \partial_\xi [\xi \chi_1] \theta_j \Widehat{f}(\xi-\eta-\sigma) \Widehat f(\eta) \Widehat f(\sigma) H(-\xi)\, \diff \xi  \diff \eta \diff \sigma \diff s\\
&=\sum_{m,k,k_1 \gtrsim k,k_2 \lesssim k_1} \int q_m e^{-\mathrm{i} s \Phi} \psi_k(\xi) \xi^2 \partial_\xi [\xi \chi_1] \theta_j \Widehat{f_{k_1}}(\xi-\eta-\sigma) \\
&\quad\times \Widehat{f_{k_2}}(\eta) \Widehat{f_{\lesssim k_2}}(\sigma) H(-\xi) \, \diff \xi \diff \eta \diff \sigma \diff s,\quad j=1,2.
\end{aligned}
\end{equation*}
Since $F^1_{22}$ and $F^1_{21}$ can be handled in a similar way, we only  deal with  $F^1_{21}$. 
According to the Coifman-Meyer multiplier theorem, \eqref{estlih}, \eqref{estlil}, \eqref{estlilow}, \eqref{estl2l} and \eqref{estH2},
\begin{equation*}
\begin{aligned}
\abs{F^1_{21}} &\lesssim \sum_{m,k,k_1,k_2} 2^m 2^{2k} \norm{f_{k_1}}_\infty \|f_{k_2}\|_\infty \norm{f_{\lesssim k_2}}_2 \norm{\psi_k H}_2\\
&\lesssim \sum_{m,k,k_1,k_2} \epsilon_1^42^m 2^k 2^{-m/2}2^{-2(k_1)_+}2^{-\widetilde{k_1}/2} (2^{-m/2} 2^{-k_2/2} \wedge 2^{k_2}) 2^{k_2/2} 2^{p_2m} \\
&\lesssim
\sum_{m,k,k_1,k_2 \leq -m/3} \epsilon_1^4 2^{p_2m}2^m 2^k 2^{-m/2}2^{-2(k_1)_+}2^{-\widetilde{k_1}/2} 2^{k_2} 2^{k_2/2}\\
&\quad+\sum_{m,k,k_1,k_2 \geq -m/3} \epsilon_1^4 2^{p_2m}2^m 2^k 2^{-m/2}2^{-2(k_1)_+}2^{-\widetilde{k_1}/2} 2^{-m/2} 2^{-k_2/2} 2^{k_2/2}
\\
&\lesssim \sum_{m,k,k_1} \epsilon_1^4 2^{p_2m}2^k 2^{-2(k_1)_+}2^{-\widetilde{k_1}/2}+\sum_{m,k,k_1,k_2 \geq -m/3} \epsilon_1^4 2^{p_2m} 2^k 2^{-2(k_1)_+}2^{-\widetilde{k_1}/2}\\
&\lesssim \epsilon_1^4 2^{2p_2L}.
\end{aligned}
\end{equation*}

\subsection{Estimates for $F^\iota_{31}$}
First, notice that $F^\iota_{31}$ can be shown to have the desired bound when $\abs{\xi-\eta-\sigma} \geq 2^{p_1m}$. Second, if $\abs{\xi-\eta-\sigma} \leq 2^{-m/3}$, then the low frequency part $F^\iota_{31l}$ satisfies that
\begin{equation*}
\begin{aligned}
\norm{F^\iota_{31l}}_2 \lesssim \sum_m \epsilon_1^4 2^m 2^m 2^{-2m/3} 2^{-(1+2\iota)m/3} 2^{-m/6} 2^{-2m/3} 2^{m/6}\lesssim \epsilon_1^4 (2^{L/3}|1).
\end{aligned}
\end{equation*}
Therefore, it suffices to consider the case $2^{-m/3}<\abs{\xi-\eta-\sigma}< 2^{p_1m}$ for $F^\iota_{31}$. $F^\iota_{32}$, $F^\iota_{33}$, $F^\iota_{34}$ can be handled similarly.

Choose the cut-off functions $\theta_1$, $\theta_2$, $\theta_3$, $\theta_4$ such that 
\begin{align*}
&\supp \theta_1 \subset \{ \abs{\xi-\eta-\sigma} \simeq \abs \eta \gg \abs \sigma \}, \\
&\supp \theta_2 \subset \{ \abs{\xi-\eta-\sigma} \simeq \abs \sigma \gg \abs \eta \}, \\
&\supp \theta_3 \subset \{ \abs{\xi-\eta-\sigma} \gg \abs \eta, \abs \sigma \}, \\
&\supp \theta_4 \subset \{ \abs{\xi-\eta-\sigma} \simeq \abs \eta \simeq \abs \sigma \gg \abs \xi\}.
\end{align*}
Then we may decompose $F^\iota_{31}$ as follows:
\begin{equation*}
\begin{aligned}
F^\iota_{31}=\sum_{j=1}^4F^\iota_{31j}
\end{aligned}
\end{equation*}
with
\begin{equation*}
\begin{aligned}
F^\iota_{31j}&=\sum_m \int q_m e^{-\mathrm{i} s \Phi} (-\mathrm{i} s \Phi_\xi) \varphi_{(-m/3,p_1m)}(\xi-\eta-\sigma) \xi^{1+2\iota} \chi_1 \theta_j \zeta_1\\
&\quad \times \Widehat f(\xi-\eta-\sigma) \Widehat f(\eta) \Widehat f(\sigma)H(-\xi) \, \diff \xi \diff \eta \diff \sigma \diff s,\quad j=1,2,3,4. 
\end{aligned}
\end{equation*}

\subsubsection{Estimates for $F^\iota_{311}$ and  $F^\iota_{312}$} 
Since $F^\iota_{312}$ can be dealt with in a similar way to $F^\iota_{311}$, we only focus on $F^\iota_{311}$. 
We divide $F^\iota_{311}$ into two parts as follows:
\begin{equation*}
\begin{aligned}
F^\iota_{311}&=\sum_m\int q_m e^{-\mathrm{i} s \Phi} (-\mathrm{i} s \Phi_\xi) \xi^{1+2\iota} \varphi_{(-m/3,100)}(\xi-\eta-\sigma) \chi_1 \theta_1 \zeta_1 \\
&\quad\times \Widehat{f}(\xi-\eta-\sigma) \Widehat{f}(\eta) \Widehat{f}(\sigma) H(-\xi) \, \diff \xi \diff \eta \diff \sigma \diff s\\
&\quad+\sum_m\int q_m e^{-\mathrm{i} s \Phi} (-\mathrm{i} s \Phi_\xi) \xi^{1+2\iota} \varphi_{[100,p_1m)}(\xi-\eta-\sigma) \chi_1 \theta_1 \zeta_1 \\
&\quad \times\hat{f}(\xi-\eta-\sigma) \Widehat{f}(\eta) \Widehat{f}(\sigma) H(-\xi) \, \diff \xi \diff \eta \diff \sigma \diff s\\
&=L^\iota_{311}+H^\iota_{311}.
\end{aligned}
\end{equation*}\\

\noindent\boxed{\emph{Analysis of $L^\iota_{311}$}.}
To handle  $L^\iota_{311}$, we choose the cut-off functions $\rho_1$ and $\rho_2$ such that 
\begin{align*}
\supp \rho_1 \subset \{\abs \sigma \gtrsim \abs \xi \}, \quad
\supp \rho_2 \subset \{\abs \sigma \ll \abs \xi \},
\end{align*}
and then decompose 
\begin{equation*}
\begin{aligned}
L^\iota_{311}=L^\iota_{3111}+L^\iota_{3112},
\end{aligned}
\end{equation*}
where
\begin{equation*}
\begin{aligned}
L^\iota_{311j}&=\sum_m \int q_m e^{-\mathrm{i} s \Phi} (-\mathrm{i} s \Phi_\xi) \xi^{1+2\iota} \varphi_{(-m/3,100)}(\xi-\eta-\sigma) \chi_1 \theta_1 \zeta_1 \rho_j\\
&\quad\times \Widehat{f}(\xi-\eta-\sigma) \Widehat{f}(\eta) \Widehat{f}(\sigma) H(-\xi) \, \diff \xi \diff \eta \diff \sigma \diff s,\quad j=1,2. 
\end{aligned}
\end{equation*}\\

\noindent\boxed{\boxed{\emph{Analysis of $L^\iota_{3111}$}.}}
We decompose $L^\iota_{3111}$ by frequencies 
\begin{equation*}
\begin{aligned}
L^\iota_{3111}&=\sum_{m,k,k_1,k_2} \int q_m e^{-\mathrm{i} s \Phi} (-\mathrm{i} s \Phi_\xi) \xi^{1+2\iota} \chi_1 \theta_1 \zeta_1 \rho_1 \Widehat{f_{k_1}}(\xi-\eta-\sigma) \\
&\quad\times \Widehat{f_{k_1}}(\eta) \Widehat{f_{k_2}}(\sigma) H_k(-\xi) \, \diff \xi \diff \eta \diff \sigma \diff s\\
&=\sum_{m,k,k_1,k_2} I^\iota_{m,k,k_1,k_2}.
\end{aligned}
\end{equation*}
Here $H_k$ denotes $\psi_k H$.
Integrating by parts with respect to $\sigma$, one has
\[
I^\iota_{m,k,k_1,k_2}=I_1+I_2+I_3,
\]
where
\begin{gather*}
I_1=\int q_m e^{-\mathrm{i} s \Phi} m_0 \partial_\sigma \Widehat{f_{k_1}}(\xi-\eta-\sigma) \Widehat{f_{k_1}}(\eta) \Widehat{f_{k_2}}(\sigma) H_k(-\xi) \, \diff \xi \diff \eta \diff \sigma \diff s;\\
I_2=\int q_m e^{-\mathrm{i} s \Phi} m_0 \Widehat{f_{k_1}}(\xi-\eta-\sigma) \Widehat{f_{k_1}}(\eta) \partial_\sigma \Widehat{f_{k_2}}(\sigma)  H_k(-\xi) \, \diff \xi \diff \eta \diff \sigma \diff s;\\
I_3=\int q_m e^{-\mathrm{i} s \Phi} \partial_\sigma m_0 \Widehat{f_{k_1}}(\xi-\eta-\sigma) \Widehat{f_{k_1}}(\eta) \Widehat{f_{k_2}}(\sigma)  H_k(-\xi) \, \diff \xi \diff \eta \diff \sigma \diff s,
\end{gather*}
and
\[
m_0=\frac{-\xi^{1+2\iota} \Phi_\xi \chi_1 \theta_1 \zeta_1 \rho_1}{\Phi_\sigma}.
\]
Since
\[
\norm{m_0}_{S_{k,k_1,k_1,k_2}} \lesssim 2^{(1+2\iota)k},
\]
it follows from Theorem \ref{multi}, \eqref{weightasp}, \eqref{weightasp2}, \eqref{estl2l}, \eqref{estlil}, \eqref{estH} that
\begin{equation*}
\begin{aligned}
\abs{I_1} &\lesssim 2^m 2^{(1+2\iota)k} (\norm{(xf)_{k_1}}_2+\epsilon_12^{-k_1/2}) \norm{Ef_{k_1}}_\infty \norm{Ef_{k_2}}_\infty \norm{H_k}_2\\
&\lesssim
\epsilon_1^2 2^m 2^k \norm{\partial^\iota (xf)_{k_1}}_2 \norm{\xi^\iota H_k}_2 2^{-m/2}2^{-k_1/2} 2^{-m/2} 2^{-k_2/2}\\
&\quad+\epsilon_1^4 2^{(1+2\iota)k} 2^{-k_1/2} 2^{-k_1/2} 2^{-k_2/2} (2^{m/6} \wedge 2^{p_2m}2^{-k}).
\end{aligned}
\end{equation*}
By Young's inequality, $I_1$ can be estimated.

Similarly, $I_2$ can be estimated. For $I_3$, due to
\[
\norm{\partial_\sigma m_0}_{S_{k.k_1.k_1,k_2}} \lesssim 2^{(1+2\iota)k}2^{-k_1},
\]
$I_3$ can be estimated.\\

\noindent\boxed{\boxed{\emph{Analysis of $L^\iota_{3112}$}.}}
First, we deal with $L^0_{3112}$.
Integrating by parts with respect to $\sigma$ yields
\begin{equation*}
\begin{aligned}
L^0_{3112}=I_1+I_2+I_3,
\end{aligned}
\end{equation*}
where
\begin{gather*}
I_1=\sum_m\int q_m e^{-\mathrm{i} s \Phi} m_0 \partial_\sigma \Widehat{f}(\xi-\eta-\sigma) \Widehat{f}(\eta) \Widehat{f}(\sigma) H(-\xi) \, \diff \xi \diff \eta \diff \sigma \diff s;\\
I_2=\sum_m\int q_m e^{-\mathrm{i} s \Phi} m_0 \Widehat{f}(\xi-\eta-\sigma)  \Widehat{f}(\eta) \partial_\sigma \Widehat{f}(\sigma) H(-\xi) \, \diff \xi \diff \eta \diff \sigma \diff s;\\
I_3=\sum_m\int q_m e^{-\mathrm{i} s \Phi} \partial_\sigma m_0 \Widehat{f}(\xi-\eta-\sigma) \Widehat{f}(\eta) \Widehat{f}(\sigma) H(-\xi) \, \diff \xi \diff \eta \diff \sigma \diff s,
\end{gather*}
with
\begin{equation*}
\begin{aligned}
m_0=\frac{-\Phi_\xi \xi  \chi_1\theta_1\zeta_1 \rho_2 \varphi_{(-m/3,100)}}{\Phi_\sigma}.
\end{aligned}
\end{equation*}

By the Coifman-Meyer multiplier theorem, \eqref{decayl}, \eqref{decayh}, \eqref{decayl6}, \eqref{weightasp} and \eqref{estH1}, one may have the following estimates:
\begin{equation*}
\begin{aligned}
\abs{I_2} \lesssim \sum_m 2^m \big\|\abs{\partial}^{1/2}u\big\|^2_\infty \norm{xf}_2 \norm{H}_2 \lesssim \sum_m \epsilon_1^4 2^{m/3} \lesssim \epsilon_1^4 2^{L/3},
\end{aligned}
\end{equation*}
and
\begin{equation*}
\begin{aligned}
\abs{I_3} \lesssim \sum_m 2^m \|u\|^3_6 \norm{H}_2 \lesssim \sum_m \epsilon_1^4 2^{m/3} \lesssim \epsilon_1^4 2^{L/3}.
\end{aligned}
\end{equation*}

To handle $I_1$, we further decompose it in frequencies 
\begin{equation*}
\begin{aligned}
I_1=\sum_m\int q_m e^{-\mathrm{i} s \Phi} m_0 \Widehat{(xf)_{k_1}}(\xi-\eta-\sigma) \Widehat{f_{k_1}}(\eta) \Widehat{f_{\ll k_1}}(\sigma) H_k(-\xi) \, \diff \xi \diff \eta \diff \sigma \diff s
\end{aligned}
\end{equation*}
and then estimate
\begin{equation*}
\begin{aligned}
\abs{I_1} &\lesssim \sum_{m,k,k_1} 2^m(\norm{(xf)_{k_1}}_2 \norm{\xi H_k}_2I_{2^{m/3}}(u_{k_1},u_{\ll k_1})\\
&\quad+\langle 2^{k_1}2^{m/3}\rangle^{-N}\norm{(xf)_{k_1}}_2 \norm{\xi H_k}_2 \norm{u}^2_\infty)\\
&\lesssim \sum_{m,k,k_1}( \epsilon_1^22^{k-k_1}\norm{(xf)_{k_1}}_2 \norm{H_k}_2\\
&\quad+ \epsilon_1^42^m(2^{k_1}2^{m/3})^{-N}2^{m/6}2^k2^{m/6}2^{-2m/3})\\
&\lesssim \epsilon_1^4 2^{L/3},
\end{aligned}
\end{equation*}
where Lemma \ref{itmulti}, \eqref{weightasp}, \eqref{estH1}, \eqref{decayinter} and \eqref{estli} have been used.

Next, we deal with $L^1_{3112}$
and decompose it into
\begin{equation*}
\begin{aligned}
L^1_{3112}&=\sum_m\int q_m e^{-\mathrm{i} s \Phi} (-\mathrm{i} s \Phi_\xi) \xi^3  \chi_1 \theta_1 \zeta_1 \rho_2 \Widehat{f_{k_1}}(\xi-\eta-\sigma) \Widehat{f_{k_1}}(\eta) \\
&\quad\times \Widehat{f_{\ll k}}(\sigma) H_k(-\xi) \, \diff \xi \diff \eta \diff \sigma \diff s\\
&=\sum_{m,k,k_1} I^\iota_{m,k,k_1}.
\end{aligned}
\end{equation*}
Integrating by parts with respect to $\eta$ gives
\[
I^\iota_{m,k,k_1}=I^\iota_{m,k,k_1;1}+I^\iota_{m,k,k_1;2}+I^\iota_{m,k,k_1;3},
\]
where
\begin{gather*}
I^\iota_{m,k,k_1;1}=\int q_m e^{-\mathrm{i} s \Phi} m_0 \partial_\eta \Widehat{f_{k_1}}(\xi-\eta-\sigma) \Widehat{f_{k_1}}(\eta) \Widehat{f_{\ll k}}(\sigma) H_k(-\xi) \, \diff \xi \diff \eta \diff \sigma \diff s;\\
I^\iota_{m,k,k_1;2}=\int q_m e^{-\mathrm{i} s \Phi} m_0 \Widehat{f_{k_1}}(\xi-\eta-\sigma) \partial_\eta \Widehat{f_{k_1}}(\eta) \Widehat{f_{\ll k}}(\sigma) H_k(-\xi) \, \diff \xi \diff \eta \diff \sigma \diff s;\\
I^\iota_{m,k,k_1;3}=\int q_m e^{-\mathrm{i} s \Phi} \partial_\eta m_0 \Widehat{f_{k_1}}(\xi-\eta-\sigma) \Widehat{f_{k_1}}(\eta) \Widehat{f_{\ll k}}(\sigma) H_k(-\xi) \, \diff \xi \diff \eta \diff \sigma \diff s,
\end{gather*}
with
\begin{equation*}
\begin{aligned}
m_0=\frac{-\Phi_\xi \xi^3 \chi_1 \theta_1 \zeta_1 \rho_2}{\Phi_\eta}.
\end{aligned}
\end{equation*}

Note that $k_1 \gg k$ and
\begin{equation*}
\begin{aligned}
\Phi_\eta=(\xi-\sigma)\tilde m
\end{aligned}
\end{equation*}
where
\begin{equation*}
\begin{aligned}
\norm{\tilde m}_{S_{k,k_1,k_1,\ll k}} \lesssim 2^{k_1}.
\end{aligned}
\end{equation*}
Then by Lemma \ref{itmulti}, \eqref{weightasp}, \eqref{weightasp2}, \eqref{estH}, \eqref{estl2l}, \eqref{decayinter} and \eqref{estli}, one obtains
\begin{equation*}
\begin{aligned}
&\quad\sum_{m,k,k_1} \abs{I^\iota_{m,k,k_1;1}}\\
&\lesssim \sum_{m,k,k_1} 2^m 2^{2k}2^{k_1}(\big\|\partial_\eta \Widehat{f_{k_1}}\big\|_2\norm{H_k}_2I_{2^{m/3}}(u_{k_1},u_{\ll k})\\
&\quad+\langle 2^k2^{m/3}\rangle^{-N}\big\|\partial_\eta \Widehat{f_{k_1}}\big\|_2\norm{H_k}_2\norm{u}^2_\infty)\\
&\lesssim \sum_{m,k,k_1}( \epsilon_1^2 2^{2k} (\norm{(xf)_{k_1}}_2+\epsilon_12^{-k_1/2}) \norm{H_k}_2\\
&\quad+ \epsilon_1^4 2^m 2^{2k} 2^{k_1} \langle 2^k 2^{m/3}\rangle^{-N} 2^{m/3} 2^{-2m/3})\\
&\lesssim \sum_{m,k,k_1} \epsilon_1^2 2^{k-k_1} \norm{\partial (xf)_{k_1}}_2 \norm{\xi H_k}_2+\epsilon_1^4 2^{2p_2L}\\
&\lesssim \epsilon_1^4 2^{2p_2L}.
\end{aligned}
\end{equation*}

$I^\iota_{m,k,k_1;2}$ can be handled similarly as $I^\iota_{m,k,k_1;1}$.
$I^\iota_{m,k,k_1;3}$ is easy to handle. \\

\noindent\boxed{\emph{Analysis of $H^\iota_{311}$}.}
We decompose $H^\iota_{311}$ by frequencies as follows:
\begin{equation*}
\begin{aligned}
H^\iota_{311}&=\sum_{m,k,k_1,k_2} \int q_m e^{-\mathrm{i} s \Phi} (-\mathrm{i} s \Phi_\xi) \xi^{1+2\iota} \chi_1 \theta_1 \zeta_1 \Widehat{f_{k_1}}(\xi-\eta-\sigma) \\
&\quad\times \Widehat{f_{k_1}}(\eta)\Widehat{f_{k_2}}(\sigma) H_k(-\xi) \, \diff \xi \diff \eta \diff \sigma \diff s\\
&=\sum_{m,k,k_1,k_2} I^\iota_{m,k,k_1,k_2}.
\end{aligned}
\end{equation*}
Notice that $I^\iota_{m,k,k_1,k_2}$  is easy to handle when $k \leq -100m$. In the following, we only 
consider the case $k> -100m$.\\

\noindent \underline{1. $k_2 \gtrsim k$.}
Integrating by parts with respect to $\sigma$ shows
\[
I^\iota_{m,k,k_1,k_2}=I^\iota_{m,k,k_1,k_2;1}+I^\iota_{m,k,k_1,k_2;2}+I^\iota_{m,k,k_1,k_2;3},
\]
where
\begin{gather*}
I^\iota_{m,k,k_1,k_2;1}=\int q_m e^{-\mathrm{i} s \Phi} m_0 \partial_\sigma \Widehat{f_{k_1}}(\xi-\eta-\sigma) \Widehat{f_{k_1}}(\eta) \Widehat{f_{k_2}}(\sigma) H_k(-\xi) \, \diff \xi \diff \eta \diff \sigma \diff s;\\
I^\iota_{m,k,k_1,k_2;2}=\int q_m e^{-\mathrm{i} s \Phi} m_0 \Widehat{f_{k_1}}(\xi-\eta-\sigma) \Widehat{f_{k_1}}(\eta) \partial_\sigma \Widehat{f_{k_2}}(\sigma) H_k(-\xi) \, \diff \xi \diff \eta \diff \sigma \diff s;\\
I^\iota_{m,k,k_1,k_2;3}=\int q_m e^{-\mathrm{i} s \Phi} \partial_\sigma m_0 \Widehat{f_{k_1}}(\xi-\eta-\sigma) \Widehat{f_{k_1}}(\eta) \Widehat{f_{k_2}}(\sigma) H_k(-\xi) \, \diff \xi \diff \eta \diff \sigma \diff s,
\end{gather*}
with
\[
m_0=\frac{- \xi^{1+2\iota} \Phi_\xi \chi_1 \theta_1 \zeta_1}{\Phi_\sigma}.
\]

Note that
\[
\norm{m_0}_{S_{k,k_1,k_1,k_2}} \lesssim 2^{(1+2\iota)k} 2^{4k_1}.
\]
Then it follows from Theorem \ref{multi}, \eqref{weightasp2}, \eqref{estl2l}, \eqref{estH}, \eqref{estlih} and \eqref{estlil} that
\begin{equation*}
\begin{aligned}
&\quad\sum_{m,k,k_1,k_2} \abs{I^\iota_{m,k,k_1,k_2;1}}\\
&\lesssim \sum_{m,k,k_1,k_2} 2^m 2^{(1+2\iota)k} 2^{4k_1} \big\|\partial_\sigma \Widehat{f_{k_1}}\big\|_2 \norm{H_k}_2 \norm{u_{k_1}}_\infty \norm{u_{k_2}}_\infty\\
&\lesssim \sum_{m,k,k_1,k_2} \epsilon_1^2 2^{(1+\iota)k} 2^{4k_1}(\norm{(xf)_{k_1}}_2+\epsilon_1 2^{-k_1/2}) \norm{\xi^\iota H_k}_2 2^{-10k_1} 2^{-k_2/2}\\
&\lesssim \sum_{m,k,k_1,k_2}( \epsilon_1^4 2^{(1+\iota)k} 2^{-7k_1} 2^{p_2m} (2^{m/6}|2^{p_2m}) 2^{-k_2/2}\\
&\quad+\epsilon_1^4 2^{(1+\iota)k} 2^{7k_1/2} (2^{m/6}|2^{p_2m}) 2^{-10k_1} 2^{-k_2/2})\\
&\lesssim \epsilon_1^4 (2^{L/3}|2^{2p_2L}),
\end{aligned}
\end{equation*}
while Theorem \ref{multi}, \eqref{weightasp2}, \eqref{estl2l}, \eqref{estH} and \eqref{estlih} imply that
\begin{equation*}
\begin{aligned}
&\quad\sum_{m,k,k_1,k_2} \abs{I^\iota_{m,k,k_1,k_2;2}}\\
&\lesssim \sum_{m,k,k_1,k_2} 2^m 2^{(1+2\iota)k}2^{4k_1}\norm{u_{k_1}}^2_\infty \big\|\partial_\sigma \Widehat{f_{k_2}}\big\|_2 \norm{H_k}_2\\
&\lesssim \sum_{m,k,k_1,k_2} \epsilon_1^2 2^{(1+\iota)k}2^{-10k_1}(\norm{(xf)_{k_2}}_2+\epsilon_12^{-k_2/2}) \norm{\xi^\iota H_k}_2\\
&\lesssim \sum_{m,k,k_1,k_2}( \epsilon_1^4 2^{(1+\iota)k} 2^{-10k_1} 2^{p_2m}2^{-k_2} (2^{m/6}|2^{p_2m})\\
&\quad+ \epsilon_1^4 2^{(1+\iota)k}2^{-10k_1} 2^{-k_2/2} (2^{m/6}|2^{p_2m}) )\\
&\lesssim \epsilon_1^4 (2^{L/3}|2^{2p_2L}).
\end{aligned}
\end{equation*}

$I^\iota_{m,k,k_1,k_2;3}$ can be estimated easily.\\

\noindent \underline{2. $k_2 \ll k$.}
In this case, one may rewrite 
\begin{equation*}
\begin{aligned}
H^\iota_{311}&=\sum_{m,k,k_1} \int q_m e^{-\mathrm{i} s \Phi} (-\mathrm{i} s \Phi_\xi) \xi^{1+2\iota} \chi_1 \theta_1 \zeta_1 \Widehat{f_{k_1}}(\xi-\eta-\sigma)\\
&\quad\times \Widehat{f_{k_1}}(\eta)\Widehat{f_{\ll k}}(\sigma) H_k(-\xi) \, \diff \xi \diff \eta \diff \sigma \diff s\\
&=\sum_{m,k,k_1} I^\iota_{m,k,k_1}.
\end{aligned}
\end{equation*}
Integrating by parts with respect to $\eta$ leads to
\[
I^\iota_{m,k,k_1}=I^\iota_{m,k,k_1;1}+I^\iota_{m,k,k_1;2}+I^\iota_{m,k,k_1;3},
\]
where
\begin{equation*}
\begin{aligned}
I^\iota_{m,k,k_1;1}&=\sum_{m,k,k_1} \int q_m e^{-\mathrm{i} s \Phi} m_0 \partial_\eta \Widehat{f_{k_1}}(\xi-\eta-\sigma) \Widehat{f_{k_1}}(\eta)\\
&\quad\times \Widehat{f_{\ll k}}(\sigma) H_k(-\xi) \, \diff \xi \diff \eta \diff \sigma \diff s;\\
I^\iota_{m,k,k_1;2}&=\sum_{m,k,k_1} \int q_m e^{-\mathrm{i} s \Phi} m_0\Widehat{f_{k_1}}(\xi-\eta-\sigma) \partial_\eta \Widehat{f_{k_1}}(\eta)\\
&\quad\times \Widehat{f_{\ll k}}(\sigma) H_k(-\xi) \, \diff \xi \diff \eta \diff \sigma \diff s;\\
I^\iota_{m,k,k_1;3}&=\sum_{m,k,k_1} \int q_m e^{-\mathrm{i} s \Phi} \partial_\eta m_0 \Widehat{f_{k_1}}(\xi-\eta-\sigma) \Widehat{f_{k_1}}(\eta)\\
&\quad\times \Widehat{f_{\ll k}}(\sigma) H_k(-\xi) \, \diff \xi \diff \eta \diff \sigma \diff s,
\end{aligned}
\end{equation*}
with
\begin{equation*}
\begin{aligned}
m_0=\frac{-\Phi_\xi \xi^{1+2\iota} \chi_1 \theta_1 \zeta_1}{\Phi_\eta}.
\end{aligned}
\end{equation*}

Note that $k_1 \gg k$ and
\begin{equation*}
\begin{aligned}
\Phi_\eta=(\xi-\sigma)\tilde m,
\end{aligned}
\end{equation*}
where $\norm{\tilde m}_{S_{k,k_1,k_1,\ll k}} \lesssim 2^{-3k_1/2}$.
$I^\iota_{m,k,k_1;3}$ can be estimated easily. 
$I^\iota_{m,k,k_1;2}$ and $I^\iota_{m,k,k_1;1}$ can be handled in a similar way. So, it suffices to estimate $I^\iota_{m,k,k_1;1}$, which are divided into two cases depending on the size of $k$ as follows.\\

\noindent \underline{\underline{Case 1: $k \leq -m/2$.}}
By Theorem \ref{multi}, \eqref{weightasp2}, \eqref{estl2l}, \eqref{estH}, \eqref{estlih} and \eqref{estlilow}, it holds that
\begin{equation*}
\begin{aligned}
&\quad\sum_{m,k,k_1} \abs{I^\iota_{m,k,k_1;1}}\\
&\lesssim \sum_{m,k,k_1} 2^m 2^{2\iota k}2^{3k_1/2}\big\|\partial_\eta \Widehat{f_{k_1}}\big\|_2 \norm{H_k}_2 \norm{u_{k_1}}_\infty \norm{u_{\ll k}}_\infty\\
&\lesssim \sum_{m,k,k_1} \epsilon_1^4 2^m 2^{\iota k} 2^{3k_1/2} 2^{p_2m} (2^{m/6}|2^{p_2m}) 2^{-10k_1} 2^{-m/2} 2^{-m/2}\\
&\lesssim \epsilon_1^4 (2^{L/3}|2^{2p_2L}).
\end{aligned}
\end{equation*}\\

\noindent\underline{\underline{Case 2: $k>-m/2$.}}
By Lemma \ref{itmulti}, \eqref{weightasp2}, \eqref{estl2l}, \eqref{estH}, \eqref{decayhl}, \eqref{estlih} and \eqref{estli}, one can get
\begin{equation*}
\begin{aligned}
&\quad\sum_{m,k,k_1} \abs{I^\iota_{m,k,k_1;1}}\\
&\lesssim \sum_{m,k,k_1} 2^m 2^{2\iota k}2^{3k_1/2} (\big\|\partial_\eta \Widehat{f_{k_1}}\big\|_2 \norm{H_k}_2 I_{2^{m}}(u_{k_1},u_{\ll k})\\
&\quad+\langle 2^k 2^{m}\rangle^{-N} \big\|\partial_\eta \Widehat{f_{k_1}}\big\|_2 \norm{H_k}_2 \norm{u_{k_1}}_\infty\norm{u}_\infty)\\
&\lesssim
\sum_{m,k,k_1}( \epsilon_1^4 2^m 2^{\iota k}2^{3k_1/2} 2^{p_2m}(2^{m/6}|2^{p_2m}) 2^{-10k_1}2^{-m}\\
&\quad+ \epsilon_1^4 2^m 2^{\iota k}2^{3k_1/2} (2^k 2^m)^{-N} 2^{p_2m} 2^{m/6} 2^{-10k_1}2^{-m/2}2^{-m/3})\\
&\lesssim \epsilon_1^4 (2^{L/3}|2^{2p_2L}).
\end{aligned}
\end{equation*}

\subsubsection{Estimates for $F^\iota_{313}$}

Let
\[
\Phi_\xi=(\eta+\sigma)m_0,
\]
where $\norm{m_0}_{S_{k,k,\ll k,\ll k}} \lesssim 2^{\tilde k}2^{-3k_+/2}$. Since $\eta m_0$ and $\sigma m_0$ are symmetric, we only consider the terms involving $\eta m_0$.
Then we decompose $F^\iota_{313}$ by frequencies as
\begin{equation*}
\begin{aligned}
F^\iota_{313}&=\sum_{m,k} \int q_m e^{-\mathrm{i} s \Phi} (-\mathrm{i} s \eta m_0) \xi^{1+2\iota} \chi_1 \theta_3 \zeta_1 \Widehat{f_{k}}(\xi-\eta-\sigma) \\
&\quad\times\Widehat{f_{\ll k}}(\eta) \Widehat{f_{\ll k}}(\sigma) H_k(-\xi) \, \diff \xi \diff \eta \diff \sigma \diff s\\
&=\sum_{m,k} I^\iota_{m,k}.
\end{aligned}
\end{equation*}
Integrating by parts with $\eta$ leads to
\[
I^\iota_{m,k}=I^\iota_{m,k;1}+I^\iota_{m,k;2}+I^\iota_{m,k;3},
\]
where
\begin{gather*}
I^\iota_{m,k;1}= \int q_m e^{-\mathrm{i} s \Phi} m_1 \partial_\eta \Widehat{f_{k}}(\xi-\eta-\sigma) \Widehat{f_{\ll k}}(\eta) \Widehat{f_{\ll k}}(\sigma) H_k(-\xi) \, \diff \xi \diff \eta \diff \sigma \diff s;\\
I^\iota_{m,k;2}= \int q_m e^{-\mathrm{i} s \Phi} m_1 \Widehat{f_{k}}(\xi-\eta-\sigma) \partial_\eta \Widehat{f_{\ll k}}(\eta) \Widehat{f_{\ll k}}(\sigma) H_k(-\xi) \, \diff \xi \diff \eta \diff \sigma \diff s;\\
I^\iota_{m,k;3}= \int q_m e^{-\mathrm{i} s \Phi} \partial_\eta m_1 \Widehat{f_{k}}(\xi-\eta-\sigma) \Widehat{f_{\ll k}}(\eta) \Widehat{f_{\ll k}}(\sigma) H_k(-\xi) \, \diff \xi \diff \eta \diff \sigma \diff s
\end{gather*}
with
\[
m_1=\frac{-\eta m_0 \xi^{1+2\iota} \chi_1 \theta_3 \zeta_1}{\Phi_\eta}.
\]

In the following, we estimate $I^\iota_{m,k;1}, I^\iota_{m,k;2}, I^\iota_{m,k;3}$ term by term. \\

\noindent\boxed{\emph{Analysis of $I^\iota_{m,k;1}$}.} 
 
\noindent \underline{1. $k<100$.}
By Lemma \ref{itmulti}, \eqref{weightasp}, \eqref{weightasp2}, \eqref{estl2l}, \eqref{estH}, \eqref{decayinter} and \eqref{estli}, it holds that
\begin{equation*}
\begin{aligned}
&\quad\sum_{m,k} \abs{I^\iota_{m,k;1}}\\
&\lesssim
\sum_{m,k} 2^m 2^k 2^{(1+2\iota)k} 2^{-2k}( \big\|\partial_\eta \Widehat{f_k}\big\|_2 \norm{H_k}_2 I_{2^{m/3}}((u_{\ll k})_x,u_{\ll k})\\
&\quad+\langle 2^k 2^{m/3}\rangle^{-N} \big\|\partial_\eta \Widehat{f_k}\big\|_2 \norm{H_k}_2 \norm{(u_{\ll k})_x}_\infty \norm{u_{\ll k}}_\infty)
\\
&\lesssim \sum_{m,k}( \epsilon_1^2 2^{\iota k}(\norm{(xf)_k}_2+\epsilon_12^{-k/2}) \norm{\xi^\iota H_k}_2\\
&\quad+ \epsilon_1^3 2^{m}2^{\iota k} (2^k 2^{m/3})^{-N} 2^{m/6} \norm{\xi^\iota H_k}_2 2^k2^{-2m/3})\\
&\lesssim \epsilon_1^4 (2^{L/3}|2^{2p_2L}).
\end{aligned}
\end{equation*}

\noindent \underline{2. $k \geq 100$.}
We further decompose $I^\iota_{m,k;1}$ into
\begin{equation*}
\begin{aligned}
I^\iota_{m,k;1}&=\sum_{k_2} \int q_m e^{-\mathrm{i} s \Phi} m_1 \partial_\eta \Widehat{f_k}(\xi-\eta-\sigma) \Widehat{f_{k_2}}(\eta) \Widehat{f_{\ll k}}(\sigma) H_k(-\xi) \, \diff \xi \diff \eta \diff \sigma \diff s\\
&=\sum_{k_2}I^\iota_{m,k,k_2;1}
\end{aligned}
\end{equation*}
and will estimate it according to the size of $k_2$.\\

\noindent \underline{\underline{Case 1: $k_2 \geq -100$.}} It follows from Lemma \ref{itmulti}, \eqref{weightasp2}, \eqref{estl2}, \eqref{estH2}, \eqref{decayhl} and \eqref{estli} that
\begin{equation*}
\begin{aligned}
&\quad\sum_{m,k,k_2} \abs{I^\iota_{m,k,k_2;1}}\\
&\lesssim \sum_{m,k,k_2} 2^m 2^{k_2} 2^{-3k/2} 2^{(1+2\iota)k} 2^{k_2/2} (\big\|\partial_\eta \Widehat{f_k}\big\|_2 \norm{H_k}_2 I_{2^{m}}(u_{k_2},u_{\ll k})\\
&\quad+\langle 2^{k_2}2^m\rangle^{-N}\big\|\partial_\eta \Widehat{f_k}\big\|_2\norm{H_k}_2\norm{u}^2_\infty)\\
&\lesssim
\sum_{m,k,k_2}( \epsilon_1^4 2^{k_2} 2^{-3k/2} 2^{(1+2\iota)k} 2^{k_2/2} 2^{-k}2^{p_2m}2^{-k} 2^{p_2m} 2^{-15k_2}\\
&\quad+ \epsilon_1^4 2^m 2^{k_2} 2^{-3k/2} 2^{(1+2\iota)k}2^{k_2/2} (2^{k_2} 2^m)^{-N} 2^{-k}2^{p_2m}2^{-k}2^{p_2m} 2^{-2m/3})
\\
&\lesssim \epsilon_1^4 (2^{L/3}|2^{2p_2L}).
\end{aligned}
\end{equation*}

\noindent \underline{\underline{Case 2: $k_2<-100$.}}
By Lemma \ref{itmulti}, \eqref{weightasp2}, \eqref{estl2}, \eqref{estH2}, \eqref{decayinter} and \eqref{estli}, one can get
\begin{equation*}
\begin{aligned}
&\quad\sum_{m,k} \abs{I^\iota_{m,k,k_2;1}}\\
&\lesssim \sum_{m,k} 2^m 2^{-3k/2}2^{(1+2\iota)k}(\big\|\partial_\eta \Widehat{f_k}\big\|_2\norm{H_k} I_{2^{m/3}}((u_{<-100})_x,u_{\ll k})\\
&\quad+\langle 2^{m/3}\rangle^{-N} \big\|\partial_\eta \Widehat{f_k}\big\|_2\norm{H_k}_2 \norm{u}^2_\infty)\\
&\lesssim \sum_{m,k}( \epsilon_1^4 2^{-3k/2}2^{(1+2\iota)k} 2^{-k}2^{p_2m}2^{-k}2^{p_2m}\\
&\quad+\epsilon_1^4 2^{-Nm/3} 2^{-k}2^{p_2m}2^{-k}2^{p_2m} 2^{-2m/3})\\
&\lesssim \epsilon_1^4 (2^{L/3}|2^{2p_2L}).
\end{aligned}
\end{equation*}

\noindent\boxed{\emph{Analysis of $I^\iota_{m,k;2}$}.} 

\noindent \underline{1. $k<100$.} 
Integrating by parts with respect to $\sigma$ yields
\[
I^\iota_{m,k;2}=I^\iota_{m,k;21}+I^\iota_{m,k;22}+I^\iota_{m,k;23},
\]
where
\begin{gather*}
I^\iota_{m,k;21}= \int q_m e^{-\mathrm{i} s \Phi} m_2 \partial_\sigma \Widehat{f_{k}}(\xi-\eta-\sigma) \partial_\eta \Widehat{f_{\ll k}}(\eta) \Widehat{f_{\ll k}}(\sigma) H_k(-\xi) \, \diff \xi \diff \eta \diff \sigma \diff s;\\
I^\iota_{m,k;22}= \int q_m e^{-\mathrm{i} s \Phi} m_2 \Widehat{f_{k}}(\xi-\eta-\sigma) \partial_\eta \Widehat{f_{\ll k}}(\eta) \partial_\sigma \Widehat{f_{\ll k}}(\sigma) H_k(-\xi) \, \diff \xi \diff \eta \diff \sigma \diff s;\\
I^\iota_{m,k;23}= \int q_m e^{-\mathrm{i} s \Phi} \partial_\sigma m_2 \Widehat{f_{k}}(\xi-\eta-\sigma) \partial_\eta \Widehat{f_{\ll k}}(\eta) \Widehat{f_{\ll k}}(\sigma) H_k(-\xi) \, \diff \xi \diff \eta \diff \sigma \diff s
\end{gather*}
with
\[
m_2=\frac{m_1}{\mathrm{i} s \Phi_\sigma}.
\]

For $I^\iota_{m,k;21}$, due to
\[
\norm{m_2/\eta}_{S_{k,k,\ll k,\ll k}} \lesssim 2^{-m}2^{-3k} 2^{(1+2\iota)k},
\]
it follows from Theorem \ref{multi}, \eqref{weightasp}, \eqref{weightasp2}, \eqref{estl2l}, \eqref{estli} and \eqref{estH} that
\begin{equation*}
\begin{aligned}
&\quad\sum_{m,k} \abs{I^\iota_{m,k;21}}\\
&\lesssim \sum_{m,k} 2^m 2^{-m} 2^{-3k} 2^{(1+2\iota)k} \big\|\partial_\sigma \Widehat{f_k}\big\|_2 \big\|\eta \partial_\eta \Widehat{f_{\ll k}}\big\|_2 \norm{u_{\ll k}}_\infty \norm{E\mathcal{F}^{-1}(H_k)}_\infty\\
&\lesssim \sum_{m,k} \epsilon_1^4 2^{-3k} 2^{(1+\iota)k} 2^{m/6} (2^k2^{m/6}|2^{p_2m})2^{-m/3}2^{k/2} (2^{m/6}|2^{p_2m})\\
&\lesssim \epsilon_1^4(2^{L/3}|2^{2p_2L}).
\end{aligned}
\end{equation*}

$I^\iota_{m,k;22}$ can be handled similarly as for $I^\iota_{m,k;21}$, and $I^\iota_{m,k;23}$ is easy to deal with.\\

\noindent \underline{2. $k \geq 100$.} It follows from Lemma \ref{itmulti}, \eqref{weightasp2}, \eqref{energyasp}, \eqref{estH2}, \eqref{decayhl}, \eqref{estlih} and \eqref{estli} that
\begin{equation*}
\begin{aligned}
&\quad\sum_{m,k} \abs{I^\iota_{m,k;2}}\\
&\lesssim \sum_{m,k} 2^m 2^{-3k/2} 2^{(1+2\iota)k} 2^{10k} (\big\|\eta \partial_\eta \Widehat{f_{\ll k}}\big\|_2 \norm{H_k}_2 I_{2^m}(u_k,u_{\ll k})\\
&\quad+\langle 2^k 2^m\rangle^{-3} \big\|\eta \partial_\eta \Widehat{f_{\ll k}}\big\|_2 \norm{H_k}_2 \norm{u_k}_\infty \norm{u}_\infty)\\
&\lesssim
\sum_{m,k}( \epsilon_1^42^{-3k/2} 2^{(1+\iota)k} 2^{10k}2^{p_2m}2^{-k}2^{p_2m}2^{-15k}\\
&\quad+ \epsilon_1^4 2^m 2^{-3k/2} 2^{(1+2\iota)k} 2^{10k} (2^k2^m)^{-3} 2^{p_2m} 2^{-k}2^{p_2m} 2^{-15k} 2^{-m/2} 2^{-m/3})
\\
&\lesssim \epsilon_1^4(2^{L/3}|2^{2p_2L}).
\end{aligned}
\end{equation*}

\noindent\boxed{\emph{Analysis of $I^\iota_{m,k;3}$}.} 

\noindent \underline{1. $k<100$.} 
Lemma \ref{itmulti}, \eqref{estl2l}, \eqref{estH}, \eqref{decayinter} and \eqref{estli} imply
\begin{equation*}
\begin{aligned}
\sum_{m,k} \abs{I^\iota_{m,k;3}}
&\lesssim \sum_{m,k} 2^m2^{2\iota k} (\norm{f_{\ll k}}_2\norm{H_k}_2 I_{2^{m/3}}(u_k,u_{\ll k})\\
&\quad+\langle 2^k 2^{m/3}\rangle^{-N} \norm{f_{\ll k}}_2 \norm{H_k}_2 \norm{u}^2_\infty) \\
&\lesssim \sum_{m,k}( \epsilon_1^4 2^{\iota k}2^{k/2}(2^{m/6}|2^{p_2m})2^{-k}\\
&\quad+ \epsilon_1^4 2^m 2^{2\iota k}(2^k 2^{m/3})^{-N} 2^{k/2} 2^{m/6}2^{-2m/3})\\
&\lesssim \epsilon_1^4 (2^{L/3}|2^{2p_2L}).
\end{aligned}
\end{equation*}

\noindent \underline{2. $k \geq 100$.}
In this case, $I^\iota_{m,k;3}$ can be estimated similarly as for $I^\iota_{m,k;2}$.

\subsubsection{Estimates for $F^\iota_{314}$}
We decompose $F^\iota_{314}$ into
\begin{equation*}
\begin{aligned}
F^\iota_{314}&=\sum_{m,k,k_1} \int q_m e^{-\mathrm{i} s \Phi} (-\mathrm{i} s \Phi_\xi) \xi^{1+2\iota} \chi_1 \theta_4 \zeta_1 \Widehat{f_{k_1}}(\xi-\eta-\sigma) \Widehat{f_{k_1}}(\eta)\\
 &\quad\times \Widehat{f_{k_1}}(\sigma) H_k(-\xi) \, \diff \xi \diff \eta \diff \sigma \diff s\\
&=\sum_{m,k,k_1} I^\iota_{m,k,k_1}.
\end{aligned}
\end{equation*}
Note that $k_1 \gg k$. Choose cut-off functions $\rho_1$ and $\rho_2$ supporting in $\abs \eta \geq 0.9 \abs \sigma$ and $\abs \sigma \geq 0.9 \abs \eta$ respectively and add them into $I^\iota_{m,k,k_1}$. We only consider the term concerning $\rho_1$ since the one concerning $\rho_2$ is similar, which is still denoted as $I^\iota_{m,k,k_1}$ for convenience. 

Integrating by parts with respect to $\eta$ gives
\[
I^\iota_{m,k,k_1}=I^\iota_{m,k,k_1;1}+I^\iota_{m,k,k_1;2}+I^\iota_{m,k,k_1;3},
\]
where
\begin{gather*}
I^\iota_{m,k,k_1;1}=\int q_m e^{-\mathrm{i} s \Phi} m_0 \partial_\eta \Widehat{f_{k_1}}(\xi-\eta-\sigma) \Widehat{f_{k_1}}(\eta) \Widehat{f_{k_1}}(\sigma)H_k(-\xi) \, \diff \xi \diff \eta \diff \sigma \diff s;\\
I^\iota_{m,k,k_1;2}=\int q_m e^{-\mathrm{i} s \Phi} m_0 \Widehat{f_{k_1}}(\xi-\eta-\sigma) \partial_\eta \Widehat{f_{k_1}}(\eta) \Widehat{f_{k_1}}(\sigma) H_k(-\xi) \diff \xi\diff \eta \diff \sigma \diff s;\\
I^\iota_{m,k,k_1;3}=\int q_m e^{-\mathrm{i} s \Phi} \partial_\eta m_0 \Widehat{f_{k_1}}(\xi-\eta-\sigma) \Widehat{f_{k_1}}(\eta) \Widehat{f_{k_1}}(\sigma) H_k(-\xi) \diff \xi\diff \eta \diff \sigma \diff s
\end{gather*}
with
\[
m_0=\frac{-\Phi_\xi \xi^{1+2\iota} \chi_1 \theta_4 \zeta_1 \rho_1}{\Phi_\eta}.
\]

When $k_1 \leq 100$,
\[
\norm{m}_{S_{k,k_1,k_1,k_1}} \lesssim 2^{(1+2\iota)k}.
\]
It follows from Theorem \ref{multi} and \eqref{estlil} that
\[
\abs{I^\iota_{m,k,k_1}} \lesssim \epsilon_1^2 2^m 2^{(1+2\iota)k} \big\|\partial_\eta \Widehat{f_{k_1}}\big\|_2 \norm{H_k}_2 2^{-m} 2^{-k_1}.
\]
By Young's inequality, $I^\iota_{m,k,k_1;1}$ could be estimated. $I^\iota_{m,k,k_1;2}$ can be handled similarly and $I^\iota_{m,k,k_1;3}$ is easy to be estimated. 

When $k_1>100$,
\[
\norm{m}_{S_{k,k_1,k_1,k_1}} \lesssim 2^{(1+2\iota)k}2^{k_1/2}.
\]
It follows from Theorem \ref{multi} and \eqref{estlih} that
\[
\abs{I^\iota_{m,k,k_1;1}} \lesssim \epsilon_1^2 2^m 2^{(1+2\iota)k}2^{k_1/2} \big\|\partial_\eta \Widehat{f_{k_1}}\big\|_2 \norm{H_k}_2 2^{-m} 2^{-10k_1}. 
\]
Hence $I^\iota_{m,k,k_1;1}$ can be estimated. $I^\iota_{m,k,k_1;2}$ can be handled similarly and $I^\iota_{m,k,k_1;3}$ is easy to deal with.

\subsection{Estimates for $F^\iota_{32}$ and $F^\iota_{33}$}
Since $F^\iota_{32}$ and $F^\iota_{33}$ can be handled in a similar way, we only consider $F^\iota_{32}$ here. Furthermore, it suffices to estimate the following generic terms in $F^\iota_{32}$:
\begin{equation*}
\begin{aligned}
&\int q_m e^{-\mathrm{i} s \Phi}\varphi_{(-m/3,p_1m)}(\xi-\eta-\sigma) \frac{\Phi_\xi}{\Phi} \xi^{1+2\iota} \chi_1\zeta_2 \Widehat f(\xi-\eta-\sigma)\\
&\quad\times \Widehat f(\eta) \Widehat f(\sigma) H(-\xi) \, \diff \xi \diff \eta \diff \sigma \diff s
\end{aligned}
\end{equation*}
and
\begin{equation*}
\begin{aligned}
&\int q_m e^{-\mathrm{i} s \Phi}\varphi_{(-m/3,p_1m)}(\xi-\eta-\sigma)  \frac{s\Phi_\xi}{\Phi} \xi^{1+2\iota}  \chi_1\zeta_2 \partial_s \Widehat f(\xi-\eta-\sigma)\\
&\quad\times \Widehat f(\eta) \Widehat f(\sigma) H(-\xi) \, \diff \xi \diff \eta \diff \sigma \diff s,
\end{aligned}
\end{equation*}
which are divided by frequencies respectively as
\begin{equation}\label{f321}
\begin{aligned}
&\sum_{m,k} \int q_m e^{-\mathrm{i} s \Phi}\varphi_{(-m/3,p_1m)}(\xi-\eta-\sigma) \frac{\Phi_\xi}{\Phi} \xi^{1+2\iota} \chi_1\zeta_2 \Widehat{f_k}(\xi-\eta-\sigma)\\
&\quad\times \Widehat{f_k}(\eta) \Widehat{f_{\ll k}}(\sigma) H_k(-\xi) \, \diff \xi \diff \eta \diff \sigma \diff s
\end{aligned}
\end{equation}
and
\begin{equation}\label{f322}
\begin{aligned}
&\sum_{m,k} \int q_m e^{-\mathrm{i} s \Phi}\varphi_{(-m/3,p_1m)}(\xi-\eta-\sigma) \frac{s\Phi_\xi}{\Phi} \xi^{1+2\iota} \chi_1\zeta_2 \partial_s \Widehat{f_k}(\xi-\eta-\sigma) \\
&\quad\times\Widehat{f_k}(\eta) \Widehat{f_{\ll k}}(\sigma) H_k(-\xi) \, \diff \xi \diff \eta \diff \sigma \diff s. 
\end{aligned}
\end{equation}

The key here is to note from Corollary \ref{reso32} that
\[
\norm{\Phi_\xi/\Phi}_{S_{k,k,k,\ll k}} \lesssim 2^{-k}2^{5k_+}.
\]
Then for \eqref{f321}, it follows from Theorem \ref{multi}, \eqref{estlil}, \eqref{estlih}, \eqref{estl2l} and \eqref{estH} that
\begin{align*}
\abs{\eqref{f321}} &\lesssim \sum_{m,k} \epsilon_1^4 2^m 2^{-k}2^{5k_+} 2^{(1+\iota)k}2^{-k}2^{-10k_+}2^{-m}2^{k/2} (2^{m/6}|2^{p_2m})\\
&\lesssim \epsilon_1^4 (2^{L/3}|2^{2p_2L}).
\end{align*}
While for \eqref{f322}, it follows from Lemma \ref{itmulti}, \eqref{psl}, \eqref{estH}, \eqref{decayinter}, \eqref{decayhl}, \eqref{estlih} and \eqref{estli} that
\begin{equation*}
\begin{aligned}
\abs{\eqref{f322}} &\lesssim \sum_{m,k} 2^m 2^m 2^{-k} 2^{10k_+} 2^{(1+2\iota)k} (\big\|\partial_s \Widehat{f_k}\big\|_2\norm{H_k}_2 I_{2^{m/3}}(u_k,u_{\ll k})\\
&\quad+\langle 2^k2^{m/3}\rangle^{-3}\big\|\partial_s \Widehat{f_k}\big\|_2\norm{H_k}_2\norm{u_k}_\infty\norm{u}_\infty)\\
&\lesssim
\sum_{m,k}( \epsilon_1^4 2^{2m} 2^{-k} 2^{10k_+} 2^{(1+\iota)k} 2^{k/2}2^{-m} (2^{m/6}|2^{p_2m}) 2^{-k}2^{-15k_+}2^{-m}\\
&\quad+ \epsilon_1^4 2^{2m} 2^{-k} 2^{10k_+} 2^{(1+\iota)k}(2^{k}2^{m/3})^{-3} 2^{k/2}2^{-m} (2^{m/6}|2^{p_2m}) )\\
&\quad\times 2^{-15k_+}2^{-m/3}2^{-m/3}
\\
&\lesssim \epsilon_1^4(2^{L/3}|2^{2p_2L}).
\end{aligned}
\end{equation*}

\subsection{Estimates for $F^\iota_{34}$}
Consider the following generic terms in $F^\iota_{34}$:
\begin{equation*}
\begin{aligned}
& \int q_m e^{-\mathrm{i} s \Phi} \varphi_{(-m/3,p_1m)}(\xi-\eta-\sigma)\xi^{2\iota} m_1 \zeta_4 \Widehat{f^{\iota_1}}(\xi-\eta-\sigma) \\
&\quad\times \Widehat{f^{\iota_2}}(\eta) \Widehat{f^{\iota_3}}(\sigma) H(-\xi) \, \diff \xi \diff \eta \diff \sigma \diff s,
\end{aligned}
\end{equation*}
\begin{equation*}
\begin{aligned}
& \int q_m e^{-\mathrm{i} s \Phi} \varphi_{(-m/3,p_1m)}(\xi-\eta-\sigma)\xi^{2\iota}s m_2 \zeta_4 \partial_s \Widehat{f^{\iota_1}}(\xi-\eta-\sigma)\\
&\quad\times \Widehat{f^{\iota_2}}(\eta) \Widehat{f^{\iota_3}}(\sigma) H(-\xi) \, \diff \xi \diff \eta \diff \sigma \diff s
\end{aligned}
\end{equation*}
and
\begin{equation*}
\begin{aligned}
& \int q_m e^{-\mathrm{i} s \Phi} \varphi_{(-m/3,p_1m)}(\xi-\eta-\sigma)\xi^{2\iota}(-\mathrm{i} s m_3 \Phi_\eta) \zeta_4 \Widehat{f^{\iota_1}}(\xi-\eta-\sigma)\\
&\quad\times \Widehat{f^{\iota_2}}(\eta) \Widehat{f^{\iota_3}}(\sigma) H(-\xi) \, \diff \xi \diff \eta \diff \sigma \diff s,
\end{aligned}
\end{equation*}
where $\norm{m_1}_{S_{k,k,k,k}} \lesssim 1$, $\norm{m_2}_{S_{k,k,k,k}} \lesssim 1$, and $\norm{m_3}_{S_{k,k,k,k}} \lesssim 2^k$. Decompose them according to different frequencies and integrate by parts with respect to $\eta$ for the third expression above. Since the frequencies are comparable, these terms could be estimated by Theorem \ref{multi}. We omit the details here.

\section*{Acknowledgments}
This is part of the Ph.D thesis of the first author written under the supervision of the third author at the Institute of Mathematical Science at the Chinese University of Hong Kong and is supported by Zheng Ge Ru Foundation, Hong Kong RGC
Earmarked Research Grants CUHK
-14301421, CUHK-14301023, CUHK-14300819 and CUHK-
14302819, and the key projects of NSFC Grants No. 12131010 and No. 11931013. The second author's research is also supported by the Grant No.
830018 from China.

%
%

\end{document}